\documentclass[journal]{new-aiaa}
\usepackage[utf8]{inputenc}
\usepackage{textcomp}

\usepackage{graphicx}
\usepackage{amsmath}
\usepackage[version=4]{mhchem}
\usepackage{siunitx}
\usepackage{longtable,tabularx}

\usepackage{bm}
\usepackage{amsmath}

\usepackage{color,soul}
\usepackage[color=gray!30]{todonotes}
\usepackage{colortbl}
\definecolor{myblue}{rgb}{0, 0.23, 0.64}
\usepackage{hyperref}
\hypersetup{
    colorlinks=true,
    linkcolor=myblue,
    filecolor=myblue,
    citecolor=myblue,
    urlcolor=myblue,
}
\usepackage{booktabs}
\usepackage{graphicx}
\usepackage[version=4]{mhchem}
\usepackage{siunitx}
\usepackage{longtable,tabularx}
\usepackage{gensymb}
\usepackage{mathtools}
\usepackage{float}
\usepackage[ruled, vlined, linesnumbered]{algorithm2e}
\SetKwComment{Comment}{$\triangleright$\ }{}
\usepackage{subcaption}
\usepackage{enumitem}
\usepackage{multirow}
\usepackage{pdflscape}
\usepackage{lscape}
\usepackage{mdframed}
\usepackage{adjustbox}
\usepackage{dsfont}
\usepackage{xspace}
\usepackage[flushleft]{threeparttable}

\newcommand{\rhobar}{\bm{\rho}}

\newcommand{\lone}{$L_{1}$\xspace}
\newcommand{\ltwo}{$L_{2}$\xspace}
\let\svthefootnote\thefootnote
\newcommand\freefootnote[1]{%
  \let\thefootnote\relax%
  \footnotetext{#1}%
  \let\thefootnote\svthefootnote%
}

\title{Embedded State Estimation for Optimization of Cislunar Space Domain Awareness Constellation Design}

\author{Thomas H. Clareson\footnote{Master's Student, Department of Mechanical, Materials and Aerospace Engineering, West Virginia University, Student Member AIAA.}, Matthew C. Fox\footnote{Master's Student, Department of Mechanical, Materials and Aerospace Engineering, West Virginia University, Student Member AIAA.}, Dominic K. Amato\footnote{Master's Student,  Department of Mechanical, Materials and Aerospace Engineering, West Virginia University, Student Member AIAA.}, and Hang Woon Lee\footnote{Assistant Professor, Department of Mechanical, Materials and Aerospace Engineering; hangwoon.lee@mail.wvu.edu. Member AIAA (Corresponding Author).}}
\affil{West Virginia University, Morgantown, WV 26506}

\begin{document}

\freefootnote{This paper is a substantially revised version of the Paper AAS 23-189, presented at the AAS/AIAA Astrodynamics Specialist Conference, Big Sky, MT, August 13-17, 2023. It offers new results and a better description of the materials.}

\maketitle

\begin{abstract}
The traffic in cislunar space is expected to increase over the coming years, leading to a higher likelihood of conjunction events among active satellites, orbital debris, and non-cooperative satellites. This increase necessitates enhanced space domain awareness (SDA) capabilities that include state estimation for targets of interest. Both Earth surface-based and space-based observation platforms in geosynchronous orbit or below face challenges such as range, exclusion, and occlusion that hinder observation. Motivated by the need to place space-based observers in the cislunar space regime to overcome these challenges, this paper proposes a cislunar SDA constellation design and analysis framework that integrates state estimation into an optimization problem for determining the placement of observers for optimal state estimation performance on a set of targets. The proposed multi-observer placement optimization problem samples from a range of possible target orbits. Upon convergence, the optimized constellation is validated against a broader set of targets to assess its effectiveness. Two comparative analyses are presented to evaluate the effects of changes in the sensor tasking procedure and sensor fidelity on the optimized constellation, comparing these to a single observer baseline case. The results demonstrate that the optimized constellations can provide accurate state estimation for various orbit families.
\end{abstract}

\section*{Nomenclature}
{\renewcommand\arraystretch{1.0}
\noindent\begin{longtable*}{@{}l @{\quad=\quad} l@{}}
$\bm{A}$ & state function Jacobian matrix\\
            $\bm{f}$ & system dynamics\\
            $G$ & nondimensionalized gravitational constant\\
            $\bm{h}$ & measurement function\\
            $\bm{H}$ & measurement Jacobian matrix\\
            $\mathcal{J}$ & set of orbital slots \\
            $\bm{K}$ & Kalman gain\\
            $L_{n}$ & Lagrange point; $n$ represents a particular Lagrange point\\
            $\mathcal{L}$ & loss function \\
            $\bm{m}$ & measurement function output vector\\
            $N$ & number of observers \\
            $\mathcal{N}$ & Gaussian distribution \\
            $\bm{P}$ & covariance matrix\\
            $\bm{Q}$ & process noise covariance matrix\\
            $\bm{R}$ & measurement noise covariance\\
            $\mathcal{U}$ & optimization target set \\
            $\bm{v}$ & measurement function noise vector\\
            $\bm{w}$ & process noise vector\\
            $\bm{x}$ & state vector \\
            $\hat{\bm{x}}$ & estimated position and velocity vector of target\\
            $\bm{X}$ & decision variable vector\\
            $\alpha$ & azimuth, radians\\
            $\bm{\gamma}$ & relative position vector, pointing from observer to reference point\\
            $\varepsilon$ & elevation, radians\\
            $\theta$ & angle between the target, observer, and gravitational body, radians\\
            $\mu$ & mass ratio\\
            $\bm{\rho}$ & relative position vector, pointing from observer to target\\
            $\phi$ & exclusion angle, radians\\
            $\omega$ & angle between the tangent line of gravitational body's surface, the observer, and the gravitational body, radians\\
\end{longtable*}}

\section{Introduction}
Within the near future, cislunar space---the region of space between Earth and the Moon---is expected to witness a significant increase in traffic. Commercial and governmental entities, both foreign and domestic, have proposed multiple crewed and uncrewed missions. Robotic lander and rover missions, as well as orbiter missions have been launched or planned to launch by various firms, such as Intuitive Machines, Astrobotic, iSpace, and Advanced Space \cite{jones2024back,gardner2023capstone}. NASA has organized a program with international and commercial support for the human return to the lunar surface known as Artemis, with one subset of the program being Gateway, a lunar space station \cite{smith2020artemis,crusan2019nasa,lehnhardt2024}. This increase in traffic requires mission operators to have better knowledge of resident space objects (RSOs) within the regime, through detection, tracking, and characterization of targets, due to higher chances of conjunction events. The practice of gaining this knowledge is referred to as \textit{space domain awareness} (SDA). An integral part of SDA is designing observation systems (henceforth referred to as \textit{observers}) that can provide this information, which includes constructing system hardware, creating procedures and software, and determining the optimal placement of the observers that will be tasked with gaining this information, the last of which is the focus of this paper.

Determining the placement of observers for cislunar SDA systems is more challenging than for SDA systems designed to track RSOs in the regime spanning between low Earth orbit (LEO) and geostationary orbit (GEO). The literature proposes SDA systems for the LEO-GEO regime that rely on the use of observers based on Earth's surface \cite{boer2017tarot,haimerl2015space}, as well as space-based observers occupying orbits within this regime \cite{du2019tenative,tommila2024mission}. Utilizing ground-based or space-based observers in the LEO-GEO regime for cislunar SDA poses significant challenges due to the unique difficulties of the cislunar environment, as outlined in Ref.~\cite{holzinger2021primer}. The large distance between the Earth and the Moon causes observers on Earth's surface and in orbit within the LEO-GEO region to struggle to detect and maintain tracking of targets in a volume of space extending far beyond GEO. Moreover, observers have difficulties taking measurements of targets due to exclusion and occlusion from the Earth, Moon, and Sun. To address these issues, observers can be positioned closer to the RSOs they aim to monitor. Two suggested locations include observers on the lunar surface and space-based observers occupying orbits within cislunar space.

A solution to the aforementioned cislunar SDA issues is placing observers on the lunar surface. Zimmer et al. \cite{zimmer2021cislunar} discuss the physical design of observation platforms and several design constraints. Additionally, Jarrett-Izzi et al. \cite{jarrett2024investigation} perform single target tracking to demonstrate the relationship between the target's orbit and the observer's placement on the lunar surface, and multi-target tracking using a Lyapunov-based algorithm, with results showing the capability of one observer to take measurements of multiple targets. Spesard et al. \cite{spesard2024prelim} approach the placement of lunar surface-based observers as an optimization problem, aiming to maximize the coverage of three volumes of space: the general cislunar region, the low lunar orbit region, and the Moon exclusion region in a series of case studies, utilizing a Genetic Algorithm (GA). Koblick et al. \cite{koblick2022cislunar} investigate a mixed-architecture system consisting of observers on the lunar surface and in cislunar periodic orbits. They examine a few different architectures with varying placements and parameters of lunar surface-based observers, along with one cislunar space-based observer. Results suggest an improvement in position error when the space-based observer is included in the architecture.

Lunar surface-based observers can face the challenges of having a limited view of cislunar space. To circumvent this limitation, there has been research on utilizing observers placed in cislunar space to provide SDA capabilities. In Block et al. \cite{block2022cislunar}, numerical analyses are performed to compare varying numbers of low- and high-fidelity sensors in the state estimation of a target in an \lone Halo orbit. Wilmer et al. \cite{wilmer2022cislunar} investigate the ability of various observers to view a set of targets, where the apparent magnitude of targets from the observer's location is a main consideration. In providing space traffic management (STM), Cunio et al. \cite{cunio2021utilization} present a methodology for dividing cislunar space into smaller, more manageable divisions, also presenting the orbits that provide the best surveillance of said regions. Gupta et al. \cite{gupta2022long} investigates resonant orbits, and their abilities to provide surveillance for regions of cislunar space. Pole-sitter orbits are offered as a possible solution for cislunar surveillance capabilities in Ewart et al. \cite{ewart2022pole}, where architecture and mission design are outlined for a possible mission. Vendl and Holzinger \cite{vendl2021cislunar} investigate cislunar target visibility, referring to the observers' ability to see these points in space, and/or their apparent magnitude. The tradeoff between visual magnitude and orbit stability provides an initial inquiry into utilizing optimization for the placement of cislunar space-based observers.

Optimization methods have been utilized in cislunar SDA constellation design to find the best placements of observers. Fahrner et al. \cite{fahrner2022capacity} approach the observer placement problem in a two-step process, combining observability and capacity analysis. In this work, observability analysis, finding the orbit within each family (and later constellation architecture) that provides the best visibility of a discretized section of cislunar space is performed, as well as capacity analysis: sizing and scheduling the best-performing constellations. In Segal et al. \cite{segal2023optimal}, the formulation of the optimization problem considers three metrics; maximizing the percent coverage of a set of orbits, maximizing the total number of observation opportunities per observer per time step, and minimizing the total number of observers in the constellation. The decision variables are the placements of these observers within different orbits. Visonneau et al. \cite{visonneau2023opt} formulate an optimization problem placing orbits to maximize the visibility of the volume of space from the GEO belt to \ltwo. The ability of the observers to observe a target is based on the target's apparent magnitude from the observers' location.

The aforementioned research in cislunar SDA has provided meaningful results in a wide variety of areas. In Refs.~\cite{cunio2021utilization, ewart2022pole}, a variety of orbits are proposed to provide STM/SDA capabilities, providing analysis of predicted traffic, and mission design. However, they do not provide analysis of their performance tracking targets. This information could justify the use of these orbits within a cislunar SDA constellation. References~\cite{block2022cislunar, wilmer2022cislunar, koblick2022cislunar, gupta2022long} contribute visibility and/or state estimation analysis of a variety of orbits. Performing analysis of these metrics provides insight into the performance of these orbits. This analysis, however, does not leverage optimization methods, which could reveal combinations of orbits that may be overlooked by a human designer. Reference~\cite{vendl2021cislunar}  provides visibility analysis, and additionally touches on optimization, through investigation of the tradeoff between the observability a certain orbit
can provide, and its stability. This examination, however, does not present an optimal constellation. In Ref.~\cite{fahrner2022capacity}, sensor tasking is intertwined into their optimization problem, through the addition of the scheduling problem for multiple targets. There is, however, importance in approaching the sensor tasking problem for a single target. References~\cite{fahrner2022capacity, segal2023optimal, visonneau2023opt} also formulate optimization problems, and computational experiments are performed to find the optimal constellations for cislunar SDA. In these works, the goal is to provide the highest visibility of either a select group of orbits or a volume of space. To provide the best SDA capabilities in cislunar space, the optimal placement of observers within a constellation should be found. The aforementioned works have utilized visibility analysis methods to provide this placement. However, other works have shown the importance of testing the state estimation capabilities of SDA systems, and how this testing can give more information about targets in the regime. Since the goal of SDA is to provide not only the detection of targets, but additionally tracking and characterization, we recognize the importance of intertwining optimization and state estimation into one design framework, embedding state estimation directly into the optimization of observer placement.

The contributions of this paper are as follows: first, this work proposes a novel optimization framework, the \textit{Cislunar SDA Constellation Design and Analysis Framework}, for designing space-based observer constellations for cislunar SDA purposes and conducting performance analysis. Recognizing that one of the primary objectives of SDA activities is to estimate the states of targets, we embed state estimation within the optimization framework. To design a constellation aimed at a set of unknown objects, we formulate an optimization problem in a deterministic manner, where we aim to find the constellation that can provide the best state estimation performance against a sampling of orbits from representative periodic orbit families. To validate the performance of the optimized constellation's state estimation, we conduct a post-optimization analysis against a wider set of targets not part of the optimization; the results indicate that the proposed deterministic design approach can lead to effective constellation designs. We provide two comparative analyses: (1) an experiment presenting the benefits of using a constellation, with a comparison to a single observer [in a Near-Rectilinear Halo Orbit (NRHO) similar to that of the Gateway program], and comparing three unique single-target sensor tasking procedures used with the optimization framework, varying the cadence of measurement and number of tasked sensors, showing how the change in sensor tasking changes the optimal constellation, and (2) an experiment comparing the optimal constellations when the fidelity of the observers is varied.

This paper is organized as follows. First, the framework used to optimize observer placement is discussed in Sec.~\ref{sec:prelims}; the subsections describe the dynamic model, state estimation filter, measurement model, optimization problem formulation, and post-optimization analysis process. Then, Sec.~\ref{sec:comp_ex} provides an overview of the experimental setup and reports the results of the comparative analyses and optimal constellation design. Finally, Sec.~\ref{sec:conc} concludes with a discussion of the new findings and avenues for future work in this area.

This paper expands on our preliminary work \cite{clareson2023opt} by offering new analysis through the implementation of varying sensor tasking procedures, an expanded set of target and observer orbit families, and a more in-depth description of the component materials of the framework.

\section{Cislunar SDA Constellation Design and Analysis Framework}
\label{sec:prelims}
This section describes the Circular Restricted three-body Problem (CR3BP) model, the basic algorithm of the Extended Kalman Filter (EKF), and the measurement model. Next, an outline of the optimization problem formulation is presented.

\subsection{Circular Restricted Three-Body Problem}
\begin{figure}[htb]
	\centering
        \includegraphics[width=.3\linewidth]{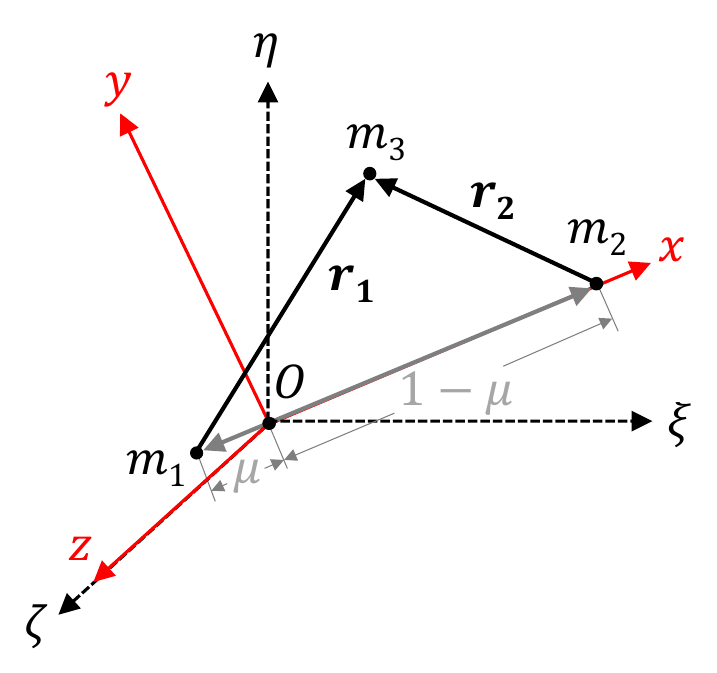}
	\caption{The CR3BP synodic frame ($x,y,z$) and the inertial frame ($\xi,\eta,\zeta$).}
	\label{fig:cr3bp}
\end{figure}

In this paper, we adopt the CR3BP as the governing dynamics for both observers and target RSOs, owing to its suitability for preliminary design and analysis purposes. The CR3BP defines the motion of a third body (\textit{e.g.}, observers and target RSOs) under the influence of two massive primary bodies, where the mass of the two primaries is large enough that the third body's mass can be considered negligible. An inertial frame is defined in the $(\xi,\eta,\zeta)$ basis, while a rotating reference frame, commonly called the synodic reference frame, is defined in the $(x,y,z)$ basis, such that the first and second primaries are stationary along the $x$-axis, rotating about the $z$-axis at a rate equal to the mean motion of the system. These frames are shown along with the primaries in Fig.~\ref{fig:cr3bp}. The locations of the first and second primaries in the synodic frame are $\mu$ and $1-\mu$ respectively, where $\mu$ is the nondimensionalized mass ratio, whose formula is $\mu=m_2/(m_1+m_2)$, with the origin being the common barycenter of these two primaries.

The equations of motion for the third body in the CR3BP model in the synodic frame are shown in Eq.~\eqref{EOM}~\cite{howell1984almost}:
\begin{subequations}
\label{EOM}
\begin{align}
    \Ddot{x}-2\dot{y} &= \frac{\partial U}{\partial x} \label{EOM_x}\\
    \Ddot{y} + 2\dot{x} &= \frac{\partial U}{\partial y} \label{EOM_y}\\
    \Ddot{z} &= \frac{\partial U}{\partial z} \label{EOM_z}
\end{align}
\end{subequations}
where the position of the third body in the synodic frame is given as $x$, $y$, and $z$, and $U$ is the pseudo-potential function defined as:
\begin{equation}
    U = \frac{1}{2}\left(x^2+y^2\right)+\frac{1-\mu}{r_1}+\frac{\mu}{r_2}
    \label{pseudo-pot}
\end{equation}
In Eq.~\eqref{pseudo-pot}, $r_1$ represents the distance from the first primary to the third body, while $r_2$ represents the distance from the second primary to the third body, respectively. They are defined as follows:
\begin{subequations}
    \begin{align}
        r_1 &= \sqrt{(x+\mu)^2+y^2+z^2} \\
        r_2 &= \sqrt{\left(x-(1-\mu)\right)^2+y^2+z^2}
    \end{align}
    \label{r1_r2}
\end{subequations}
Length, mass, and time are nondimensionalized using the characteristic terms shown in Eqs.~\eqref{nondim_l}--\eqref{nondim_t}, respectively:
\begin{subequations}
  \begin{align}
    l^\ast &= d_1+d_2 \label{nondim_l}\\
    m^\ast &= m_1+m_2 \label{nondim_m}\\
    t^\ast &= \sqrt{\frac{{(l^\ast)}^3}{Gm^\ast}} \label{nondim_t}
\end{align}  
\end{subequations}
In Eq.~\eqref{nondim_l}, the distance from the center of the first primary to the barycenter is given as $d_1$, and the distance from the center of the second primary to the barycenter is given as $d_2$. Equation~\eqref{nondim_t} uses a nondimensionalized form of the gravitational constant $G$, which allows the nondimensionalized mean motion to be unity. Table~\ref{tab:constants} shows the CR3BP constants and canonical units for the Earth-Moon system \cite{jpl3bpdata}. Additionally, within CR3BP, there are five equilibrium points, known as Lagrange points, labeled as \lone--$L_5$.

\begin{table}[htb]
\renewcommand{\arraystretch}{1}
    \caption{CR3BP constants and canonical values for the Earth-Moon system \cite{jpl3bpdata}.}
    \label{tab:constants}
      \centering
      \begin{tabular}{l r r}
        \hline \hline
        Constant & Symbol &  Value\\
        \hline
            Mass ratio & $\mu$ & 0.0121506\\ 
            Distance unit (DU) & $l^\ast$ & \SI{389703}{km}\\ 
            Time unit (TU) & $t^\ast$ & \SI{382981}{s}\\
        \hline \hline
      \end{tabular}
\end{table}

Despite the highly non-linear equations of motion governing the motion of a third body in the CR3BP, there still exist periodic orbits of interest for cislunar SDA, such as the halo orbit families, that have value in mission design. There are a wide variety of different orbit families, some occurring around the Lagrange points. These are the orbits that will be introduced in Sec.~\ref{sec:comp_ex}. Periodic orbits can be constructed through the use of differential correction methods applied to initial conditions produced through approximated analytical solutions \cite{richardson1980}.

\subsection{Extended Kalman Filter}
\label{sec:ekf} 
We perform orbit determination utilizing the EKF to track a target's position and velocity. This filter is chosen for its reduced computational cost over other filters. The process of the EKF is outlined in the following section. The state variables $\bm{x}$ of the system are shown in Eq.~\eqref{eq:statevars}:
\begin{equation}
    \renewcommand{\arraystretch}{0.8}
    \bm{x}=\left[ x,y,z,\Dot{x},\Dot{y},\Dot{z} \right]^T= \begin{bmatrix}
        \bm{x}_\text{pos} \\
        \bm{x}_\text{vel}
    \end{bmatrix}
    \label{eq:statevars}
\end{equation}
where $\bm{x}_\text{pos}=[x,y,z]^T$ and $\bm{x}_\text{vel}=[\Dot{x},\Dot{y},\Dot{z}]^T$ represent vectors for the position and velocity of the target, respectively.

The EKF requires defining the system dynamics and measurement model that constitute the following non-linear filtering problem~\cite{carpenter2018navigation}:
\begin{subequations}
    \label{eq:filter}
    \begin{alignat}{2}
    \Dot{\bm{x}}&=\bm{f}(\bm{x})+\bm{w}, \qquad  &&\text{where} \ \ \bm{w}\sim \mathcal{N}(\bm{0},\bm{Q})\label{eq:dynamics}\\
    \bm{m}_{k}&=\bm{h}(\bm{x}_{k})+\bm{v}_{k}, &&\text{where} \ \ \bm{v}_k \sim \mathcal{N}(\bm{0},\bm{R}_k) \label{eq:measurement}
    \end{alignat}
\end{subequations}
where $k$ represents the index for an iteration of the EKF, and $\mathcal{N}$ denotes the multivariate Gaussian distribution. In the filtering problem~\eqref{eq:filter}, the non-linear system dynamics with uncertainties are modeled with Eq.~\eqref{eq:dynamics}, using the CR3BP equations of motion from Eq.~\eqref{EOM}. Here, $\bm{f}(\bm{x})$ represents the system dynamics, given the state vector $\bm{x}$, while $\bm{w}$ is the dynamical system error, which is assumed to be a zero mean Gaussian with covariance $\bm{Q}$. Equation~\eqref{eq:measurement} is the measurement model wherein $\bm{m}_k$ represents the output of the model (a vector of measurements), $\bm{h}(\bm{x}_k)$ refers to the measurement function (in this paper, angles-only measurements, outlined in Sec.~\ref{sec:meas_mod}), and $\bm{v}_k$ is the measurement noise, which is also assumed to be a zero mean Gaussian with covariance $\bm{R}_k$. Note that $\bm{R}_k$ is assumed to be constant in this paper, such that $\bm{R}_k=\bm{R}$. 
    
We initialize $\bm{Q}$, the process noise covariance matrix, through the use of State Noise Compensation (SNC) \cite{schutz2004stat,caudill2021relative}, a method for estimating acceleration that is not modeled in the dynamical model. The matrix is formulated through the following equation:
\begin{equation}
    \label{eq:snc}
    \bm{Q} = \bm{\Gamma}\bm{Q}_{\text{snc}}\bm{\Gamma}^T
\end{equation}
The matrix $\bm{Q}_{\text{snc}}$ is defined as:
    \begin{equation*}
        \renewcommand{\arraystretch}{0.8}
        \bm{Q}_{\text{snc}} = \begin{bmatrix}
            \sigma_\text{dyn}^2 & 0 & 0 \\
            0 & \sigma_\text{dyn}^2 & 0 \\
            0 & 0 & \sigma_\text{dyn}^2
        \end{bmatrix}
    \end{equation*}
where $\sigma_\text{dyn}$ is the standard deviation of the unmodeled acceleration, and $\sigma_\text{dyn}^2$ is its variance. In Eq. \eqref{eq:snc}, $\bm{\Gamma}$ transforms the unmodeled accelerations into values that can be mapped onto the six states of $\bm{x}$:
\begin{equation*}
    \renewcommand{\arraystretch}{0.8}
    \bm{\Gamma} = \begin{bmatrix}
        \frac{\Delta t^2}{2} & 0 & 0 \\
        0 & \frac{\Delta t^2}{2} & 0 \\
        0 & 0 & \frac{\Delta t^2}{2} \\
        \Delta t & 0 & 0 \\
        0 & \Delta t & 0 \\
        0 & 0 & \Delta t \\
    \end{bmatrix}
\end{equation*}
where $\Delta t$ represents the time step between two consecutive iterations of the EKF. We assume this matrix is constant.

\subsubsection{Prediction Step} 
The prediction step of an iteration of the EKF simultaneously propagates the estimated state vector $\hat{\bm{x}}$ and error covariance matrix $\bm{P}$:
\begin{subequations}
\begin{align}
    \Dot{\hat{\bm{x}}}&=\bm{f}(\hat{\bm{x}}) \label{eq:state_fn} \\
    \Dot{\bm{P}}&=\bm{A}\bm{P}+\bm{P}\bm{A}^{T}+\bm{Q} \label{eq:riccati}
\end{align}
\end{subequations} 
Through integration, Eqs.~\eqref{eq:state_fn} and \eqref{eq:riccati} return $\hat{\bm{x}}=\hat{\bm{x}}_{k}^-$ and $\bm{P}=\bm{P}_{k}^-$, respectively, with $\hat{\bm{x}}^-_{k}$ and $\bm{P}^-_{k}$ denoting the \textit{a priori} estimates of the state vector and covariance matrix, respectively, before the estimate correction is performed. This paper utilizes a Riccati equation for numerical propagation of $\bm{P}$ in the EKF algorithm \cite{lefferts1982kalman,carpenter2018navigation}. In Eq.~\eqref{eq:riccati}, $\bm{A}$ is the Jacobian of the system dynamics with respect to the state vector:

\begin{equation*}
    \bm{A}=\frac{\partial{\bm{f}}}{\partial{\bm{x}}}\bigg|_{{\bm{x}=\hat{\bm{x}}}}
\end{equation*}

\subsubsection{Correction Step}
The correction step of an iteration of the EKF fuses the \textit{a priori} estimate from the prediction step with a set of measurements to obtain the \textit{a posteriori} estimate of the state of the target \cite{schutz2004stat}. First, the Kalman gain $\bm{K}$ is computed in Eq.~\eqref{eq:kalman}:
\begin{equation}
    \bm{K}_{k}=\bm{P}^{-}_{k}\bm{H}^{T}_{k}(\bm{H}_{k}\bm{P}^{-}_{k}\bm{H}^{T}_{k}+\bm{R}_{k})^{-1} \label{eq:kalman} \\
\end{equation}
In Eq.~\eqref{eq:kalman}, $\bm{P}^{-}$ represents the \textit{a priori} estimate error covariance, and $\bm{H}$ is the Jacobian of the measurement function, referred to as the measurement Jacobian hereafter, and will be discussed in further detail in Sec.~\ref{sec:meas_mod}. The matrix $\bm{H}$ is computed in Eq.~\ref{eq:hmj}.
\begin{equation}
    \label{eq:hmj}
    \bm{H}=\frac{\partial{\bm{h}}}{\partial{\bm{x}}}\bigg|_{\bm{x}=\hat{\bm{x}}^{-}}
\end{equation}
    
With the Kalman gain, Eq.~\eqref{eq:stateupdate} corrects the state estimate where $\hat{\bm{x}}^{+}_{k}$ is the \textit{a posteriori} state estimate after the measurement correction. In this step, $\bm{m}_k$ represents the actual measurement, and $\bm{h}(\hat{\bm{x}}^{-}_{k})$ represents the predicted measurement, using the \textit{a priori} state estimate. Lastly, Eq.~\eqref{eq:covarianceupdate} corrects the covariance matrix, returning $\bm{P}_{k}^{+}$, the \textit{a posteriori} covariance matrix. 
\begin{subequations}
\begin{align}
    \hat{\bm{x}}_{k}^{+}&=\hat{\bm{x}}^{-}_{k}+\bm{K}_{k}(\bm{m}_{k}-\bm{h}(\hat{\bm{x}}^{-}_{k})) \label{eq:stateupdate}\\
    \bm{P}_{k}^{+}&=(\bm{I}-\bm{K}_{k}\bm{H}_{k})\bm{P}^{-}_{k} \label{eq:covarianceupdate}
\end{align}
\end{subequations}    
where $\bm{I}$ is an identity matrix with a conformable dimension.

\subsection{Measurement Model} \label{sec:meas_mod}
In this paper, we use angles-only measurements to map observations to the states of the target RSOs.

\subsubsection{Angles-Only Measurements}
The measurement model consists of two angular measurements; azimuth and elevation, denoted as $\alpha$ and $\varepsilon$, respectively. Given an observer, a target, and a reference point, two relative position vectors are required to calculate these angles. The relative position vectors are $\bm{\gamma}$, which points from the observer to a reference point, the Moon in this paper, and $\bm{\rho}$, which points from the observer to the target. They are calculated as follows \cite{block2022cislunar}:
\begin{subequations}
    \renewcommand{\arraystretch}{0.8}
    \begin{align}
    \bm{\gamma} &= [1-\mu,0,0]^T-[x_{o},y_{o},z_{o}]^T = [\gamma_{x},\gamma_{y},\gamma_{z}]^T \label{eq:gamma} \\
    \rhobar &= [x,y,z]^T-[x_{o},y_{o},z_{o}]^T=[\rho_{x},\rho_{y},\rho_{z}]^T \label{eq:rho}
    \end{align}
\end{subequations}
where $x_{o}$, $y_{o}$, and $z_{o}$ denote the coordinates of the observer.

\begin{figure}[htb]
	\centering
        \includegraphics[width=.35\linewidth]{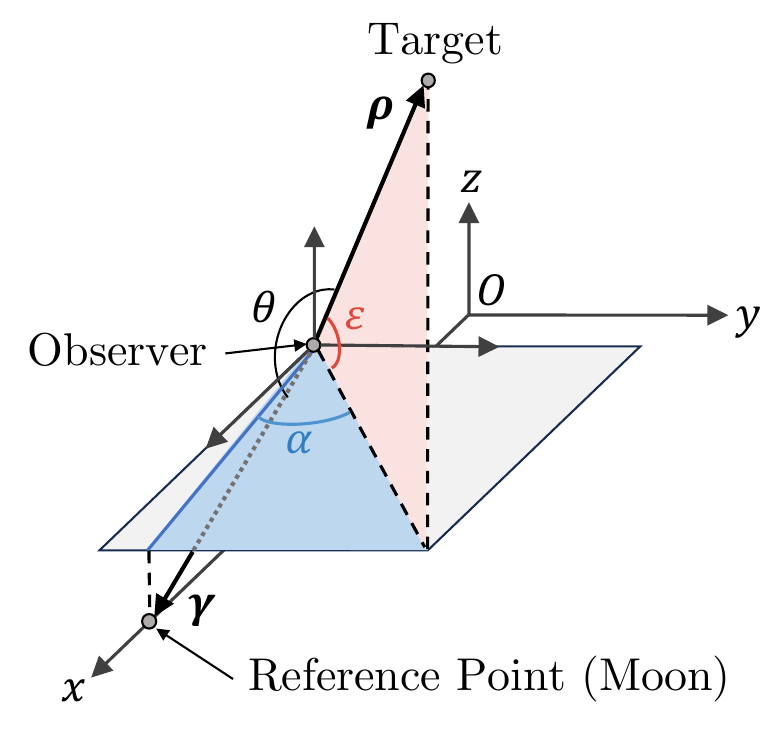}
        \caption{Visualization of the angles and vectors in the measurement model.}
	\label{fig:meas_model}
\end{figure}

The azimuth, $\alpha$, and the elevation, $\varepsilon$, of the target RSO with respect to the observer are then computed as follows:
\begin{subequations}
    \begin{align*}
        \alpha &= \arccos \left(\frac{\gamma_{x}\rho_{x}+\gamma_{y}\rho_{y}}{\sqrt{\gamma_{x}^{2}+\gamma_{y}^{2}}\sqrt{\rho_{x}^{2}+\rho_{y}^{2}}}\right)  \\
        \varepsilon &= \arccos \left(\frac{\gamma_{x}\rho_{x}+\gamma_{z}\rho_{z}}{\sqrt{\gamma_{x}^{2}+\gamma_{z}^{2}}\sqrt{\rho_{x}^{2}+\rho_{z}^{2}}}\right) 
    \end{align*}
\end{subequations}
These two angular measurements constitute the output of the measurement function such that $\bm{h}=[\alpha, \varepsilon]^T$. These angles are shown in Fig.~\ref{fig:meas_model}.

To compute the matrix $\bm{H}$, one must calculate the Jacobian of $\bm{h}$ with respect to $\bm{x}$ as shown in Eq.~\eqref{eq:hmj}:
\begin{equation*}
        \bm{H}=\frac{\partial{\bm{h}}}{\partial{\bm{x}}}=\frac{\partial{\bm{h}}}{\partial{\bm{\rho}}}\frac{\partial{\bm{\rho}}}{\partial{\bm{x}}}
\end{equation*}
where
\begin{equation*}
    \renewcommand{\arraystretch}{0.8}
    \frac{\partial{\bm{h}}}{\partial{\bm{\rho}}}
                        =\begin{bmatrix}
                            \frac{\partial{\alpha}}{\partial{\rho_{x}}} & \frac{\partial{\alpha}}{\partial{\rho_{y}}} & \frac{\partial{\alpha}}{\partial{\rho_{z}}} \\
                            \frac{\partial{\varepsilon}}{\partial{\rho_{x}}} & \frac{\partial{\varepsilon}}{\partial{\rho_{y}}} & \frac{\partial{\varepsilon}}{\partial{\rho_{z}}}
                        \end{bmatrix} \quad \text{and} \quad
                        \frac{\partial{\rhobar}}{\partial{\bm{x}}}=\begin{bmatrix}
            1 & 0 & 0 & 0 & 0 & 0 \\
            0 & 1 & 0 & 0 & 0 & 0 \\
            0 & 0 & 1 & 0 & 0 & 0
    \end{bmatrix}
\end{equation*}
\subsubsection{Measurement Constraints}
A few constraints are imposed on the observers' ability to make these measurements: (1) the maximum sensor range constraint and (2) the lunar and terrestrial exclusion/occlusion constraint.
 
The first constraint checks that the range (\textit{i.e.}, the norm of $\bm{\rho}$) does not exceed a maximum value, modeling hardware limitations of the sensor used by an observer. The second ensures that the target is observable by checking the exclusion and occlusion of the target by the two primary bodies in the system (the Earth and the Moon). The angle $\theta_p$, depicted in Fig.~\ref{fig:meas_model}, represents the angle between $\bm{\gamma}_p$, the relative position vector from the observer to primary body \textit{p} and the $\rhobar$ vector. It is computed from the following dot product formula:
\begin{equation*}
    \theta_p=\arccos \left(\frac{\bm{\gamma}_p\cdot\bm{\rho}}{\lVert \bm{\gamma}_p \rVert \lVert \bm{\rho} \rVert}\right)
\end{equation*}
The angle $\omega_p$ represents the angle between the relative position vector from the observer to the center of primary body \textit{p} and the relative position vector pointing from the observer that is tangent to the surface of the same primary body. It is defined as:
\begin{equation*}
    \omega_p=\arctan \left(\frac{\frac{r_{p}}{\gamma_p}\sqrt{\gamma_p^{2}-r_{p}^{2}}}{\gamma_p-\frac{r_{p}^2}{\gamma_p}}\right)
\end{equation*}
where $r_{p}$ represents the radius of primary body \textit{p} (in the following experiments, the Earth and the Moon, with radii of \SI{6378.1}{km} and \SI{1737.1}{km}, respectively), and $\gamma_p =||\bm{\gamma}_p||$ is the Euclidean norm of the relative position vector pointing from the observer to primary body \textit{p}. Therefore, the lunar and terrestrial exclusion/occlusion constraint is consequently an inequality that checks if $\theta_p \ge \omega_p$ for all $p$.

If an observer does not satisfy either of these constraints during a time step of the EKF, the measurements made by that observer are discarded from $\bm{h}$ and $\bm{H}$, and the correction step is performed with partial information. If no observers can make a measurement, the correction step is not performed.

\subsection{Sensor Tasking Procedures}
\label{sec:task}
To operate a constellation of space-based observers, operators need to determine when their observers take measurements, and of which targets they make observations, a procedure known as sensor tasking, taking into account such factors as each sensor's individual cadence. In this paper, we test three different sensor tasking procedures for a single target to assess how this affects the optimization problem proposed in Sec.~\ref{sec:optimization}.

\begin{figure}[htbp]
\centering
\begin{subfigure}{0.24\textwidth}
\includegraphics[width=\linewidth]{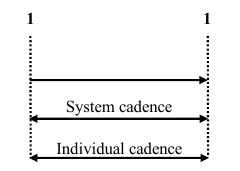}
\caption{Single observer.}
\end{subfigure}
\begin{subfigure}{0.24\textwidth}
\includegraphics[width=\linewidth]{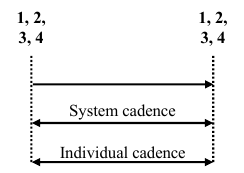}
\caption{STP-A.}
\end{subfigure}
\begin{subfigure}{0.24\textwidth}
\includegraphics[width=\linewidth]{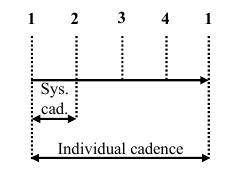}
\caption{STP-B.}
\end{subfigure}
\begin{subfigure}{0.24\textwidth}
\includegraphics[width=\linewidth]{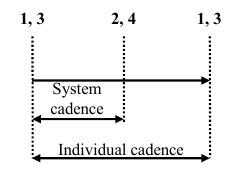}
\caption{STP-C.}
\end{subfigure}
\caption{Visual representation of sensor tasking procedures, where each number represents an individual observer.}
\label{fig:tasking_procedures}
\end{figure}

Three sensor tasking procedures are used in this paper, and the assumption is made that all sensors have identical individual cadences. Sensor tasking procedure A (STP-A) tasks all observers to take measurements in unison, at the start of the simulation, and then after every cadence of their sensors. The system cadence is equal to an observer's individual cadence. Sensor tasking procedure B (STP-B) has each observer take a measurement sequentially, at equally spaced times. Therefore, the system cadence is calculated by dividing the individual cadence of a sensor by the number of sensors in the system. Sensor tasking procedure C (STP-C) splits the constellation into two different groups. In the first time step, the first grouping takes a measurement of the target. The second group then takes a measurement of the target halfway through the first grouping's individual cadences. Thus, the system cadence is calculated as half of the individual sensor cadence. Figure~\ref{fig:tasking_procedures} is an illustrative case of how a single observer and the three tasking procedures operate when four observers compose a constellation.

\subsection{Multi-observer Placement Optimization}
\label{sec:optimization}
The accuracy of the state estimation of a target RSO by a constellation, as compared to its ground truth state, is contingent upon not only the system dynamics model-plant mismatch, but also the quality and quantity of the measurements provided. In turn, this is largely determined by the relative positioning of observers. In this paper, the aim is to optimize the placement of multiple observers such that the constellation provides the best state estimation performance.

The orbits in which each target will reside are unknown, therefore the optimization can be formulated as a stochastic optimization problem, however, this requires knowledge of the probabilities of which orbits will be utilized, but data regarding this is not readily available. In response to this, we formulate the problem deterministically, selecting a small, albeit representative, sampling of targets from each orbit family of interest, referred to as the \textit{optimization target set}. Following a convergence by the algorithm used to process the optimization problem, a post-optimization analysis is performed through the simulation of the optimal constellation against a larger set of target orbits, referred to as the \textit{validation target set}.

In designing a multi-observer constellation system for cislunar SDA, a set of initial conditions for a possible observer is treated as an orbital slot, and the existence of an observer in this slot is binary. Consequently, the problem of placing observers in a given set of orbital slots is posed as a combinatorial optimization problem. As the number of candidate orbits increases linearly, the number of possible constellations grows exponentially. Performing an exhaustive search within such a vast solution space becomes computationally challenging. Furthermore, the process of getting from the binary decision variables to an RMSE value requires non-linear calculations, meaning that the optimization problem itself can be considered non-linear. This requires the use of algorithms that are suited for this category of problem, such as a GA or particle swarm optimization.

Let $\mathcal{J}$ be the set of orbital slots that observers can occupy, indexed by $j$, and let $\mathcal{U}$ be the optimization target set, indexed by $u$. The objective of the multi-observer placement problem is to position $N$ observers in $\mathcal{J}$ such that the average state estimation performance of a constellation against the optimization target set $\mathcal{U}$ is maximized. The decision variables are defined as follows:
\begin{equation*}
    X_{j} = \begin{cases}
1, &\text{if an observer occupies orbital slot $j$ in $\mathcal{J}$} \\
0, &\text{otherwise}
\end{cases}
\end{equation*}

Root-mean-square-error (RMSE) of a target RSO is chosen as the figure of merit to assess how well a cislunar SDA constellation performs state estimation, or how accurately the constellation estimates a target's trajectory. It is computed for the position of the target with the notation according to Eq.~\eqref{eq:statevars}. Note that this loss function is a non-linear function of the decision variables $\bm{X}$ because the EKF directly affects it. The position RMSE loss function, $\mathcal{L}$, is defined in Eq.~\eqref{eq:loss_fn}.
\begin{equation}
    \mathcal{L}=\sqrt{\frac{\sum_{k=1}^{n}\left(\lVert \hat{\bm{x}}_{\text{pos},k}^+ \rVert-\lVert \bm{x}_{\text{pos}, k} \rVert\right)^2}{n}}
    \label{eq:loss_fn}
\end{equation}
where $n$ is the number of time steps within the simulation, $\hat{\bm{x}}_{\text{pos}, k}^+$ refers to the \textit{a posteriori} position vector of the target at each time step $k$ of simulation, generated by the EKF, and $\bm{x}_{\text{pos}, k}$ is the ground truth position of the target at each time step.

With the preceding notations, the multi-observer placement problem can be formulated as a non-linear program as follows:
\begin{subequations}
\label{eq:opt_problem}
\begin{align}
	\underset{\bm{X}}{\text{minimize}} \quad & \frac{1}{|\mathcal{U}|}\sum_{u\in \mathcal{U}} \mathcal{L}_u(\bm{X}) \label{eq:RMSE_opt} \\
	\text{subject to} \quad
	& \sum_{j\in\mathcal{J}} X_j = N\label{eq:c1} \\
	& X_j \in \{0, 1\}, \quad \forall j\in \mathcal{J} \label{eq:c2}
\end{align}
\end{subequations}
In Problem~\eqref{eq:opt_problem}, the goal is to minimize the objective function~\eqref{eq:RMSE_opt}, which represents the average position RMSE of the optimization target set (thus maximizing the state estimation performance). Here, $\mathcal{L}_u(\bm{X})$ serves as the position loss function for target $u \in \mathcal{U}$. Constraint~\eqref{eq:c1} is a cardinality constraint that limits the number of observers to $N$. Constraints~\eqref{eq:c2} define the domain of the decision variables.

\subsection{Post-optimization Analysis}
\label{sec:postoptanalysis}
The multi-observer placement optimization problem, in Problem~\eqref{eq:opt_problem}, optimizes the placement of observers by only considering the state estimation results of the optimization target set. The optimization target set is created by sampling different orbit families, so it is important to provide a post-optimization analysis of how the optimized constellation performs against a wider set of targets, the validation target set, to ensure that the optimized constellation can perform against targets it is not optimized to track. The validation target set contains orbits within the same orbit families as the optimization target set, as well as the optimization target set itself. All of these targets are tracked using the EKF, and their results are reviewed on the statistical level, providing an analysis of how a constellation performs against different sets of targets.

\subsection{Framework Outline}
\begin{figure}[htb]
	\centering
        \includegraphics[width=\linewidth]{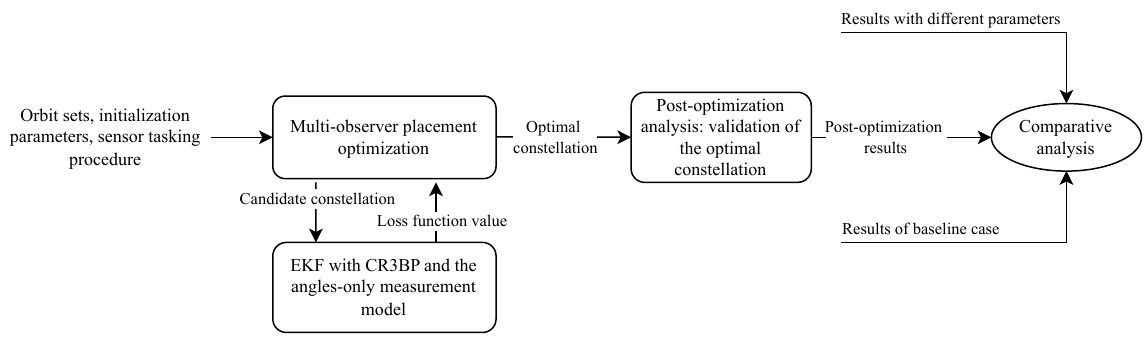}
	\caption{Block diagram outlining the Cislunar SDA Constellation Design and Analysis Framework.}
	\label{fig:blockd}
\end{figure}

Through the use of these components, the Cislunar SDA Constellation Design and Analysis Framework operates as outlined in Fig.~\ref{fig:blockd}. The framework requires the definition of a set of parameters, including the orbits of the targets and candidate observers, the parameters of the sensors (\textit{e.g.}, their individual cadences and their uncertainty), the sensor tasking procedure used by the constellation, and the parameters of the EKF and CR3BP. The optimization algorithm attempts to solve the multi-observer placement problem, as shown in Problem~\eqref{eq:opt_problem}, by selecting possible constellations to use the EKF to track the optimization target set. Once there is convergence to a solution, the constellation is tested against the validation target set, testing the robustness of the solution against a wide variety of targets. It is important to note that the goal of this work is to demonstrate the effectiveness of embedding state estimation into the preliminary design of cislunar SDA constellations; the design can be adjusted by changing the fidelity of the modules within the framework and modifying the parameters used to initialize it.

\section{Computational Experiments}
\label{sec:comp_ex}

In Sec.~\ref{sec:comp_ex}, we present computational experiments to compare a single observer, the baseline, to three solutions of the Cislunar SDA Constellation Design and Analysis Framework that vary their sensor tasking procedure. Section~\ref{sec:ex_res} presents the results of the four different cases, while Sec.~\ref{sec:parameters} describes another set of experiments to test how the optimized constellation changes with varied sensor fidelity. 

\subsection{Experimental Setup}
\label{sec:ex_set}
The framework presented in the previous section is used with three different sensor tasking procedures outlined in Sec.~\ref{sec:task}. Additionally, a single observer in an orbit similar to that of Gateway is used as a baseline for comparison against multi-observer systems \cite{zimovan2023baseline}. Each procedure is run through the framework individually to test how the optimized constellation changes. The optimization algorithm used is MATLAB's built-in GA function, \texttt{ga()} \cite{matlab2022stats}, considering the non-linear nature of the problem, as discussed in Sec.~\ref{sec:optimization}. The default parameters for this function are used (crossover fraction of 0.8, maximum generation of $100*numberOfVariables$, maximum stall generations set to 50, no maximum stall or total time, and \textit{selectionstochunif} for the selection function). The population size is set to 50 for computational efficiency.

There are 12 orbit families in the target and observer sets. They are the Butterfly Northern (BNO), Butterfly Southern (BSO), Distant Retrograde (DRO), \lone Northern Halo (L1NHO), \lone Southern Halo (L1SHO), \ltwo Northern Halo (L2NHO), \ltwo Southern Halo (L2SHO), Low Prograde Eastern (LPEO), Low Prograde Western (LPWO), Resonant 1:1 (R1:1O), Resonant 2:1 (R2:1O) and Resonant 4:1 (R4:1O) orbit families. Additionally, transfer trajectories from Earth to an L1NHO (L1TT) are tracked, as it is important to provide SDA capabilities for spacecraft traveling to their final orbit. The number of orbits within each family is listed in Table~\ref{tab:familynums}. To generate these transfer trajectories, 19 linearly spaced points along a single \lone Northern Halo orbit are backpropagated to meet the $yz$-plane of the synodic frame. The state vector at this point is then treated as the initial condition for the simulation. A note here is made that these transfers are stable invariant manifolds, requiring an extra section of the trajectory to be found in order to connect them to an Earth parking orbit.

\begin{table}[htb]
\renewcommand{\arraystretch}{1}
\caption{Number of orbits investigated in each family.}
\label{tab:familynums}
\centering
\begin{tabular}{lr}
\hline \hline
Orbit family & Number of orbits \\
\hline
BNO          & 11               \\
BSO          & 11               \\
DRO          & 668              \\
L1NHO         & 21               \\
L1SHO        & 21               \\
L2NHO        & 170              \\
L2SHO        & 170              \\
LPEO          & 447              \\
LPWO          & 924              \\
R1:1O         & 276              \\
R2:1O         & 696              \\
R4:1O         & 539              \\
L1TT         & 19               \\
\hline \hline
\end{tabular}
\end{table}

A few metrics are used to size the validation target and observer orbit sets. A threshold for the orbits' stability index (SI) is set at \num{1.3}, as high SI values indicate greater instability and require a higher station-keeping cost to maintain the orbit \cite{guzzetti2017stationkeeping}. Hence orbits with high SI values are discarded from the orbit sets. Another constraint imposed on the target and observer sets requires that they have a maximum period of no greater than \SI{6.28}{TU}, selected through testing to size the orbit sets to be computationally efficient.

\begin{figure}[htbp]
     \centering
     \begin{subfigure}[b]{0.65\textwidth}
         \centering
         \includegraphics[width=\textwidth]{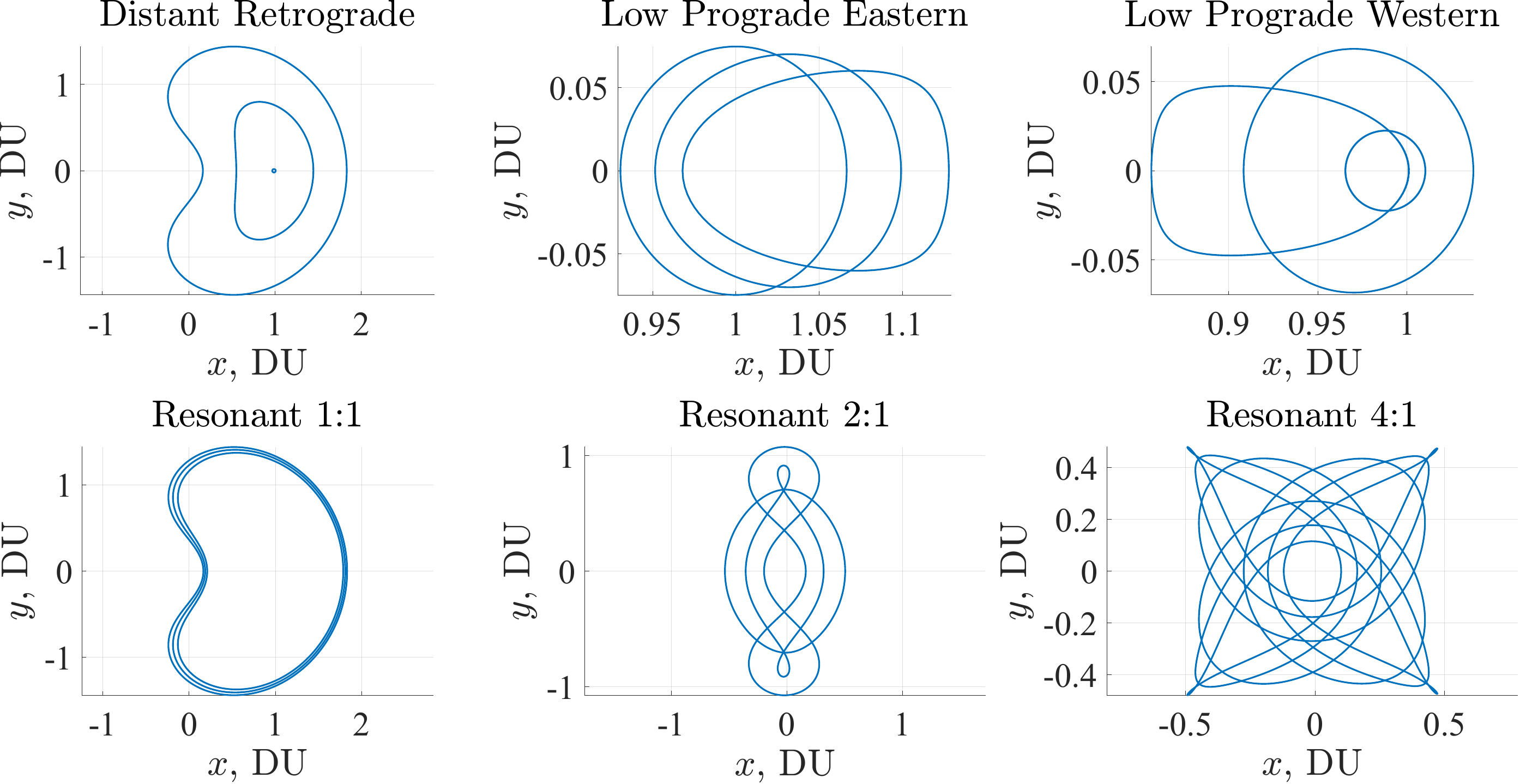}
         \caption{Planar orbit families.}
         \label{fig:PlanarOrbits}
     \end{subfigure}
     \begin{subfigure}[b]{0.65\textwidth}
         \centering
         \includegraphics[width=\textwidth]{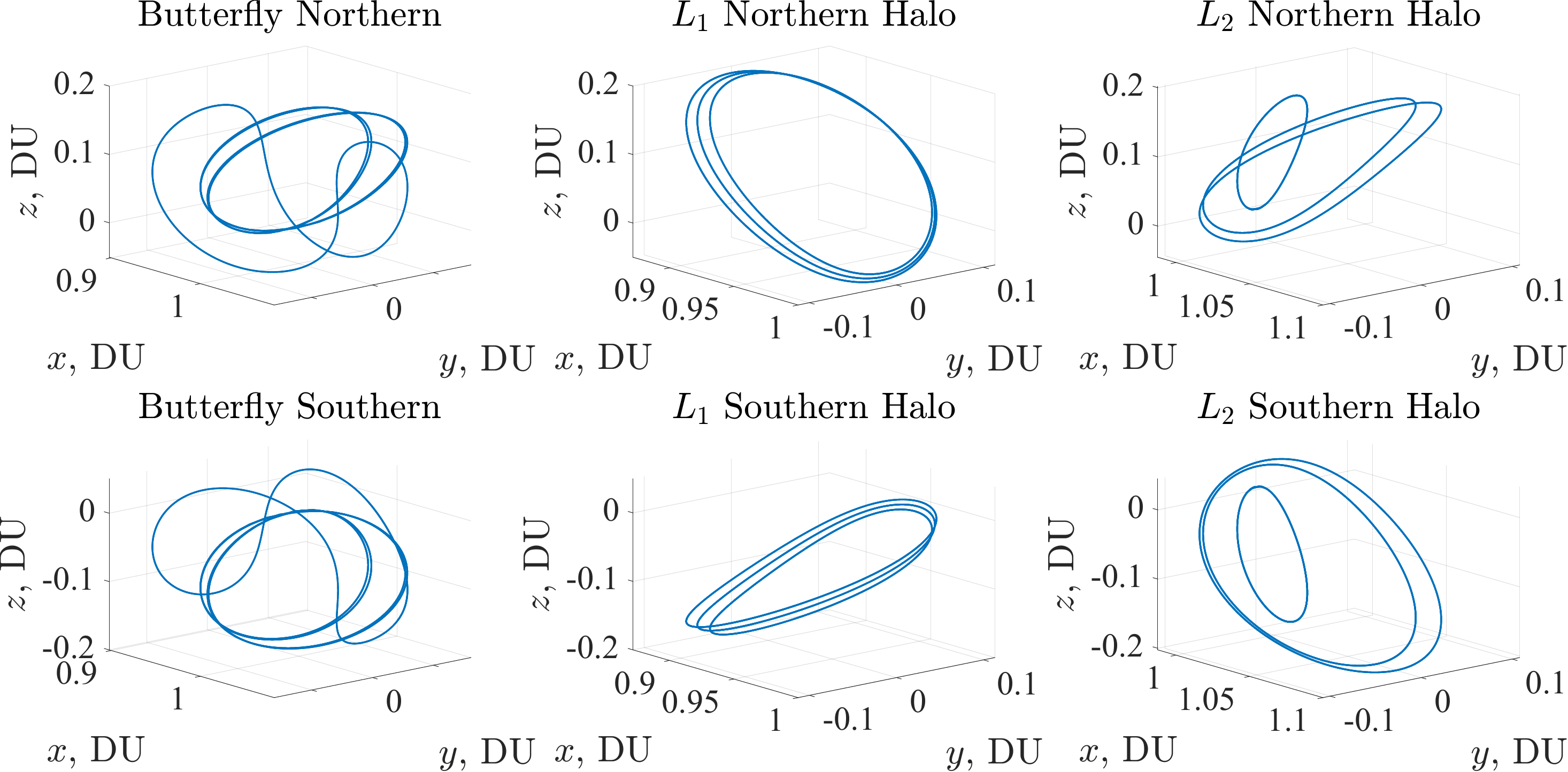}
         \caption{Non-planar orbit families.}
         \label{fig:NonPlanarOrbits}
     \end{subfigure}
     \begin{subfigure}[b]{.4\textwidth}
        \centering
        \includegraphics[width=0.7\textwidth]{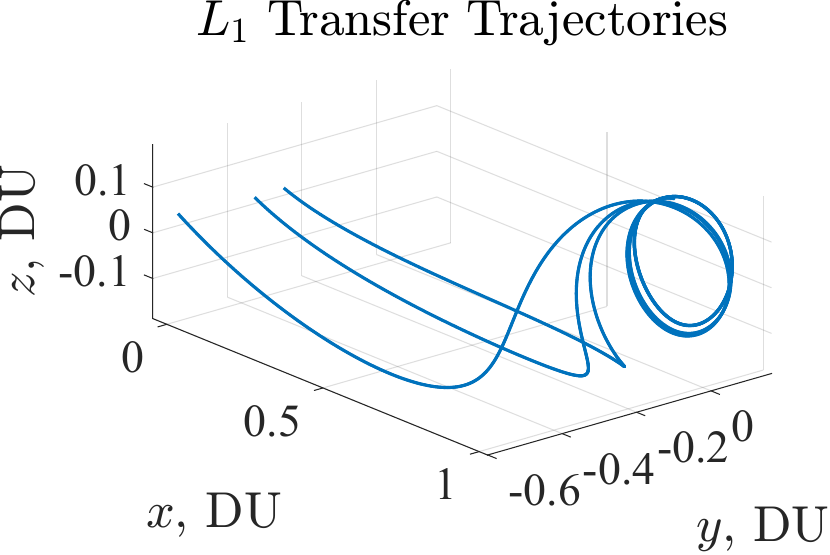}
         \caption{\lone Northern Halo Orbit Transfer Trajectories.}
         \label{fig:TransferTrajectories}
     \end{subfigure}
     \caption{Orbits used in the optimization target set.}
     \label{fig:FamPlots}
\end{figure}

The optimization target set samples the orbits within each family in the validation target set with the longest, shortest, and median period and three of the transfer trajectories, thus the total optimization target set has \num{39} orbits. The orbit families and initial conditions of all orbits in the optimization target set are included in Table~\ref{tab:OTS_ICs}, within Appendix~A. The optimization target set can be seen in Fig.~\ref{fig:FamPlots}. Due to how the targets are sampled, this figure is also representative of the validation target set (consisting of \num{3973} targets, across the 12 orbit families and trajectory set). For these experiments, the candidate observer orbit set is the same as the validation target set, but without the transfer trajectories, and with five equally spaced slots in each orbit, allowing for the GA to choose a constellation with multiple observers in one orbit, but phased differently. The number of observers is set to four.

The initial conditions for all target and observer orbits are obtained from Jet Propulsion Laboratory's (JPL) Three-Body Periodic Orbit repository \cite{jpl3bpdata}, which includes each orbit's initial conditions, orbital period, and SI values. Certain orbits from the JPL repository are removed from the orbit sets, as they intersect with the primary bodies, or do not exhibit the proper periodicity.  For both the optimization and validation target set, these initial conditions are propagated to a random point throughout their period, to simulate realistic initial conditions of an RSO in the regime. The state vector at this random point is then treated as the initial condition for the simulation, propagated using  MATLAB's ode45, an adaptive step-size ordinary differential equation solver, relying on a fourth and fifth order Runge-Kutta formula, the Dormand-Prince pair. The absolute and relative tolerances are both set as $\num{1e-12}$. The simulation horizon length is set to $\SI{8}{TU}$, which is equivalent to about $\SI{35.4}{days}$. This value is chosen, as it is greater than the maximum period of the orbits in the validation target set.

The initial conditions of a target given to the EKF are randomized by adding noise to each state, which is derived from a Gaussian distribution of each diagonal of the $\bm{Q}$ matrix. This added noise allows for initialization of the $\bm{P}$ matrix through the formulation outlined in Ref.~\cite{schneider2013how}:
    \begin{equation*}
        \bm{P}_0 = \text{diag}\left((\bm{\hat{x}}_0-\bm{x}_0)(\bm{\hat{x}}_0-\bm{x}_0)^T\right)
    \end{equation*}
where $\text{diag}(\cdot)$ returns a diagonal matrix with the values on the target matrix's main diagonal.

The individual sensor cadence is set to $\SI{0.02}{TU}$, dictating that each sensor is limited to taking a measurement every $\SI{127.7}{minutes}$, similar to other works \cite{hinga2023autonomous,greaves2021observation,qi2023investigation}. With this individual cadence limitation, STP-A, STP-B, and STP-C have respective system cadences of \SI{.02}{TU}, \SI{.005}{TU}, and \SI{.01}{TU}, respectively, as outlined in Fig.~\ref{fig:tasking_procedures}. Each sensor has a given maximum sensing distance of $\SI{500000}{km}$, and uncertainties of \SI{192.0118}{arcseconds} for the low-fidelity sensor experiments and \SI{26.7518}{arcseconds} for the high-fidelity sensor experiments, which is used to set $\bm{R}$, both originally used in Ref.~\cite{block2022cislunar}. To compute $\bm{Q}$, $\sigma_\text{dyn}$ is set to equal $\SI{1e-5}{DU/TU^2}$ to represent unmodeled acceleration due to solar radiation pressure \cite{vallado2001fundamentals}. When this value of $\sigma_\text{dyn}$ is utilized in creating the additive noise $\bm{w}$, it is negligible in the propagation and state estimation performance, so it is ignored and $\bm{w}$ is set to zero.

\subsection{Experimental Results}
\label{sec:ex_res}

\begin{figure}[htb]
     \centering
     \includegraphics[width=\textwidth]{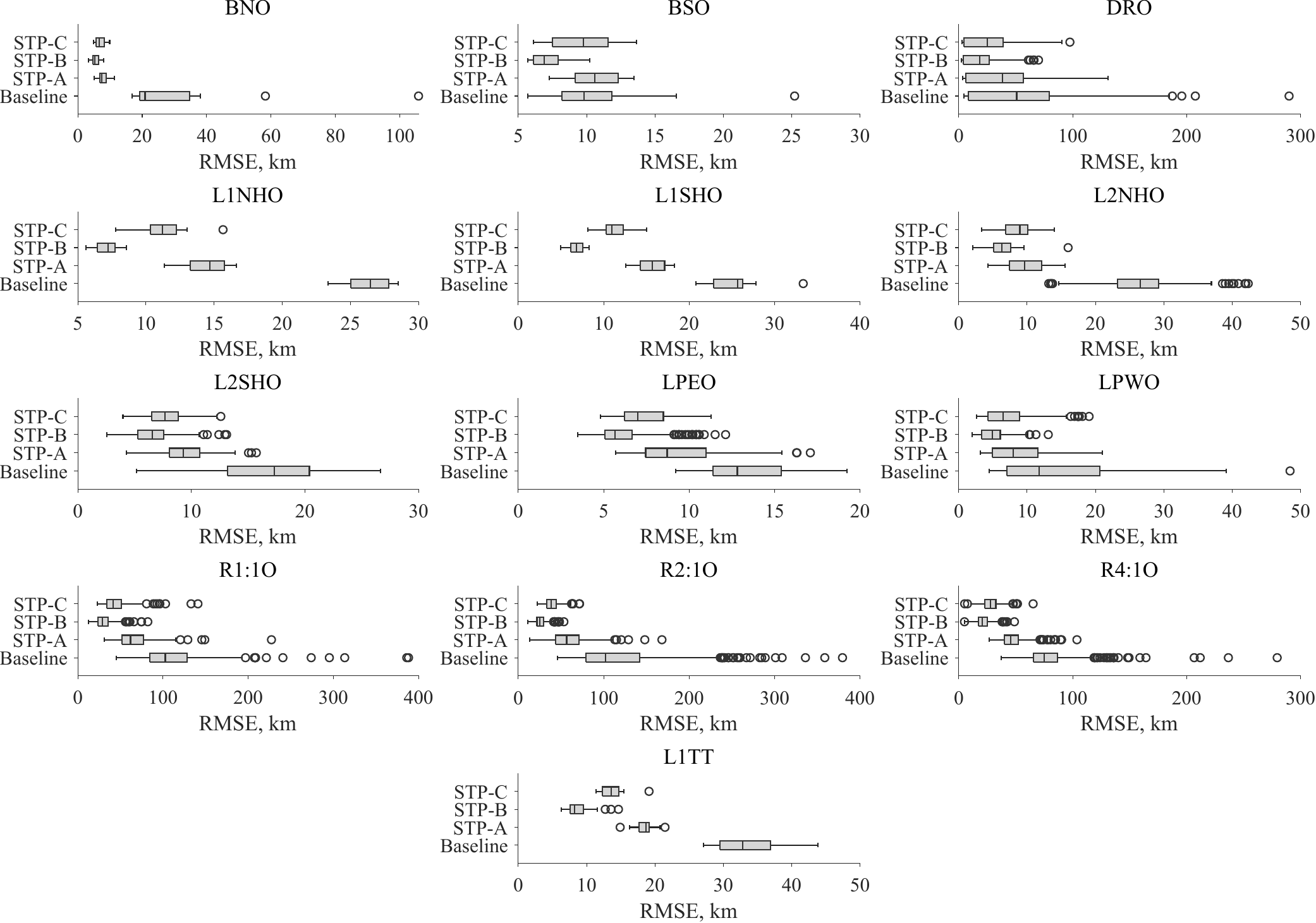}
     \caption{Box plots for each family and experimental case, for the low-fidelity sensors.}
     \label{fig:box_plots}
\end{figure}

In the first computational experiment comparing sensor tasking procedures, the GA returns loss function values of $\SI{20.81}{km}$, $\SI{10.98}{km}$, and $\SI{15.58}{km}$ for STP-A, STP-B and STP-C, respectively. A note is made here that these values are the average of the RMSEs of the optimization target set, and do not reflect the RMSE that is returned for any one target. In performing post-optimization analysis, each sensor tasking procedure's statistical performance against each orbit family can be seen in the box plots shown in Fig.~\ref{fig:box_plots}. These plots show the median, quartiles, maximum, and minimum RMSE values of the post-optimization validation. All three optimized constellations return better position RMSE averages and ranges than the baseline for all orbit families, except the BSO family, where the baseline performs slightly better than STP-A. The best-performing sensor tasking procedure is STP-B across all 13 orbit families, and the worst-performing is STP-A. While the GA is not guaranteed to find the global optimum, this could show that it is better to have a shorter system cadence, leading to more frequent measurements of a target.

\begin{figure}[htbp]
    \centering
    \begin{subfigure}[b]{0.7\textwidth}
        \centering
        \includegraphics[width=\linewidth]{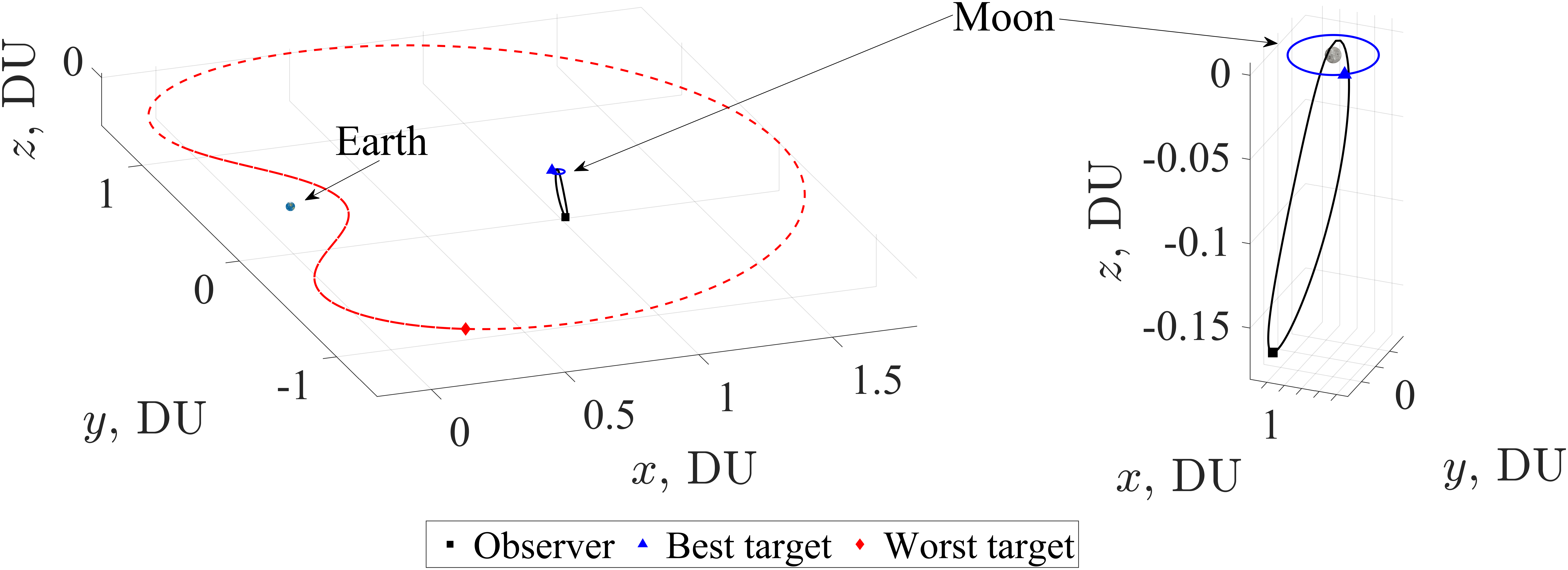}
        \caption{Baseline orbit.}
	    \label{fig:single_obs}
    \end{subfigure}
    \begin{subfigure}[b]{0.7\textwidth}
        \centering
        \includegraphics[width=\linewidth]{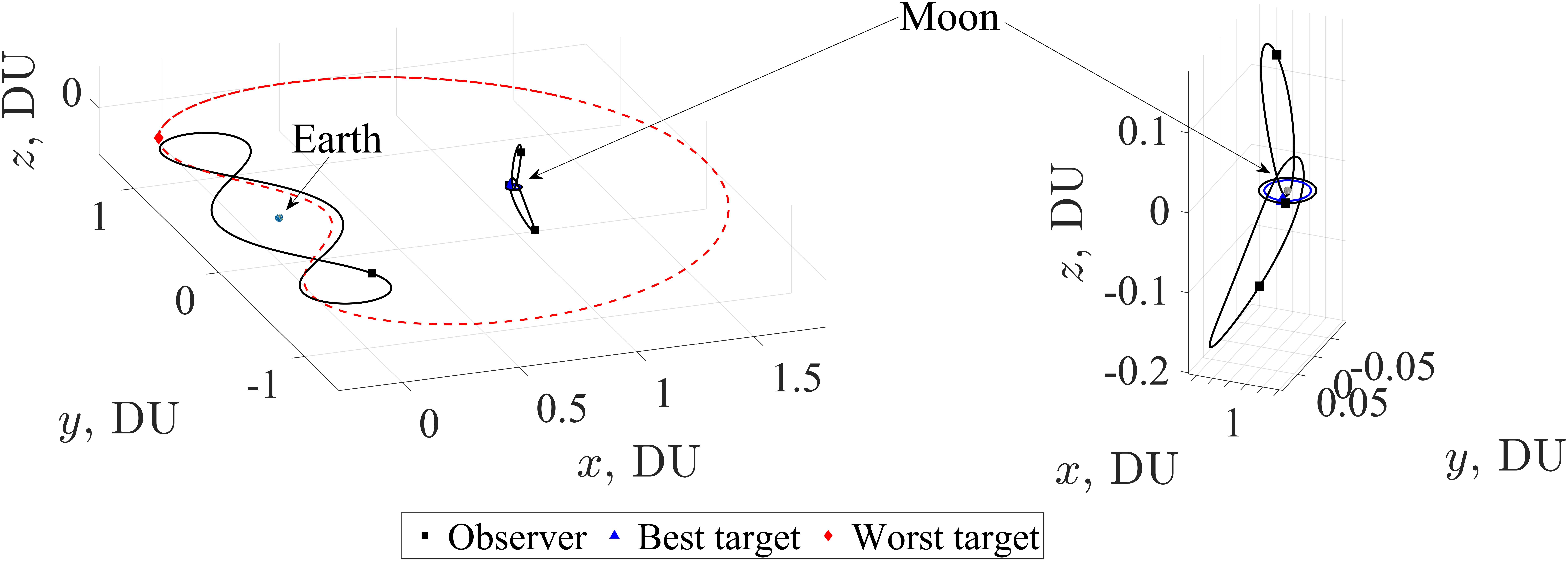}
        \caption{STP-A.}
	  \label{fig:TPA_constellation}
    \end{subfigure}
    \begin{subfigure}[b]{0.7\textwidth}
        \centering
        \includegraphics[width=\linewidth]{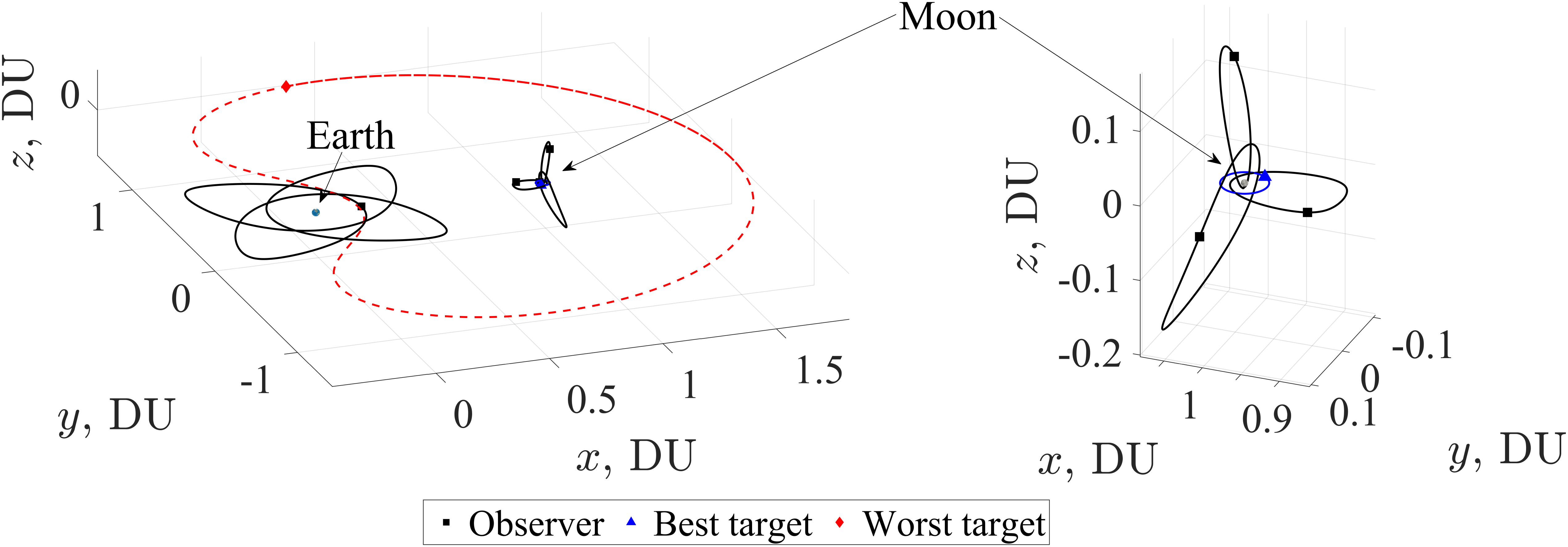}
        \caption{STP-B.}
	    \label{fig:TPB_constellation}
    \end{subfigure}
     \begin{subfigure}[b]{0.7\textwidth}
         \centering
        \includegraphics[width=\linewidth]{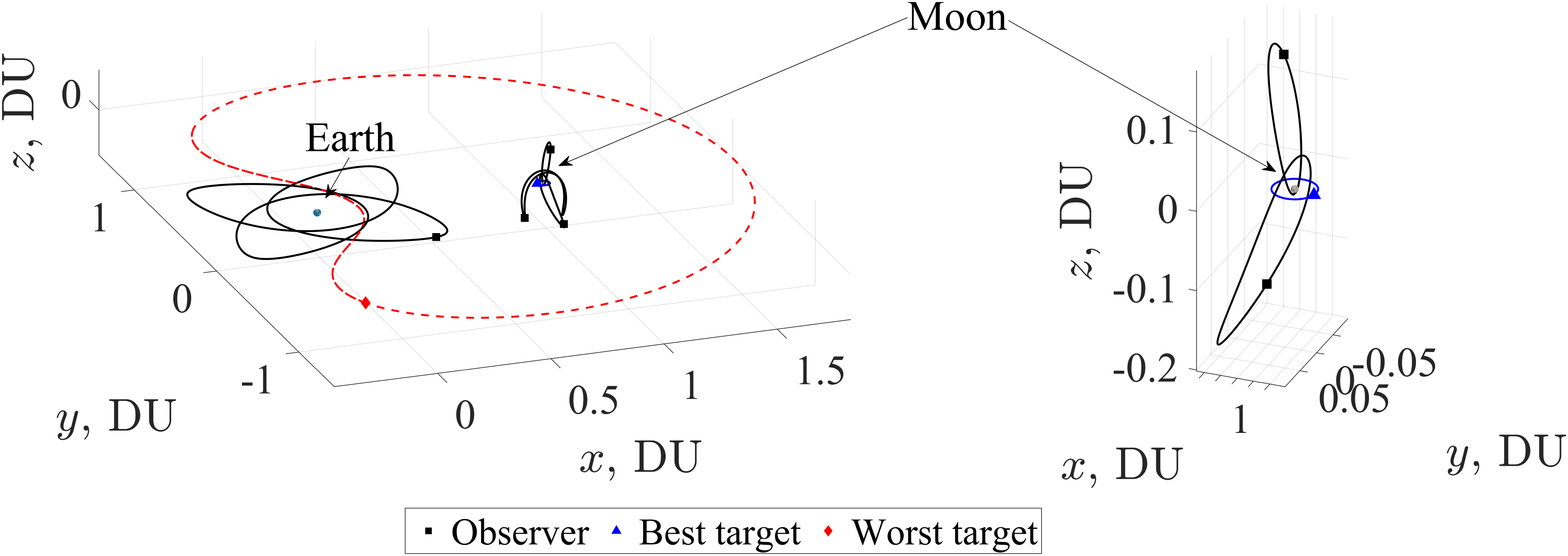}
        \caption{STP-C.}
	\label{fig:TPC_constellation}
     \end{subfigure}
        \caption{Observer orbits for each case, along with worst- and best-performing targets, and a close-up of near the Moon (right), for the low-fidelity observers.}
        \label{fig:optimal_constellations}
\end{figure}

The make-up of each optimized constellation can be seen in Fig.~\ref{fig:optimal_constellations}, while each observer's initial conditions and period/family information are found in Table~\ref{tab:optimal_ICs}, within Appendix~B. Across these optimized constellations, certain orbital slots are shared: all three tasking procedures share one orbital slot, in the L2NHO family. STP-A and STP-C also share another orbital slot in the L2SHO family. While STP-B does not share this same slot in the L2SHO family, it does share an observer in that family. The optimized constellations all additionally contain a resonant orbit. STP-B and STP-C contain observers in an R4:1O, the two orbits having periods within $\SI{0.0001}{TU}$ of each other (\textit{i.e.}, they are very similar orbits), and STP-A has an R2:1O.

Figure~\ref{fig:optimal_constellations} also visually shows the targets that each constellation performs the best and worst against (their pertinent information is in Table~\ref{tab:targ_ics_lofi}, within Appendix~B). The worst-performing cases are all in the R1:1O family. The best-performing targets are all in the LPWO family, returning values of \SI{4.47}{km}, for the baseline, \SI{3.15}{km} for STP-A, \SI{1.97}{km} for STP-B, and \SI{2.62}{km} for STP-C.

\begin{figure}[htbp]
     \centering
     \begin{subfigure}[b]{0.37\textwidth}
            \centering
            \includegraphics[width=\linewidth]{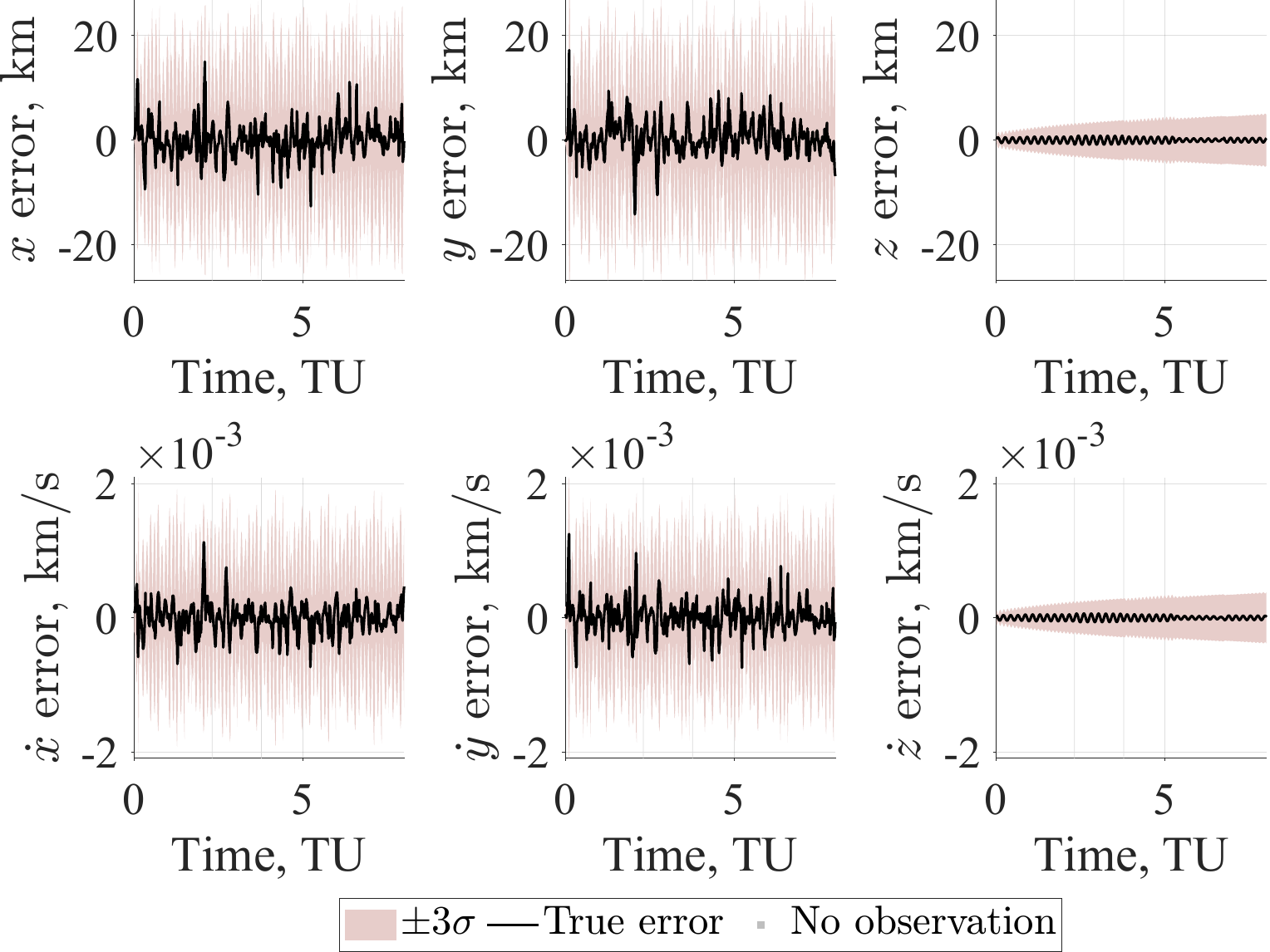}
            \caption{Baseline best case, an LPWO.}
	    \label{fig:3sigma_SO_best}
     \end{subfigure}
    \hspace{1em}
     \begin{subfigure}[b]{0.37\textwidth}
         \centering
        \includegraphics[width=\linewidth]{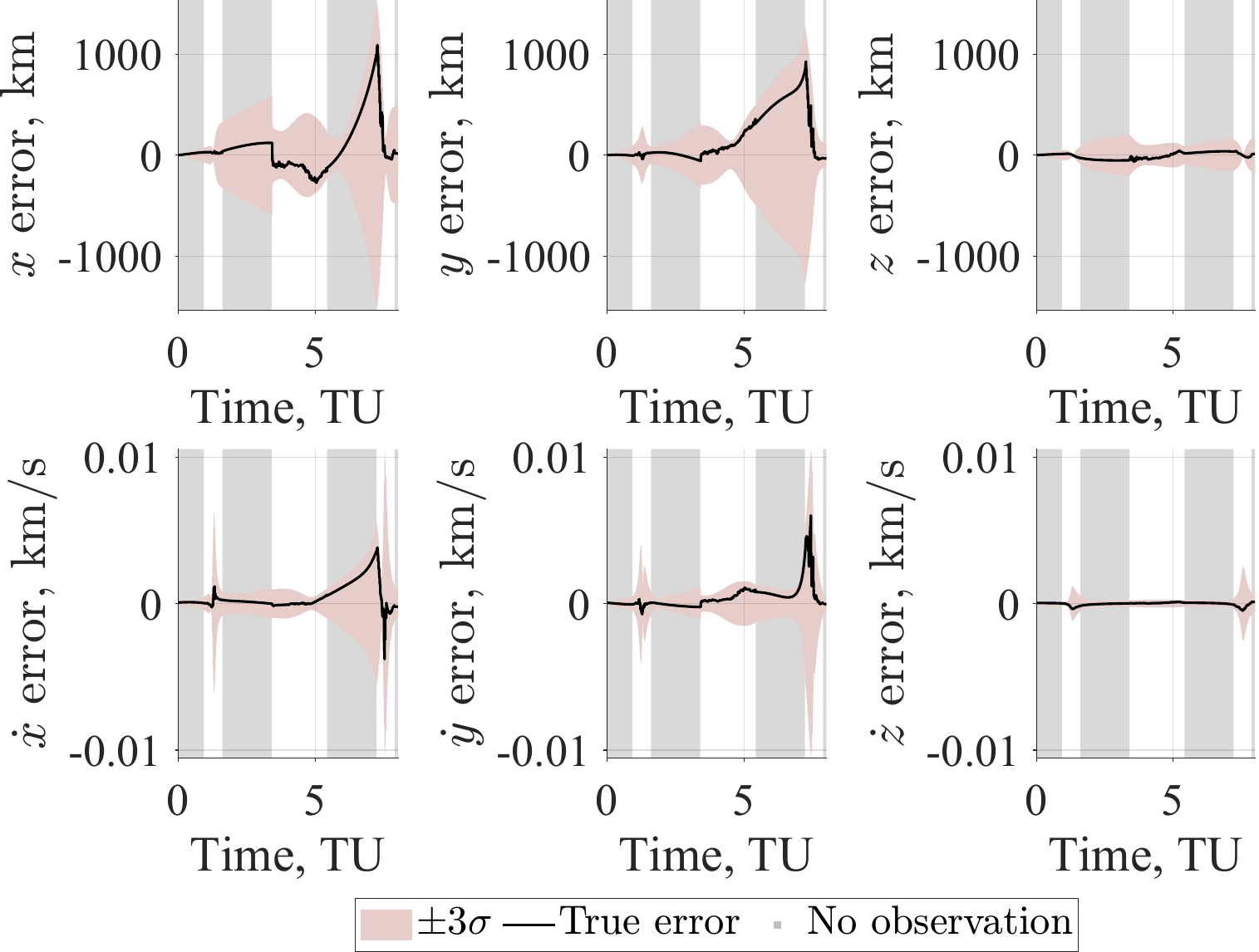}
        \caption{Baseline worst case, an R1:1O.}
	\label{fig:3sigma_SO_worst}
     \end{subfigure}
    \hspace{1em}
     \begin{subfigure}[b]{0.37\textwidth}
         \centering
         \includegraphics[width=\linewidth]{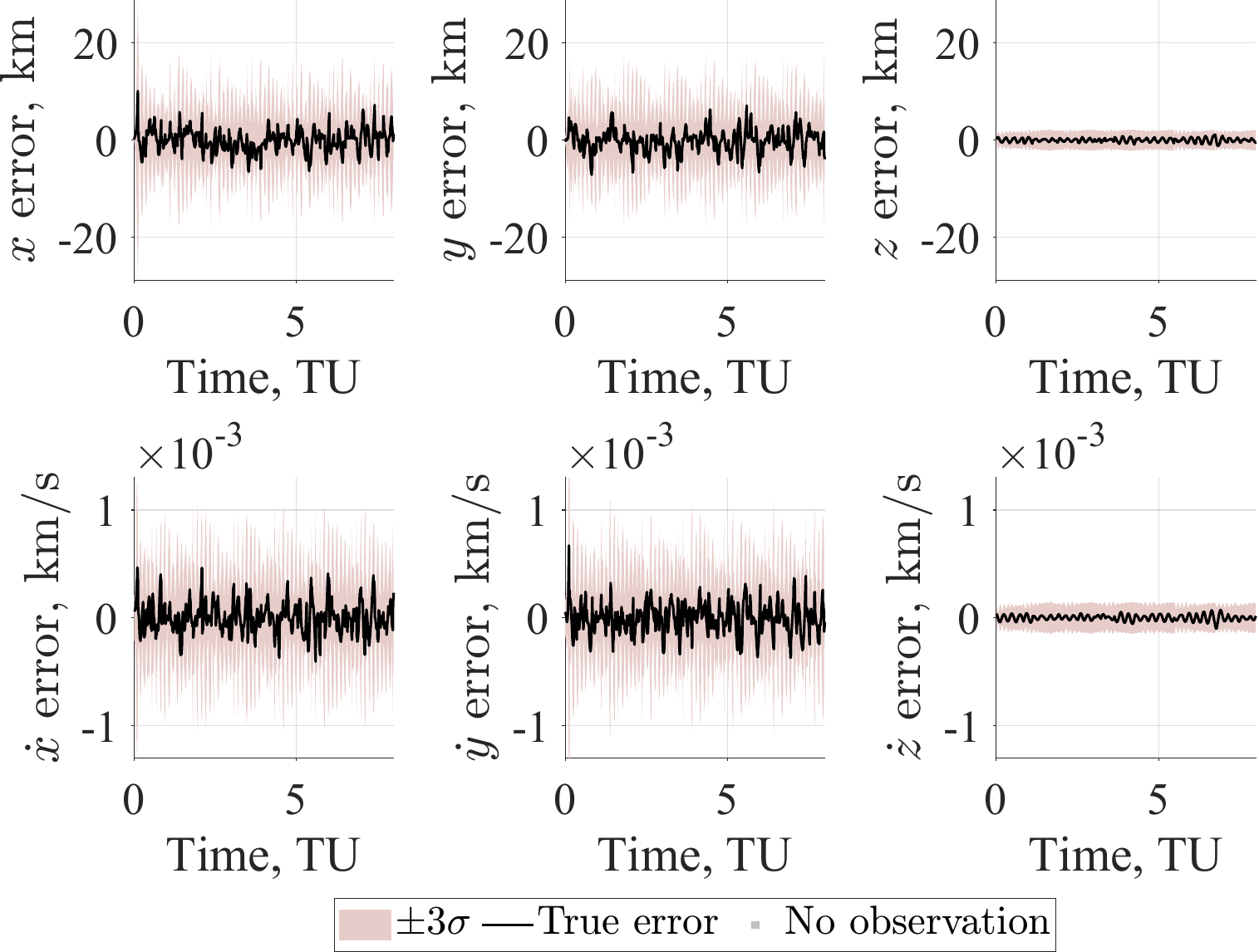}
            \caption{STP-A best case, an LPWO.}
	   \label{fig:3sigma_TPA_best}
     \end{subfigure}
    \hspace{1em}
     \begin{subfigure}[b]{0.37\textwidth}
         \centering
         \includegraphics[width=\linewidth]{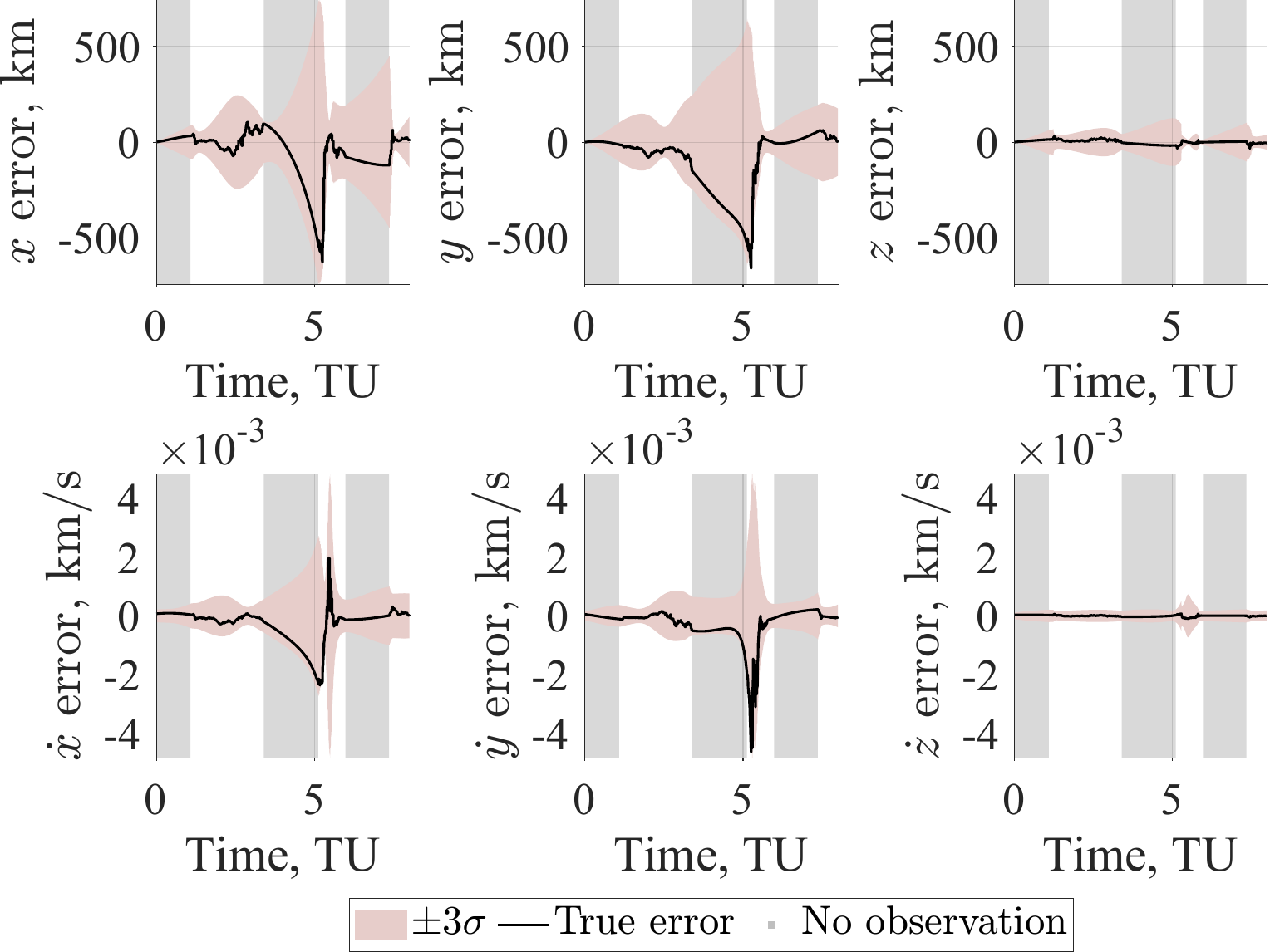}
        \caption{STP-A worst case, an R1:1O.}
	\label{fig:3sigma_TPA_worst}
     \end{subfigure}
    \hspace{1em}
     \begin{subfigure}[b]{0.37\textwidth}
            \centering
            \includegraphics[width=\linewidth]{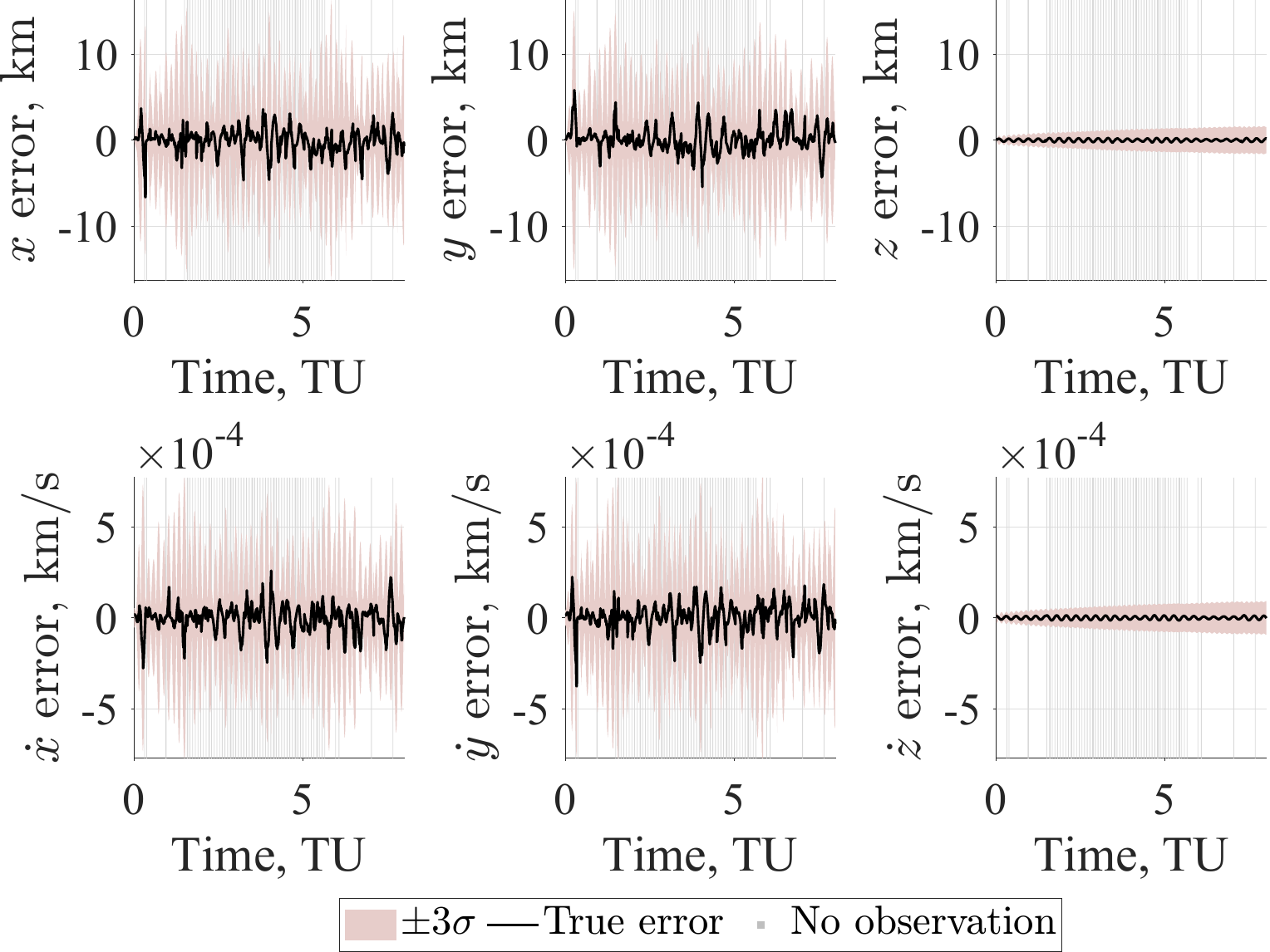}
        \caption{STP-B best case, an LPWO.}
	\label{fig:3sigma_TPB_best}
     \end{subfigure}
    \hspace{1em}
     \begin{subfigure}[b]{0.37\textwidth}
         \centering
            \includegraphics[width=\linewidth]{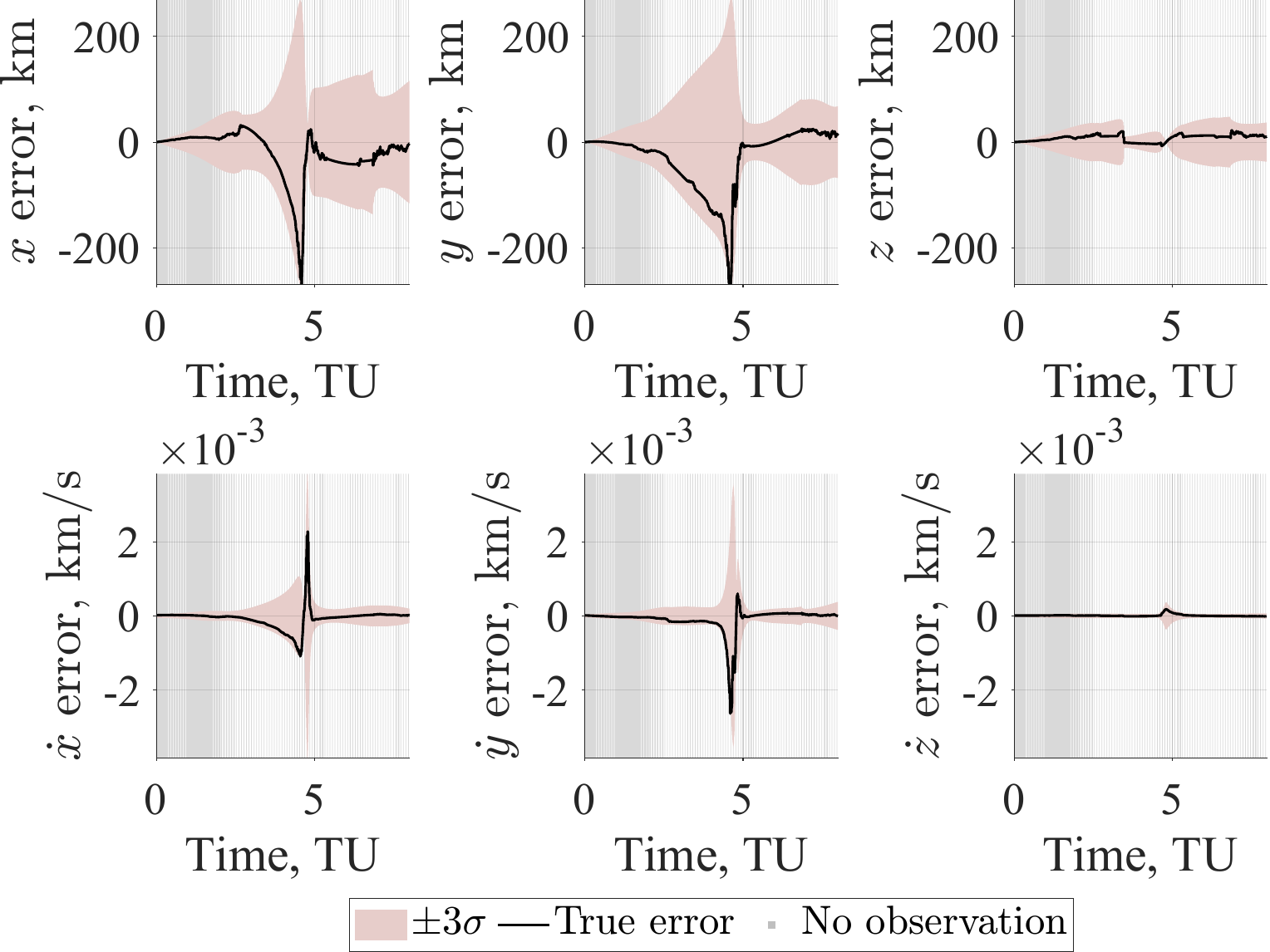}
            \caption{STP-B worst case, an R1:1O.}
	   \label{fig:3sigma_TPB_worst}
     \end{subfigure}
    \hspace{1em}
     \begin{subfigure}[b]{0.37\textwidth}
            \centering
            \includegraphics[width=\linewidth]{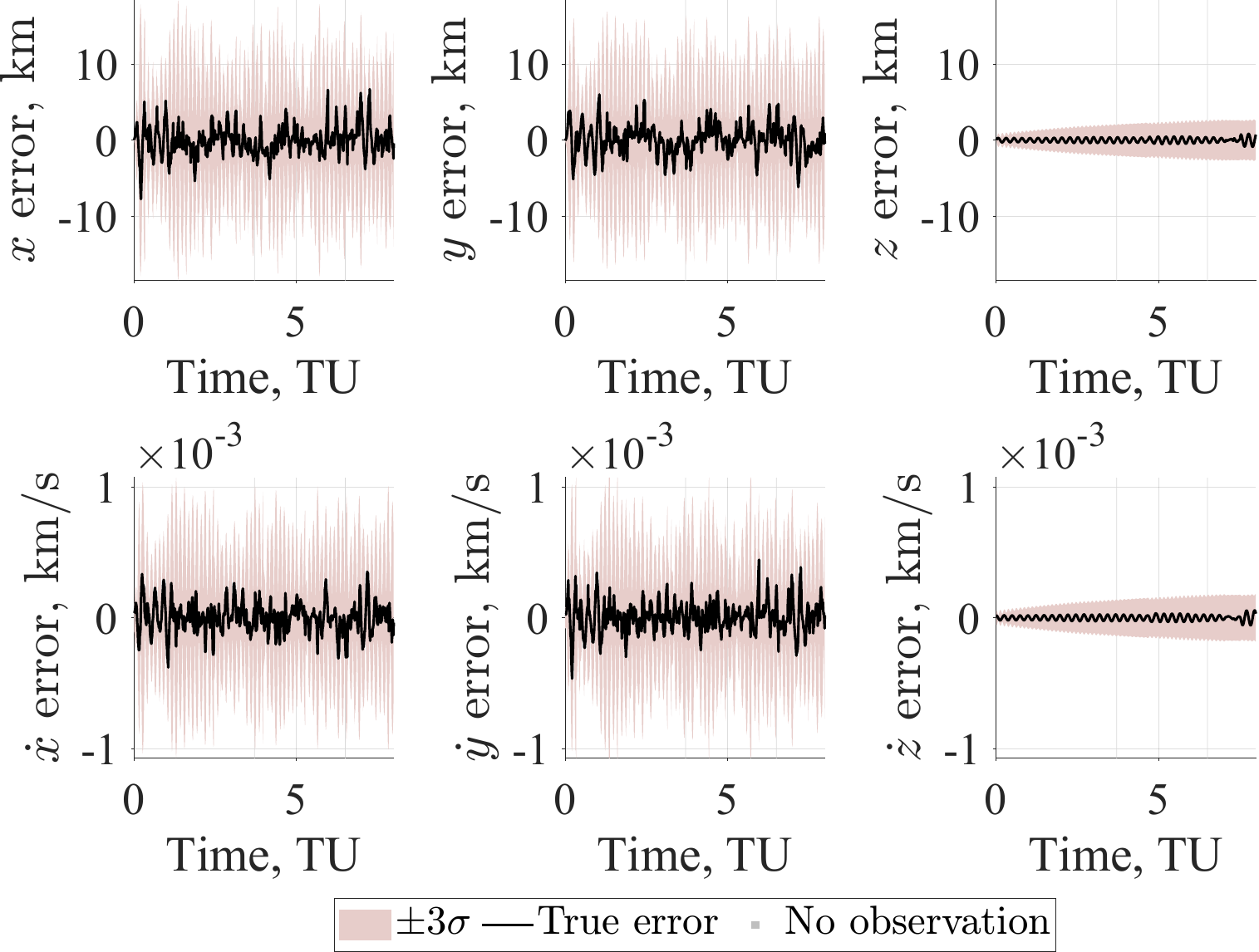}
        \caption{STP-C best case, an LPWO.}
	\label{fig:3sigma_TPC_best}
     \end{subfigure}
    \hspace{1em}
     \begin{subfigure}[b]{0.37\textwidth}
         \centering
            \includegraphics[width=\linewidth]{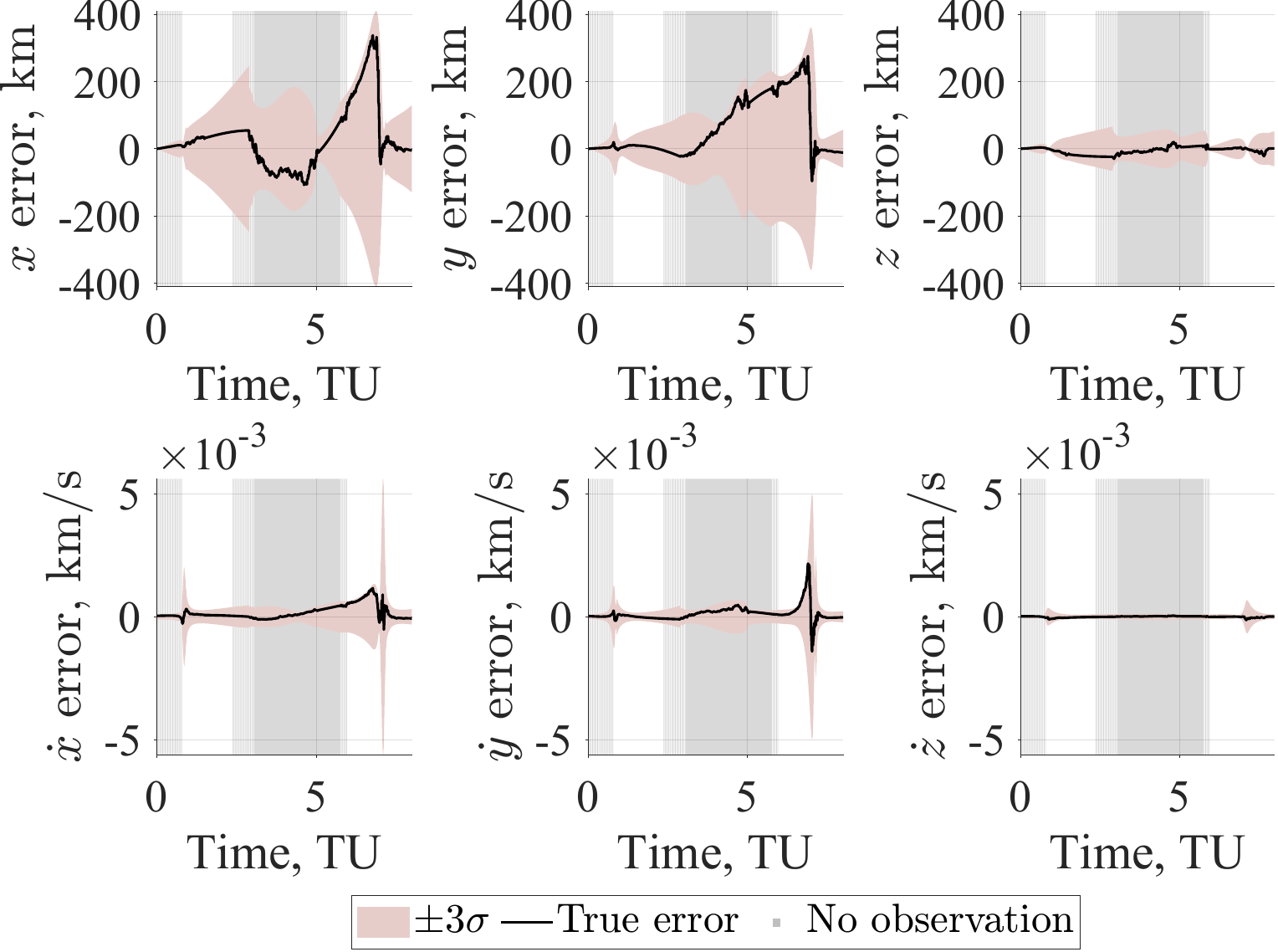}
        \caption{STP-C worst case, an R1:1O.}
	\label{fig:3sigma_TPC_worst}
     \end{subfigure}
        \caption{3$\sigma$ error plots of each case's best and worst performance, for the low-fidelity sensors.}
        \label{fig:3sigma_of}
\end{figure}

To better visualize the performance against each of these targets, Fig.~\ref{fig:3sigma_of} shows the error between the ground truth as well as the three-sigma envelopes. In these figures, the gray shading represents a lack of target visibility by the tasked observer(s). In all three STPs, the worst-performing cases have more missed observations than the best-performing cases.

\begin{figure}
    \centering

\begin{subfigure}[htbp]{0.6\textwidth}
        \centering
        \includegraphics[width=\linewidth]{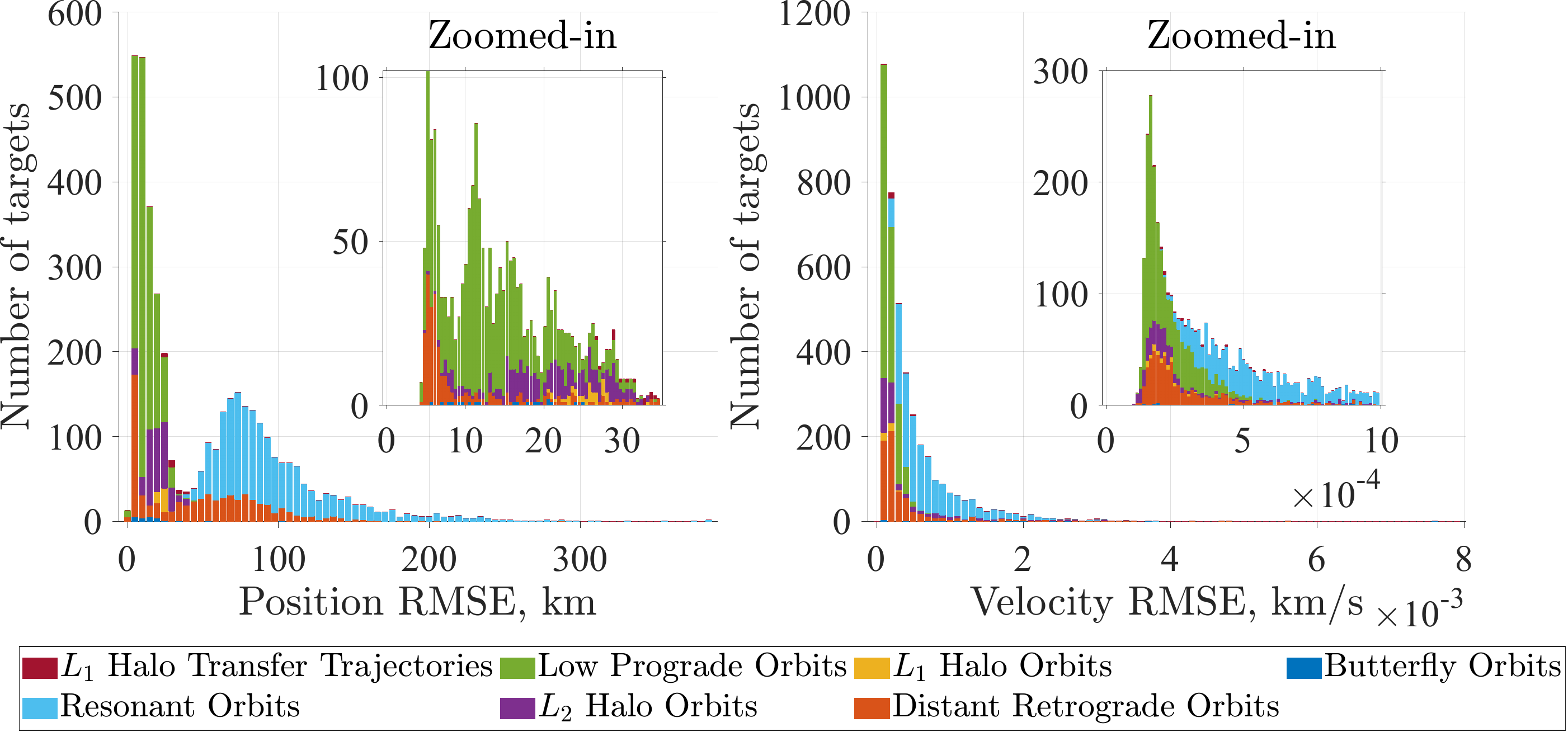}
        \caption{Baseline orbit.}
	    \label{fig:SO_hist}
\end{subfigure}

\begin{subfigure}[htbp]{0.6\textwidth}
        \centering
        \includegraphics[width=\linewidth]{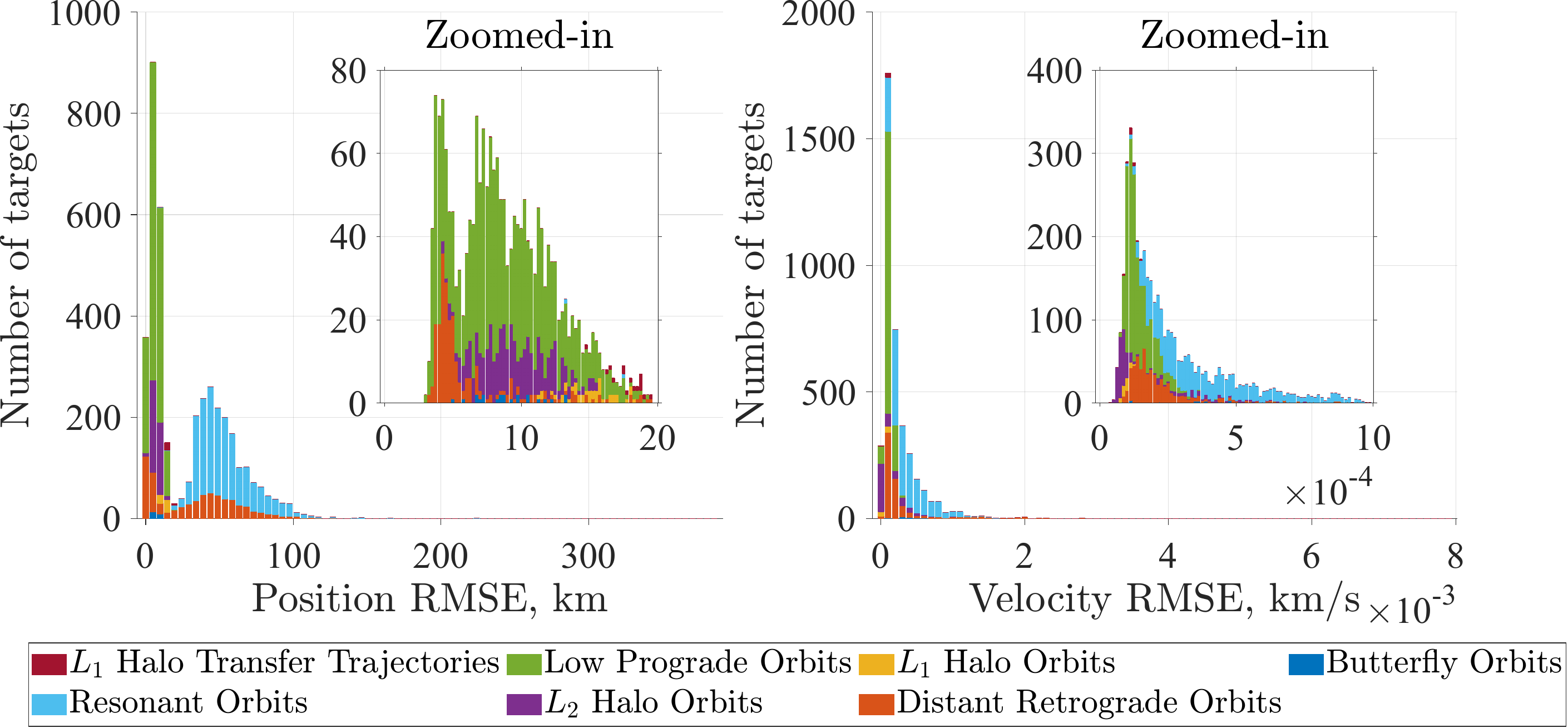}
        \caption{STP-A.}
	  \label{fig:TPA_hist}
\end{subfigure}

\begin{subfigure}[htbp]{0.6\textwidth}
        \centering
        \includegraphics[width=\linewidth]{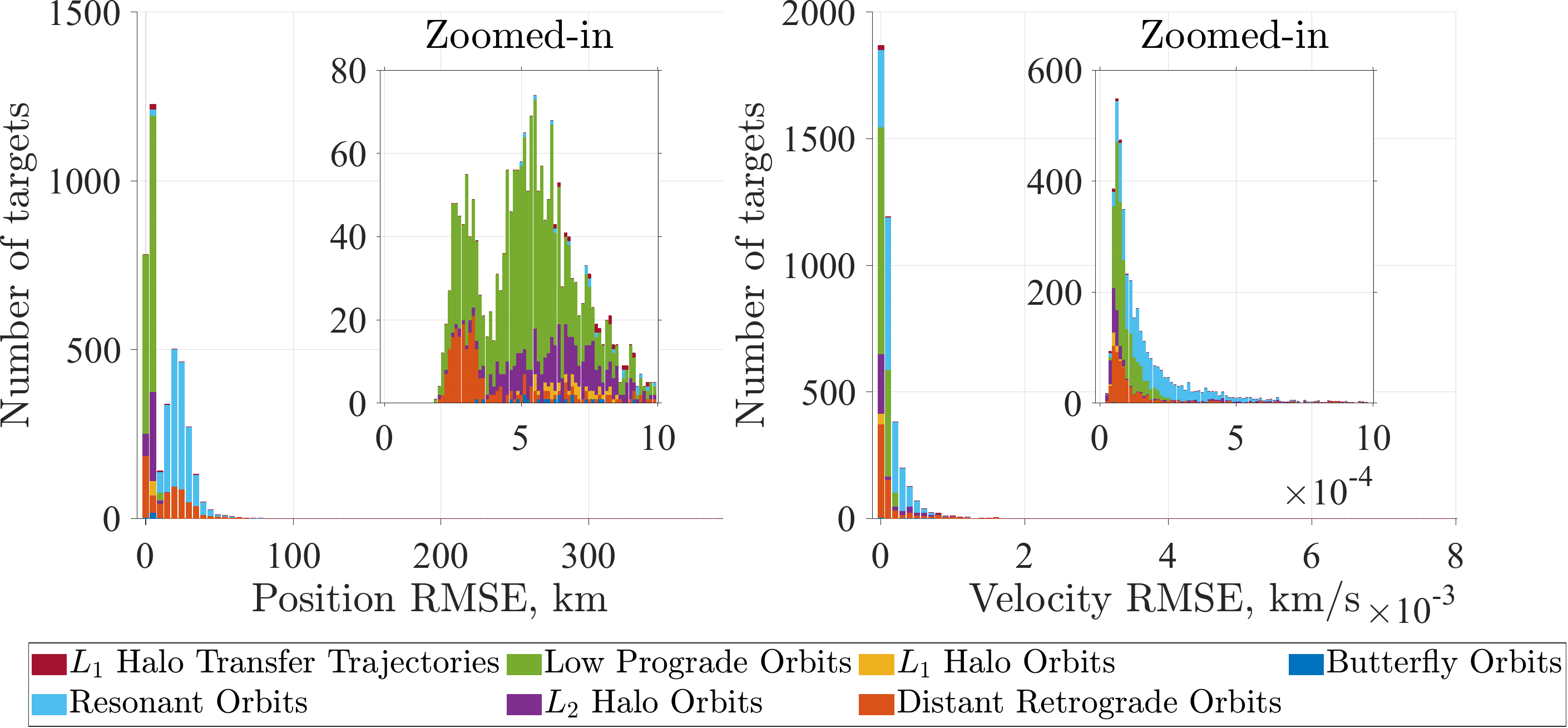}
        \caption{STP-B.}
	    \label{fig:TPB_hist}
\end{subfigure}

\begin{subfigure}[htbp]{0.6\textwidth}
         \centering
        \includegraphics[width=\linewidth]{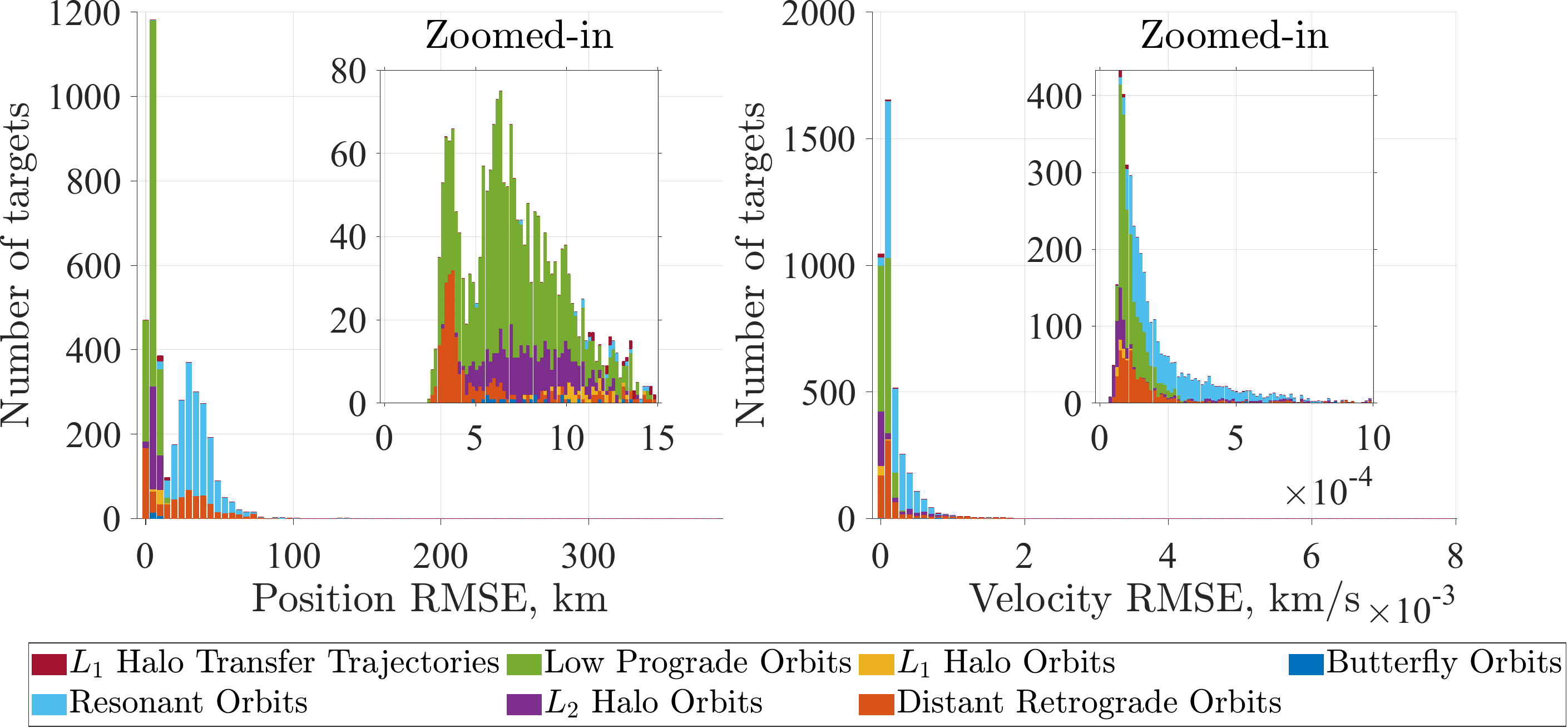}
        \caption{STP-C.}
	   \label{fig:TPC_hist}
\end{subfigure}

    \caption{Histograms showing the distribution of targets separated into similar orbit groupings, for the low-fidelity observers.}
    \label{fig:lowfi_histograms}
\end{figure}

While the box plots showing the statistical distribution across each family can be seen in Fig.~\ref{fig:box_plots}, it is helpful to look at the distribution of target RMSEs across the validation target set. Histograms dividing the orbit families into similar groups (for example, the L1NHO and L1SHO are combined) show this distribution in Fig.~\ref{fig:lowfi_histograms}. The velocity RMSEs are all skewed to the left, indicating that the constellation performs well. However, in the position RMSE histograms there is a secondary grouping of targets (resonant orbits and the DRO families) with higher RMSEs than the majority of targets. These orbits have longer periods, and therefore travel a farther distance than the rest of the orbit families.

A trend that can be seen in Figs.~\ref{fig:box_plots}~and~\ref{fig:lowfi_histograms} is that the range of RMSE values for the DRO family is larger than that of the other families. The DRO family is not homogeneous, as seen in the period length of the orbits (demonstrated in how far the orbit traverses cislunar space). The difference between the longest and shortest period in the DRO family is more than \SI{6}{TU}, with the LPWO family having the next highest difference, at only about \SI{2.4}{TU}. The other eight orbit families that are not resonant orbits or DROs have zero orbits with periods greater than \SI{3.75}{TU}.

To investigate further, we split the DRO family into two groups, the short-period DROs (having a period less than \SI{3.75}{TU}) and the long-period DROs (having a period greater than \SI{3.75}{TU}). The short-period DROs have average position RMSEs of only \SI{8.70}{km}, \SI{6.40}{km}, and \SI{6.27}{km} for STP-A, STP-B, and STP-C, respectively. This is closer to the best-performing orbit families. The long-period DROs have average position RMSEs of \SI{55.29}{km}, \SI{27.15}{km}, and \SI{38.71}{km} for STP-A, STP-B, and STP-C, respectively, which is closer to the results of the worst-performing families.

\begin{figure}[htbp]
        \centering
        \includegraphics[width=0.8\linewidth]{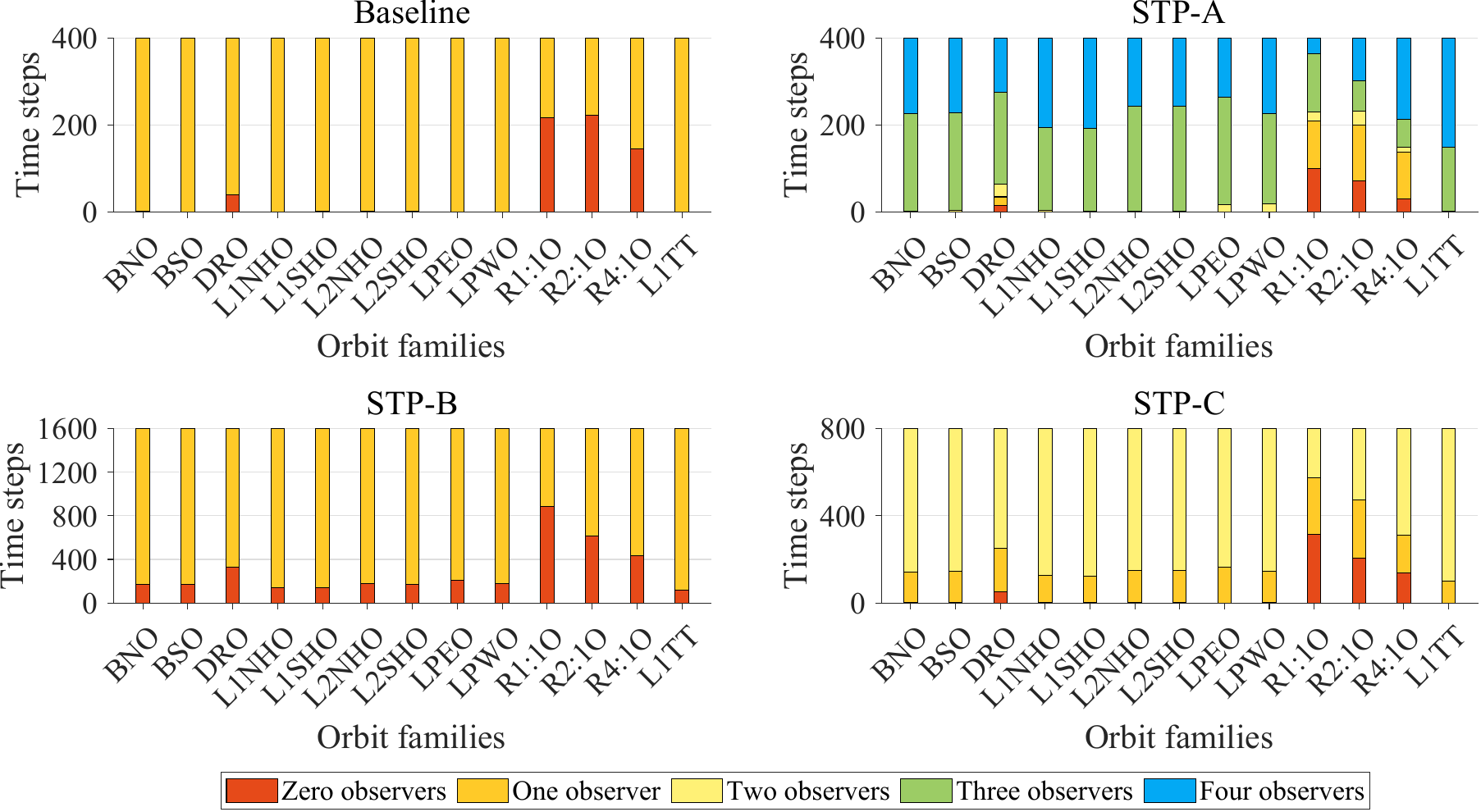}
        \caption{Bar charts showing the average visibility of targets by orbit family, for low-fidelity observers.}
        \label{fig:const_vis}
\end{figure}

As the worst-performing targets in Fig.~\ref{fig:3sigma_of} are all in the R1:1O family, another family with long periods, one possibility deserving investigation for the poorer performance against these targets in comparison to other orbit families is the visibility of the targets. Figure~\ref{fig:const_vis} shows the average visibility of targets by family for all four cases. The orbits in the R1:1O, R2:1O, R4:1O, and DRO  families have a lower average visibility than the other orbit families in all sensor tasking procedures.

    Testing the baseline and three sensor tasking procedures against a wide variety of targets demonstrates a few trends that should be noted. First, the worst-performing families in the experiments have periods much longer than those of the best-performing targets. This means that they travel further distances in cislunar space, which leads to poorer visibility on the part of the constellation. This poorer visibility means that the EKF has fewer steps to update the state estimate. One solution to this issue that can be easily tested is raising the fidelity of the sensors, which will be tested in Sec.~\ref{sec:parameters}.

\subsection{Comparative Analysis with High-Fidelity Sensors}
\label{sec:parameters}
\begin{figure}[htb]
     \centering
     \includegraphics[width=\textwidth]{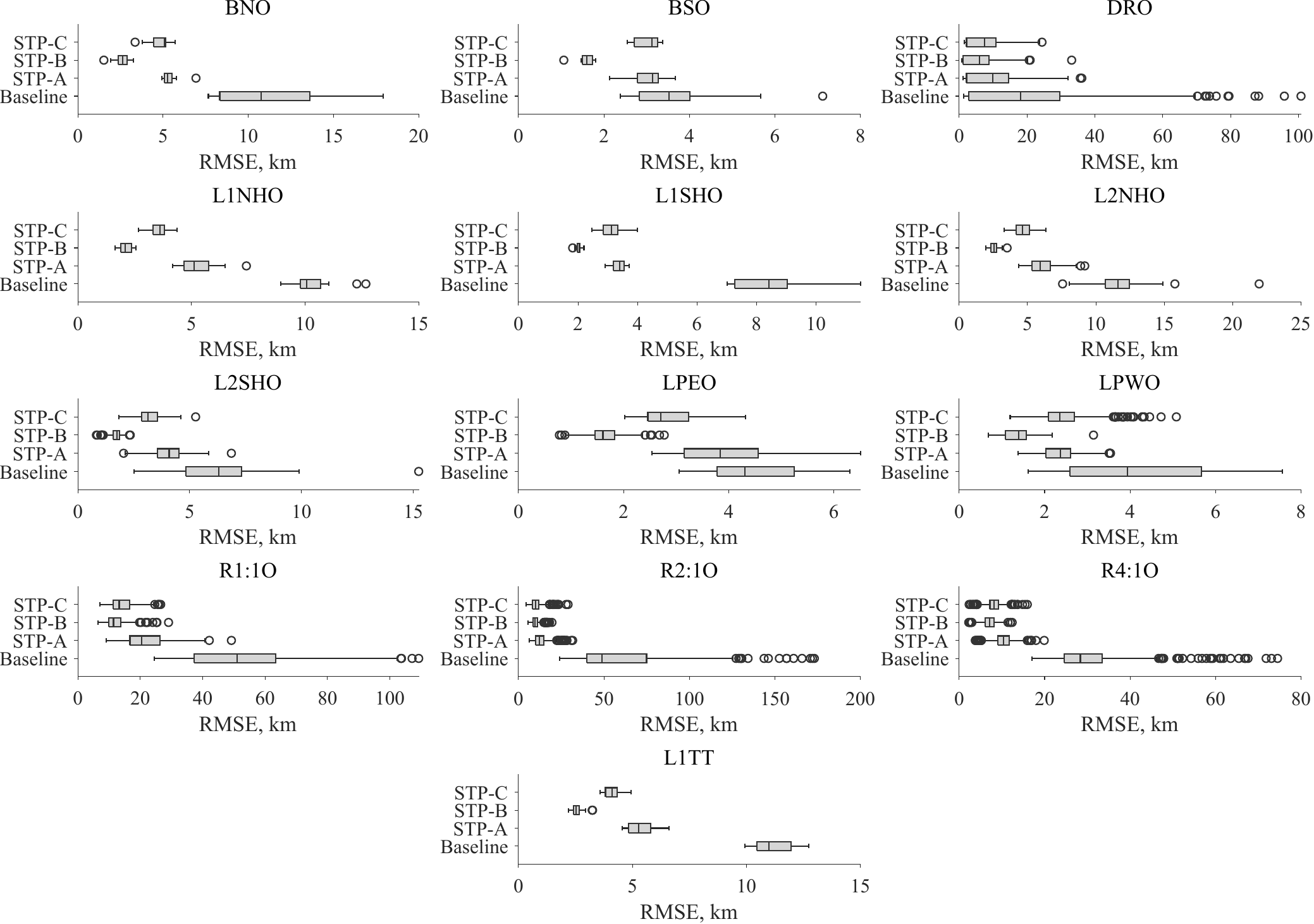}
     \caption{Box plots for each family and experimental case, with high-fidelity sensors.}
     \label{fig:hifi_box_plots}
\end{figure}

The previous experiments only use the low-fidelity sensor uncertainty value from Ref.~\cite{block2022cislunar}. A comparative analysis is presented here to discuss how the optimization changes when high-fidelity sensors are used, and if the RMSE results can be lowered across the whole set of targets over Sec.~\ref{sec:ex_res}. The new uncertainty of the sensors is changed to \SI{26.7518}{arcseconds} \cite{block2022cislunar}. All other parameters are the same, meaning the optimization and validation target set are identical to the previous experiments. The GA returns loss function values of \SI{6.59}{km}, \SI{3.43}{km}, and \SI{5.07}{km} for STP-A, STP-B, and STP-C, respectively. For all three of these sensor tasking procedures, these values are about a $66\%$ decrease over the values from Sec.~\ref{sec:ex_res}. Figure~\ref{fig:hifi_box_plots} shows box plots showing the median, quartiles, maximum, and minimum RMSE values for all four cases across the 13 orbit families of the targets. As in the low-fidelity experiment, STP-B performs the best across all the orbit families, while STP-A is the worst procedure. The same conclusion can be drawn here; that it could be ideal to make measurements as often as possible, no matter how many observers are taking measurements.

\begin{figure}[htbp]
    \centering
    \begin{subfigure}[b]{0.75\textwidth}
        \centering
        \includegraphics[width=\linewidth]{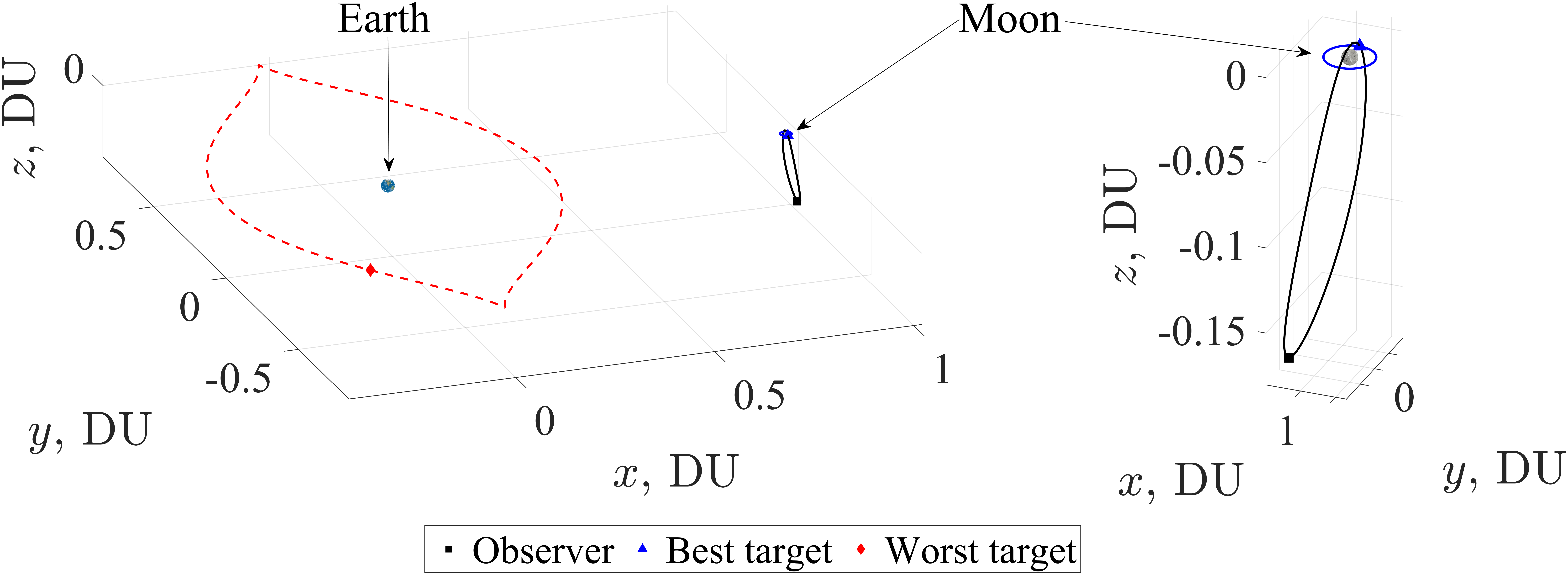}
        \caption{Baseline orbit.}
	    \label{fig:hifi_single_obs}
    \end{subfigure}
    \begin{subfigure}[b]{0.75\textwidth}
        \centering
        \includegraphics[width=\linewidth]{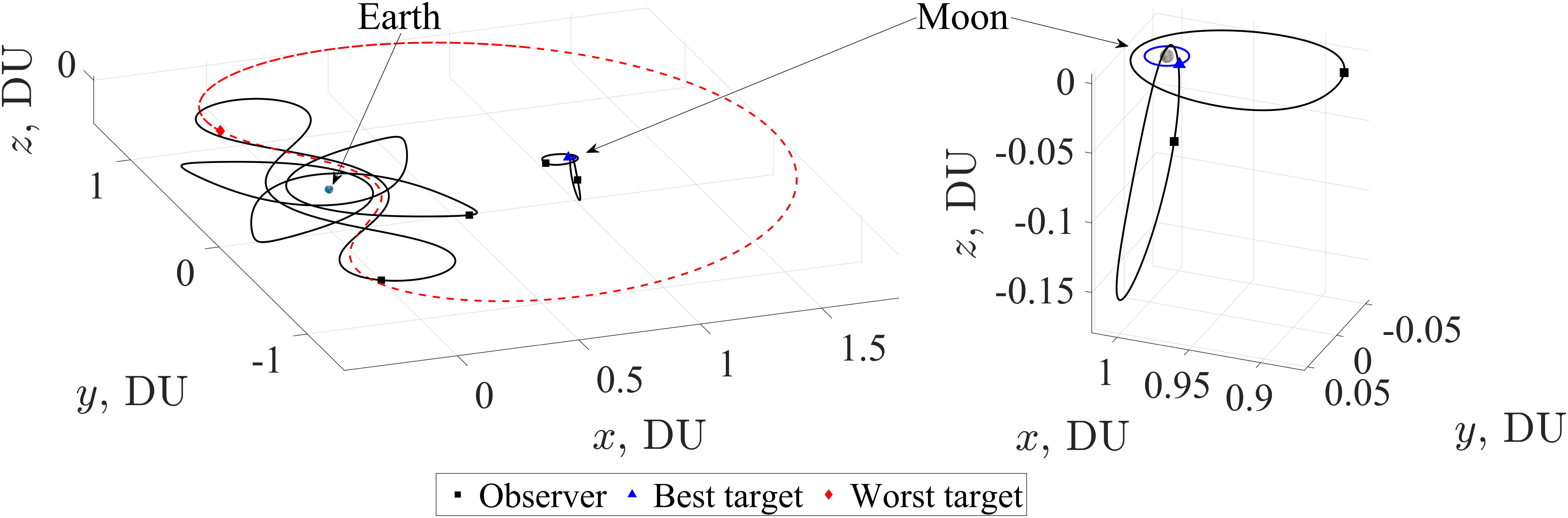}
        \caption{STP-A.}
	  \label{fig:hifi_TPA_constellation}
    \end{subfigure}
    \begin{subfigure}[b]{0.75\textwidth}
        \centering
        \includegraphics[width=\linewidth]{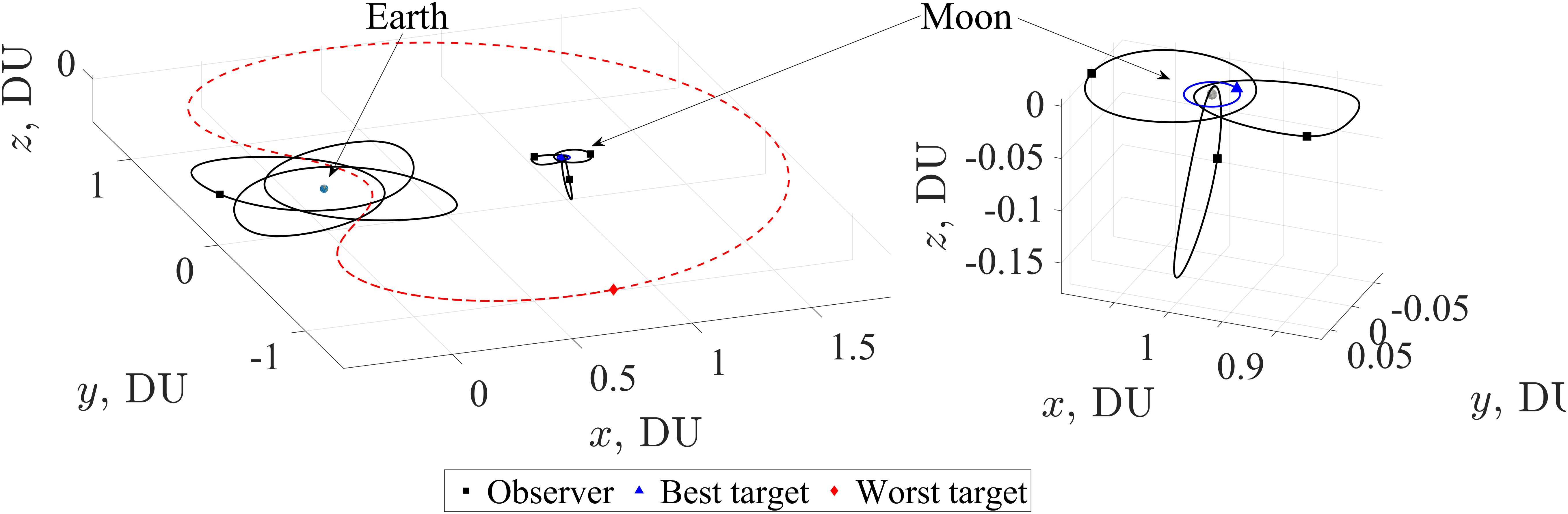}
        \caption{STP-B.}
	    \label{fig:hifi_TPB_constellation}
    \end{subfigure}
     \begin{subfigure}[b]{0.75\textwidth}
         \centering
        \includegraphics[width=\linewidth]{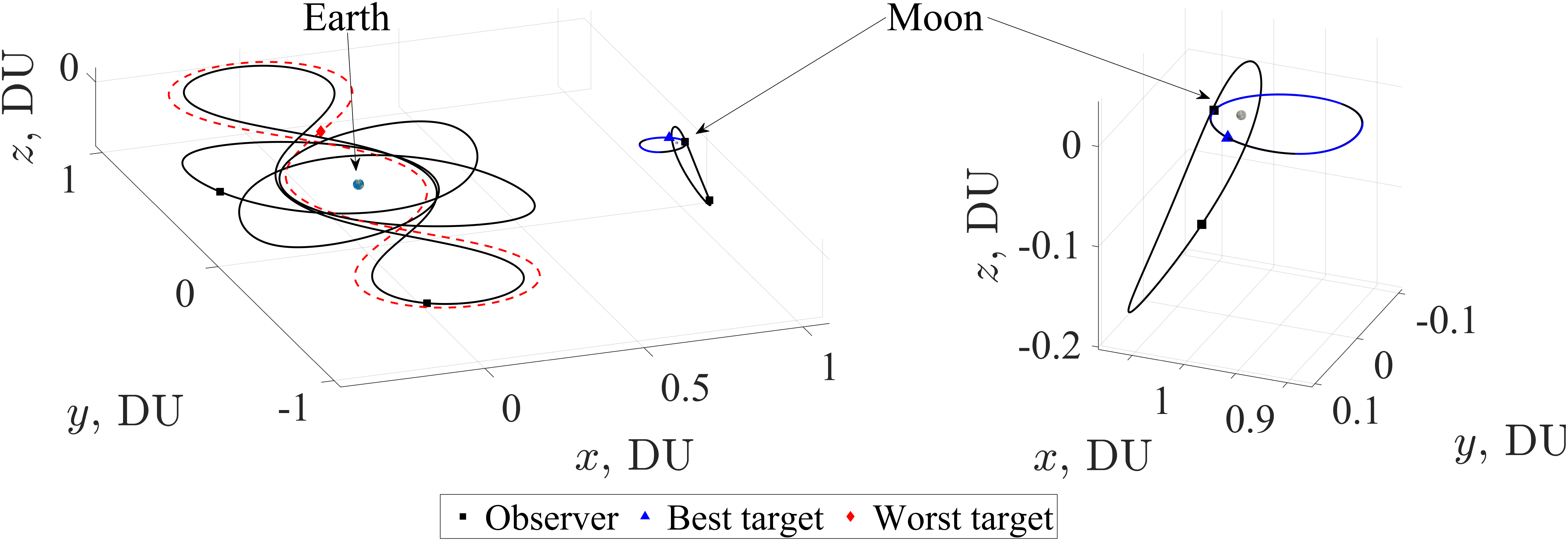}
        \caption{STP-C.}
	\label{fig:hifi_TPC_constellation}
     \end{subfigure}
        \caption{Observer orbits for each high-fidelity case, along with worst- and best-performing targets, and a close-up of near the Moon (right), for high-fidelity observers.}
        \label{fig:hifi_optimal_constellations}
\end{figure}

The initial conditions of the constituent observers in each of these constellations can be seen in Table~\ref{tab:opt_ICs_hifi}, within Appendix~C. STP-A and STP-B share an observer in an L2SHO orbit, while STP-B and STP-C share an observer in an R4:1O orbit. In both of these cases, the third constellation has an observer in that orbit family, albeit in a different orbit (STP-C has an L2SHO, and STP-A has an R4:1O). Each constellation has observers in different orbits of the LPWO family. These orbits can be seen in Fig.~\ref{fig:hifi_optimal_constellations}. STP-A and STP-C having two observers in resonant orbits is a new phenomenon not seen in the previous experiment.

\begin{figure}[htbp]
     \centering
     \begin{subfigure}[b]{0.37\textwidth}
            \centering
            \includegraphics[width=\linewidth]{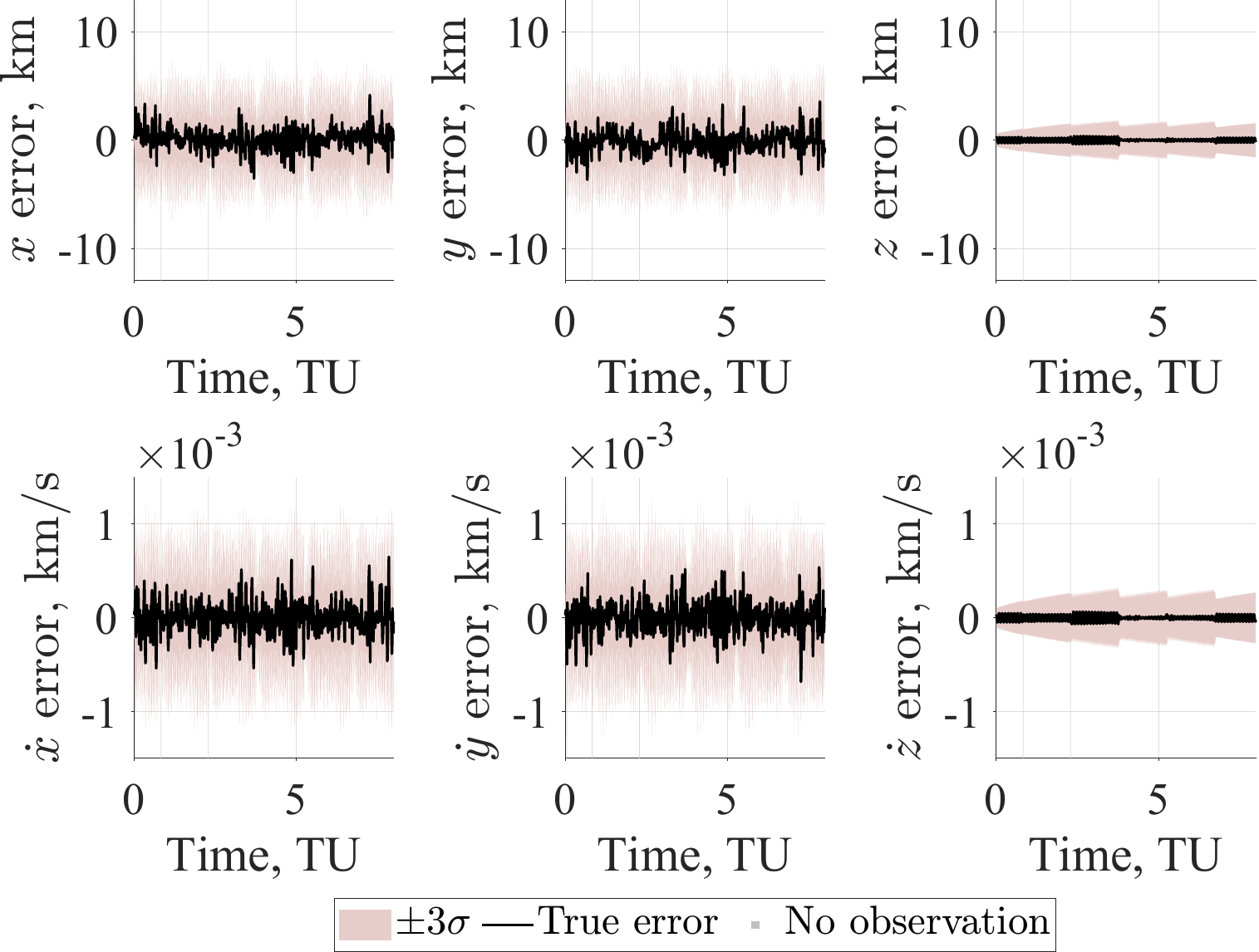}
            \caption{Baseline best case, a DRO.}
	    \label{fig:hifi_3sigma_SO_best}
     \end{subfigure}
    \hspace{1em}
     \begin{subfigure}[b]{0.37\textwidth}
         \centering
        \includegraphics[width=\linewidth]{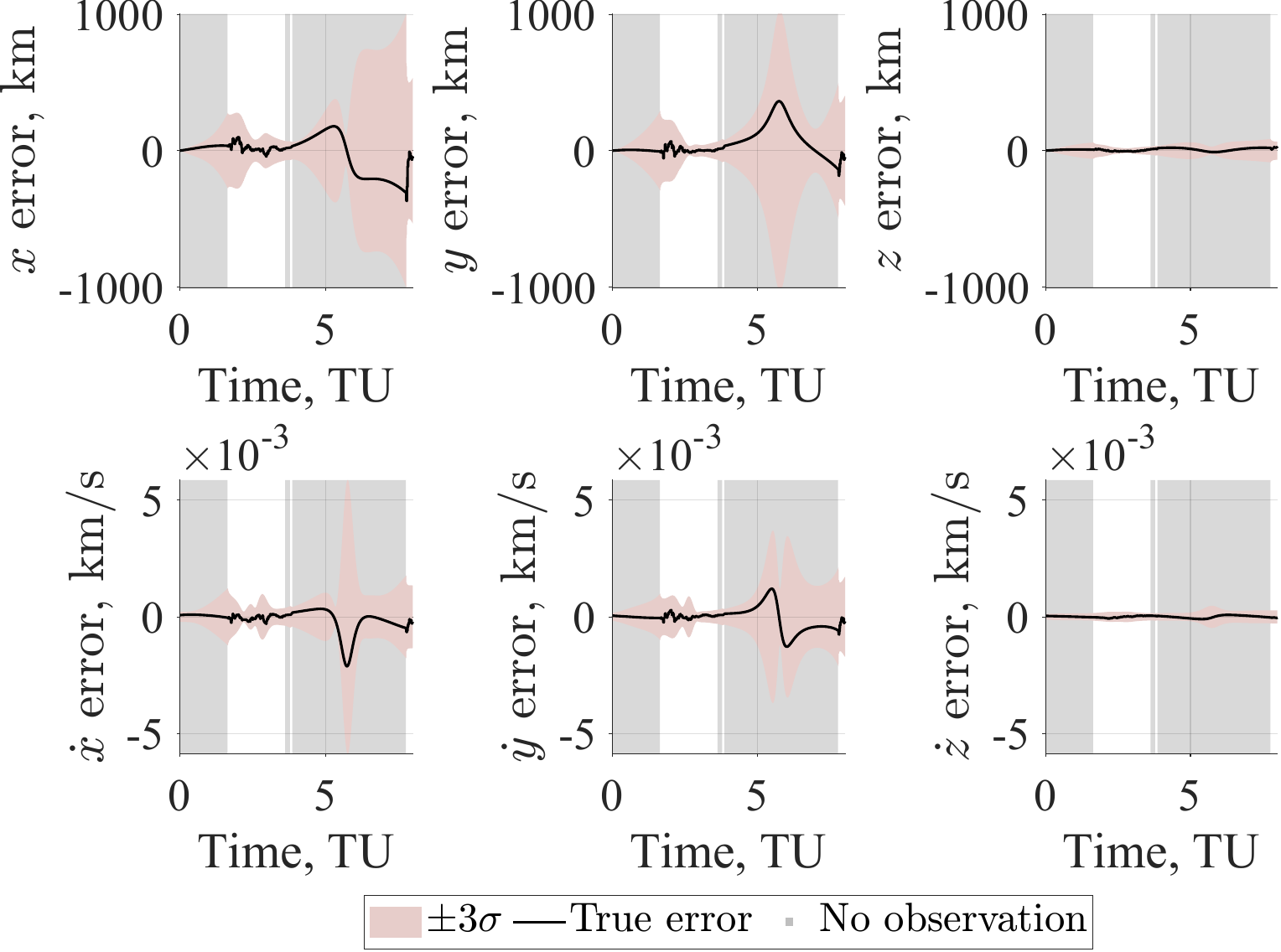}
        \caption{Baseline worst case, an R2:1O.}
	\label{fig:hifi_3sigma_SO_worst}
     \end{subfigure}
     \hspace{1em}
     \begin{subfigure}[b]{0.37\textwidth}
         \centering
         \includegraphics[width=\linewidth]{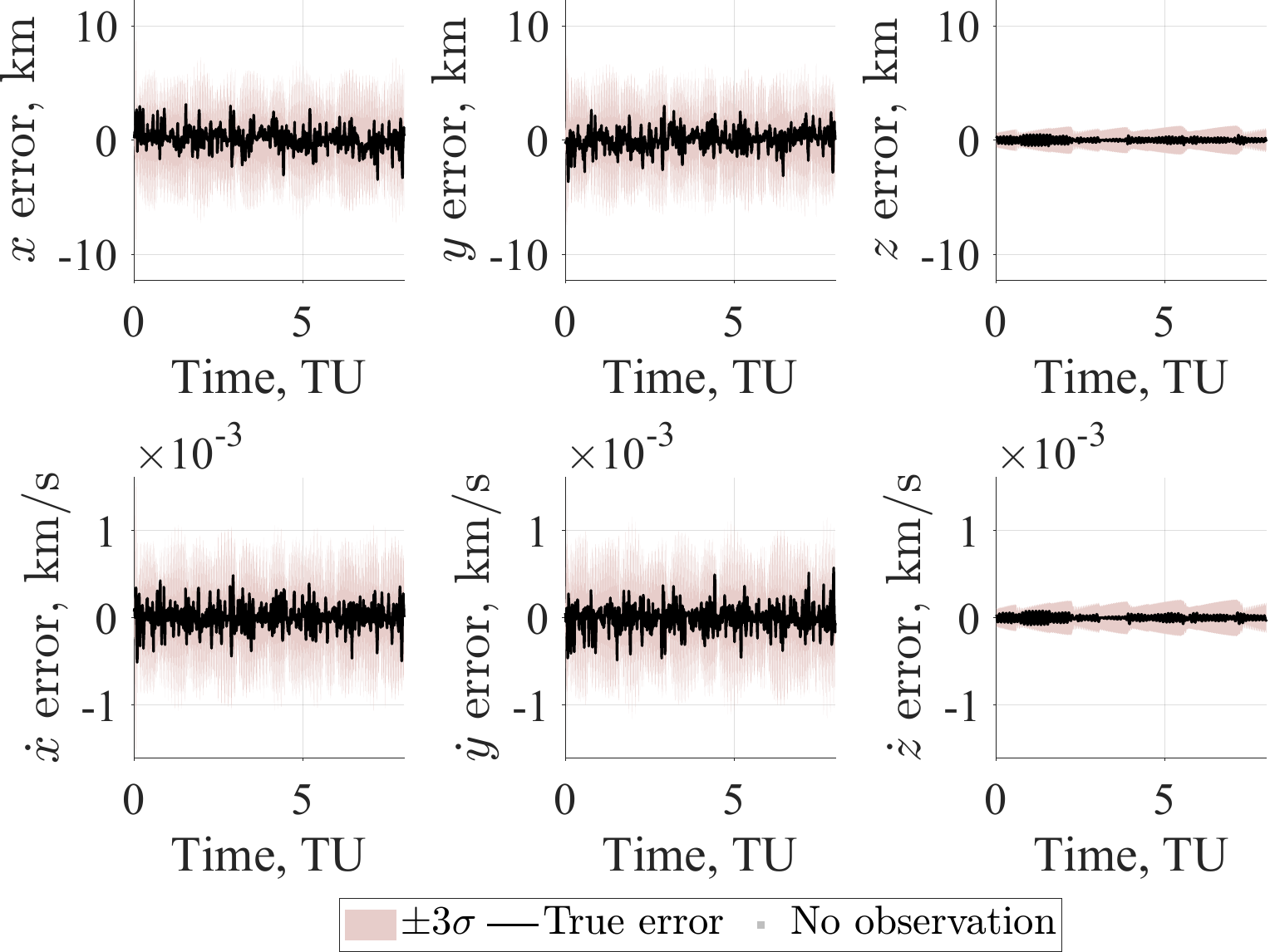}
            \caption{STP-A best case, a DRO.}
	   \label{fig:hifi_3sigma_TPA_best}
     \end{subfigure}
    \hspace{1em}
     \begin{subfigure}[b]{0.37\textwidth}
         \centering
         \includegraphics[width=\linewidth]{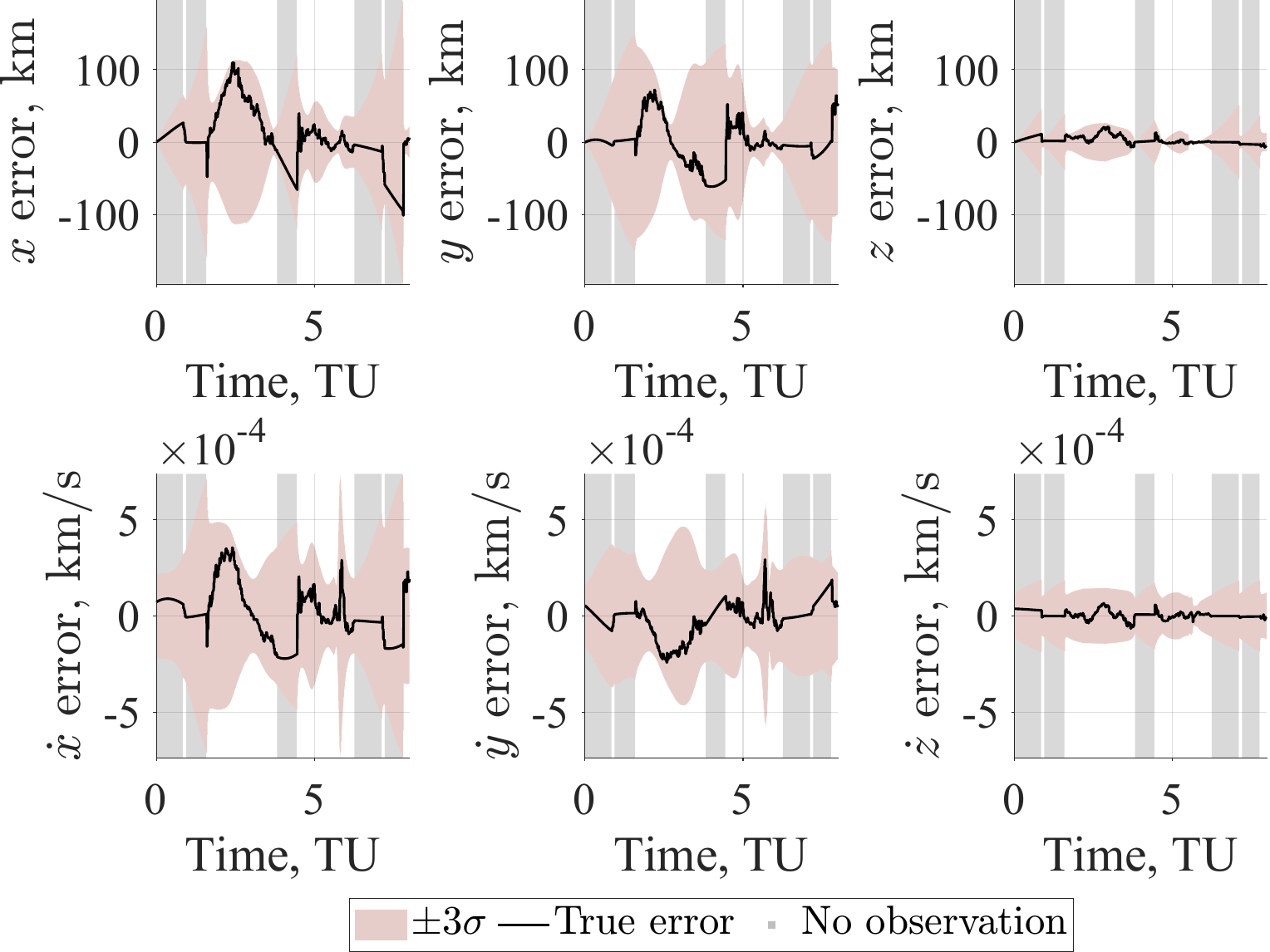}
        \caption{STP-A worst case, an R1:1O.}
	\label{fig:hifi_3sigma_TPA_worst}
     \end{subfigure}
     \hspace{1em}
     \begin{subfigure}[b]{0.37\textwidth}
            \centering
            \includegraphics[width=\linewidth]{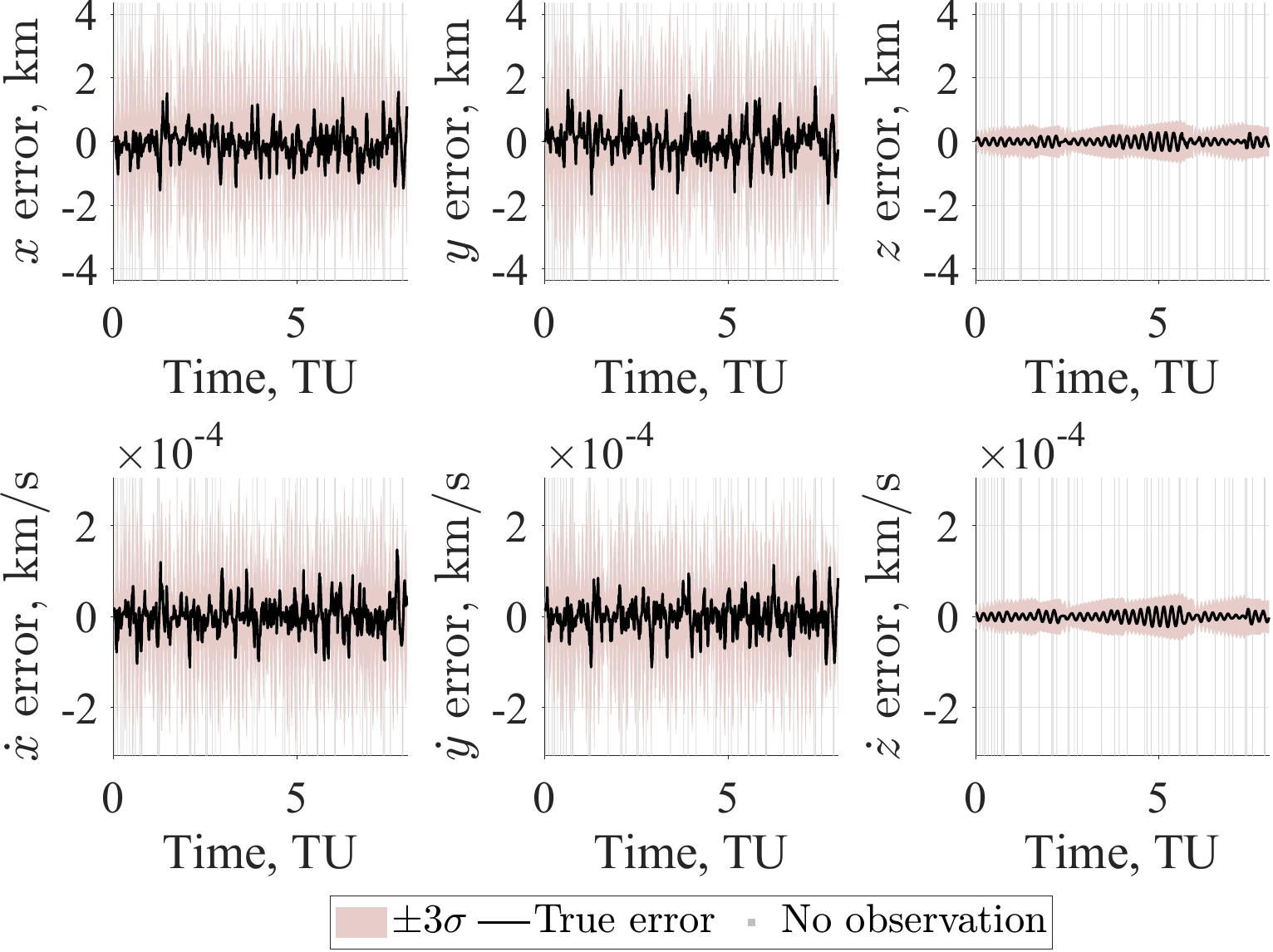}
        \caption{STP-B best case, an LPWO.}
	\label{fig:hifi_3sigma_TPB_best}
     \end{subfigure}
    \hspace{1em}
     \begin{subfigure}[b]{0.37\textwidth}
         \centering
            \includegraphics[width=\linewidth]{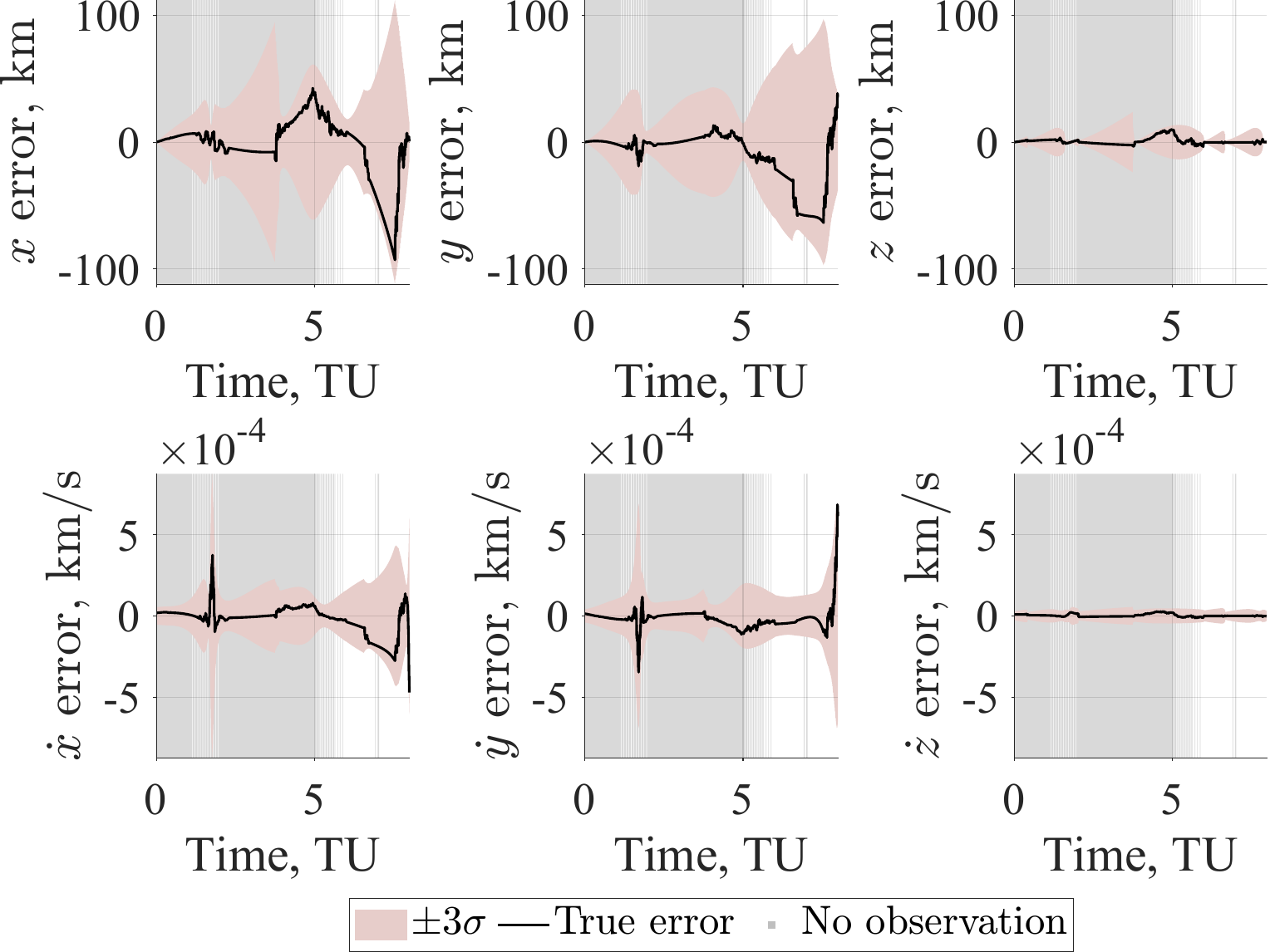}
            \caption{STP-B worst case, a DRO.}
	   \label{fig:hifi_3sigma_TPB_worst}
     \end{subfigure}
     \hspace{1em}
     \begin{subfigure}[b]{0.37\textwidth}
            \centering
            \includegraphics[width=\linewidth]{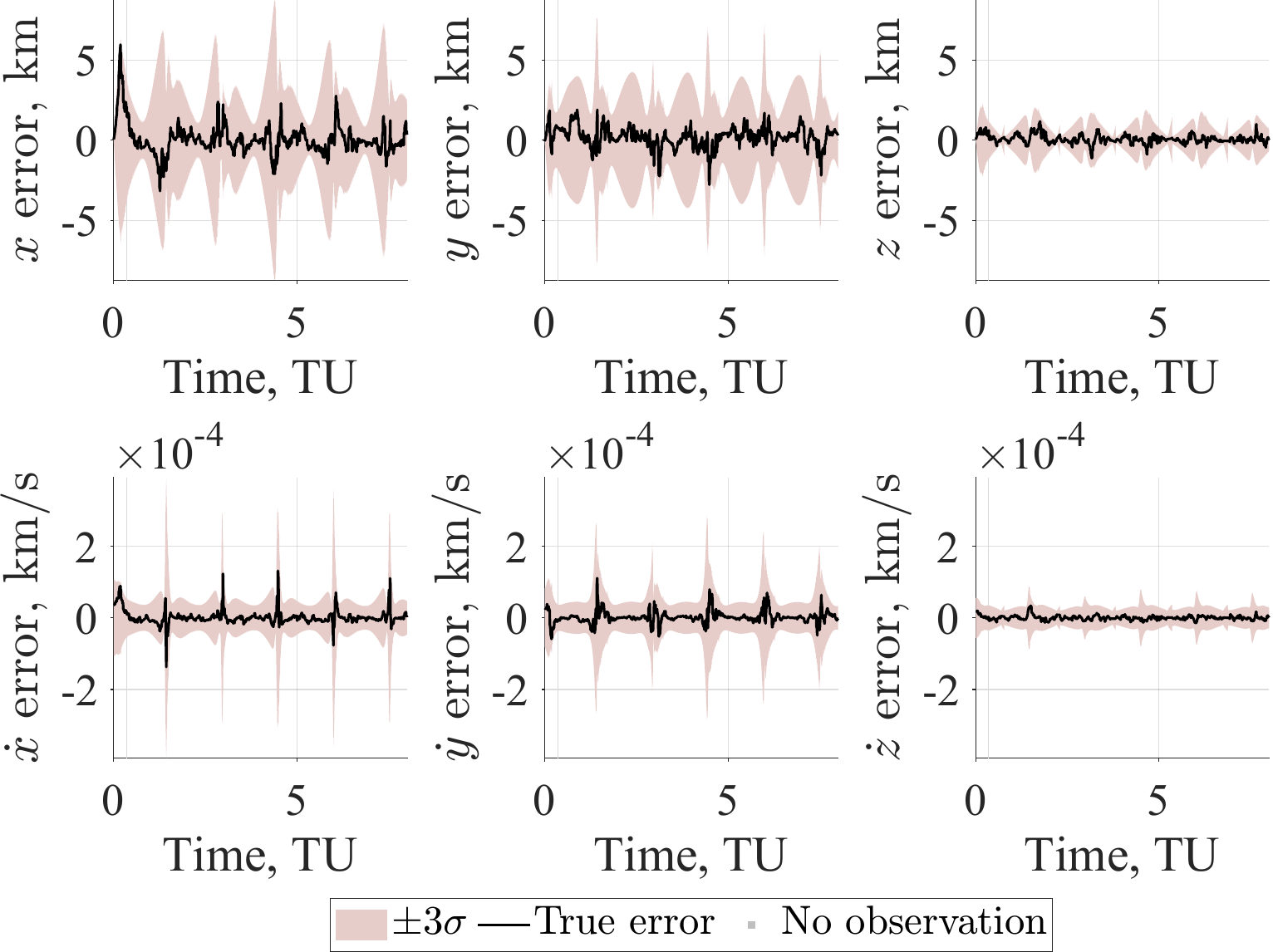}
        \caption{STP-C best case, an LPWO.}
	\label{fig:hifi_3sigma_TPC_best}
     \end{subfigure}
    \hspace{1em}
     \begin{subfigure}[b]{0.37\textwidth}
         \centering
            \includegraphics[width=\linewidth]{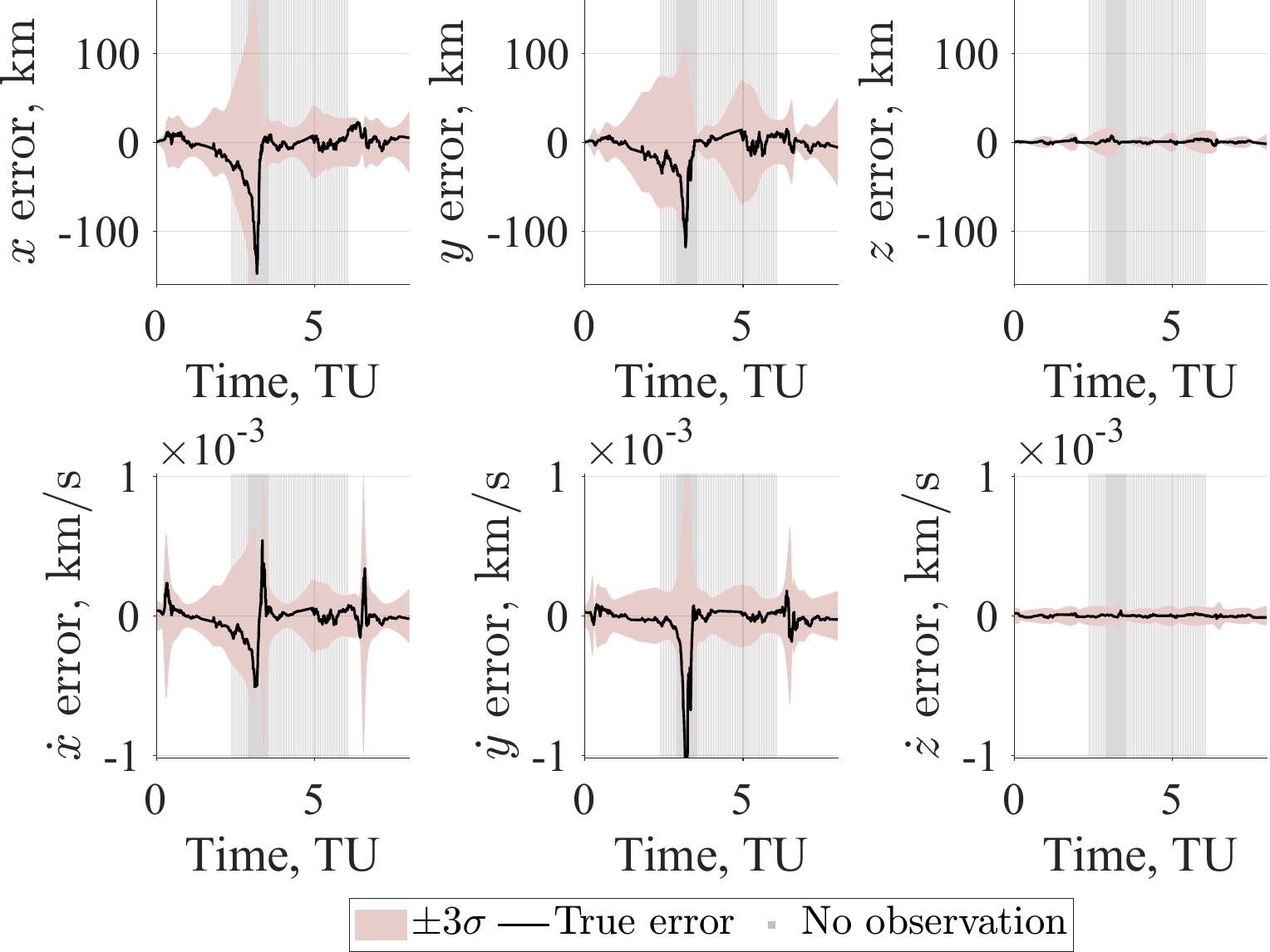}
        \caption{STP-C worst case, an R2:1O.}
	\label{fig:hifi_3sigma_TPC_worst}
     \end{subfigure}
    \caption{3$\sigma$ error plots of each case's best and worst performance, for the high-fidelity sensors.}
        \label{fig:hifi_3sigma_of}
\end{figure}

As in the low-fidelity experiments, Fig.~\ref{fig:hifi_optimal_constellations} shows the targets each case performs the best and worst against. For the baseline, STP-A, STP-B, and STP-C, the worst-performing are in the R2:1O, R1:1O, DRO, and R2:1O families, respectively. The best performances came against DROs for the baseline (an RMSE of \SI{1.34}{km}) and STP-A (\SI{1.23}{km}), while STP-B (\SI{0.68}{km}) and STP-C (\SI{1.19}{km}) perform best against LPWOs. The initial conditions of these targets can be seen in Table~\ref{tab:targ_ics_hifi}, within Appendix~C. Their true error and three-sigma envelope can be seen in Fig.~\ref{fig:hifi_3sigma_of}, and again, the best-performing targets have a much higher visibility than the worst-performing cases.

\begin{figure}[htbp]
    \centering

    \begin{subfigure}[htbp]{0.55\textwidth}
        \centering
        \includegraphics[width=\linewidth]{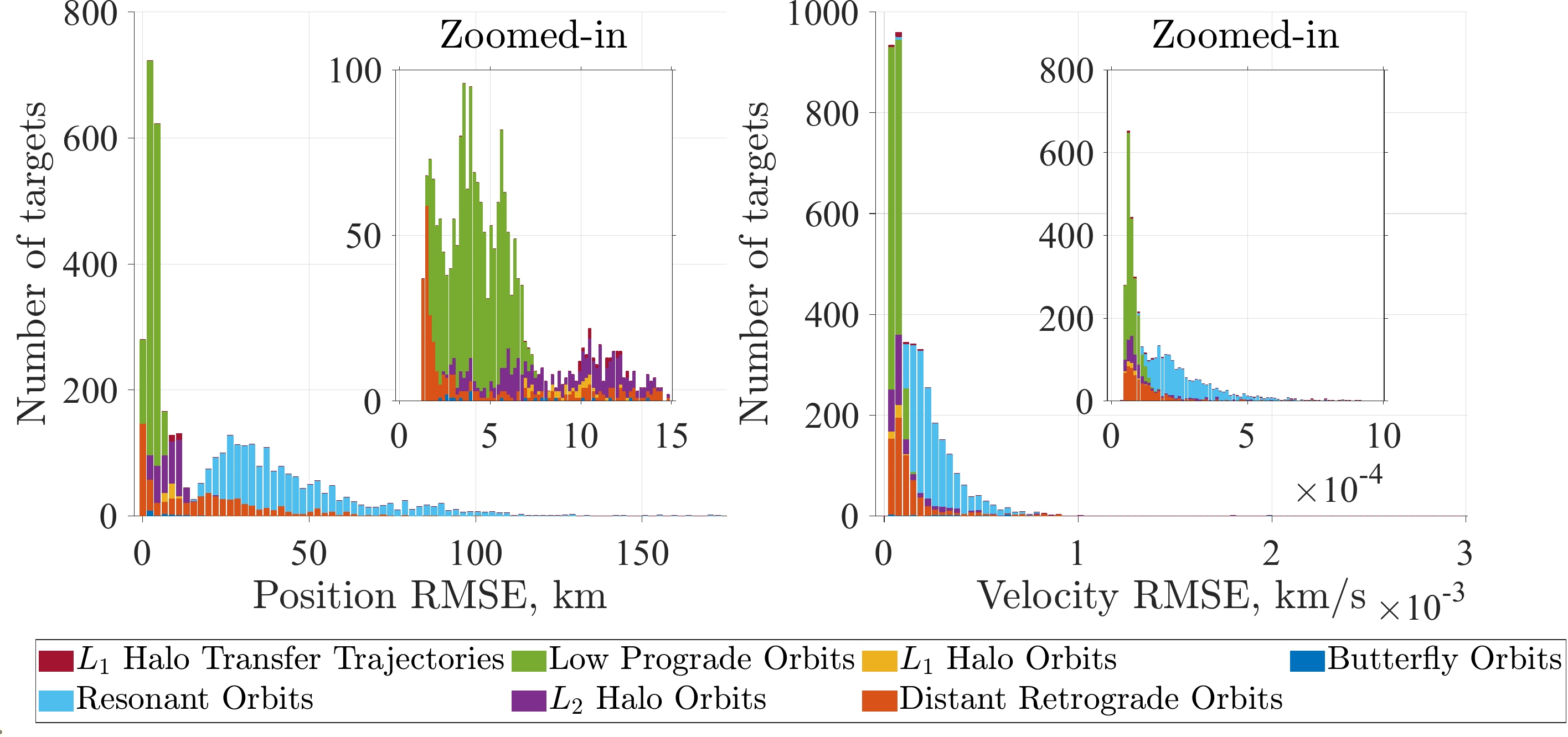}
        \caption{Baseline orbit.}
	    \label{fig:HF_SO_H}
    \end{subfigure}

    \begin{subfigure}[htbp]{0.55\textwidth}
        \centering
        \includegraphics[width=\linewidth]{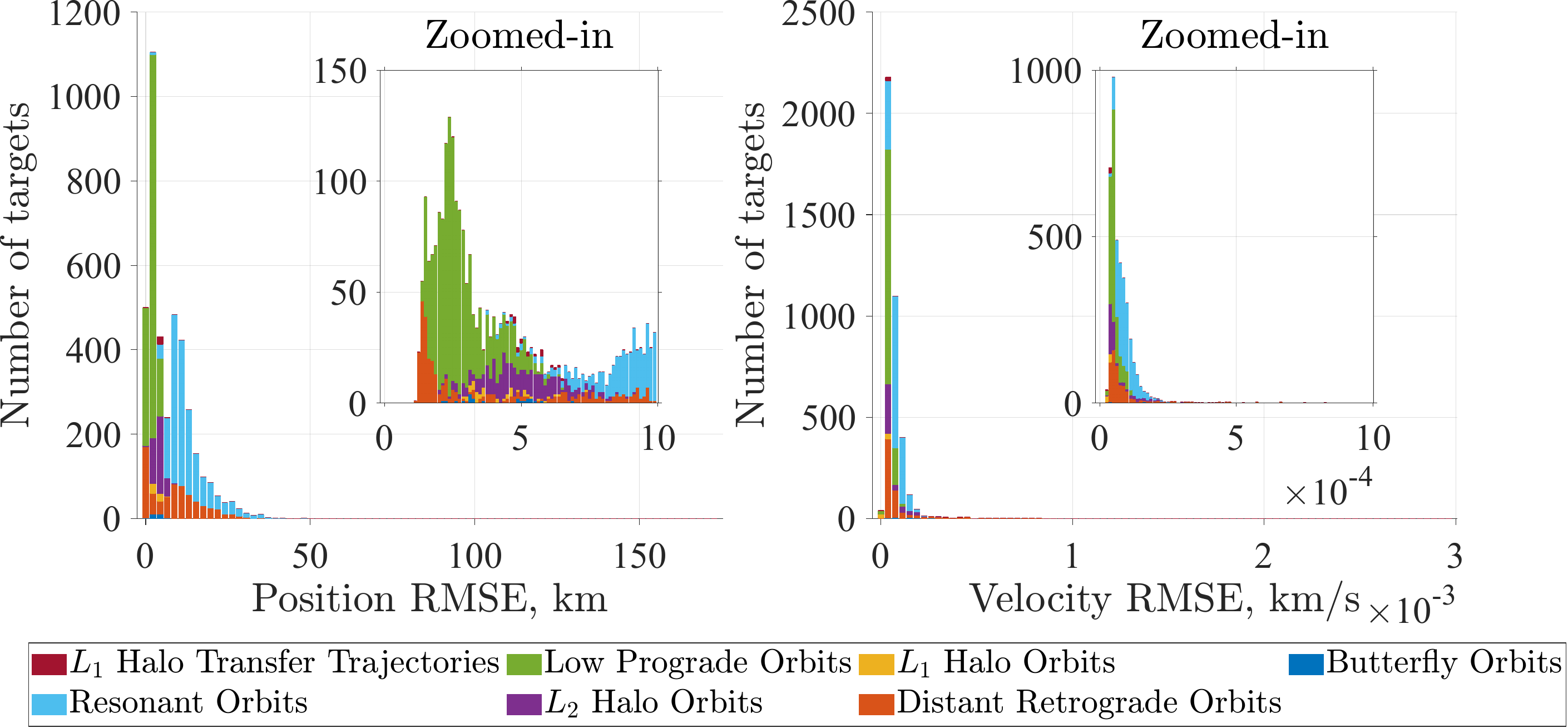}
        \caption{STP-A.}
	  \label{fig:HF_TPA_H}
    \end{subfigure}

    \begin{subfigure}[htbp]{0.55\textwidth}
        \centering
        \includegraphics[width=\linewidth]{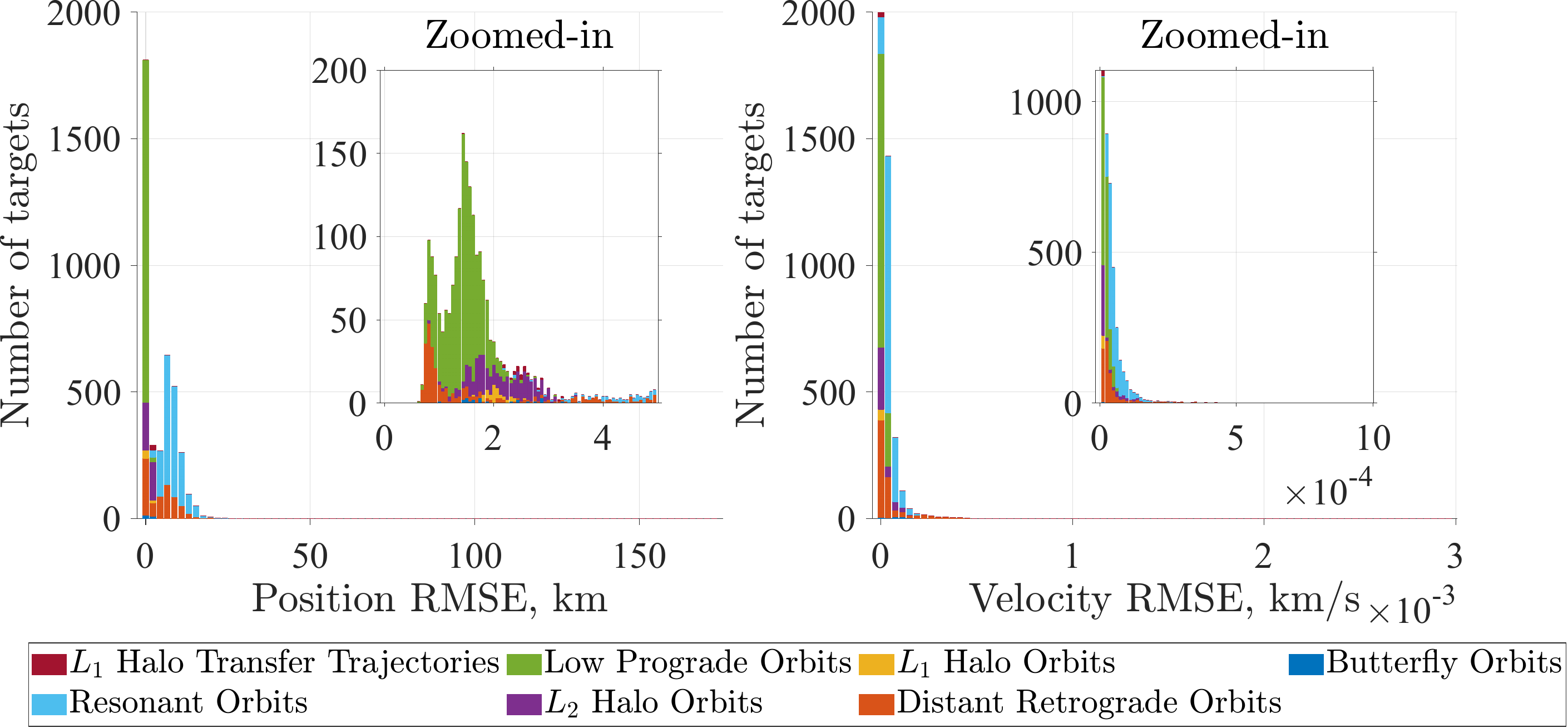}
        \caption{STP-B.}
	    \label{fig:HF_TPB_H}
    \end{subfigure}

    \begin{subfigure}[htbp]{0.55\textwidth}
         \centering
        \includegraphics[width=\linewidth]{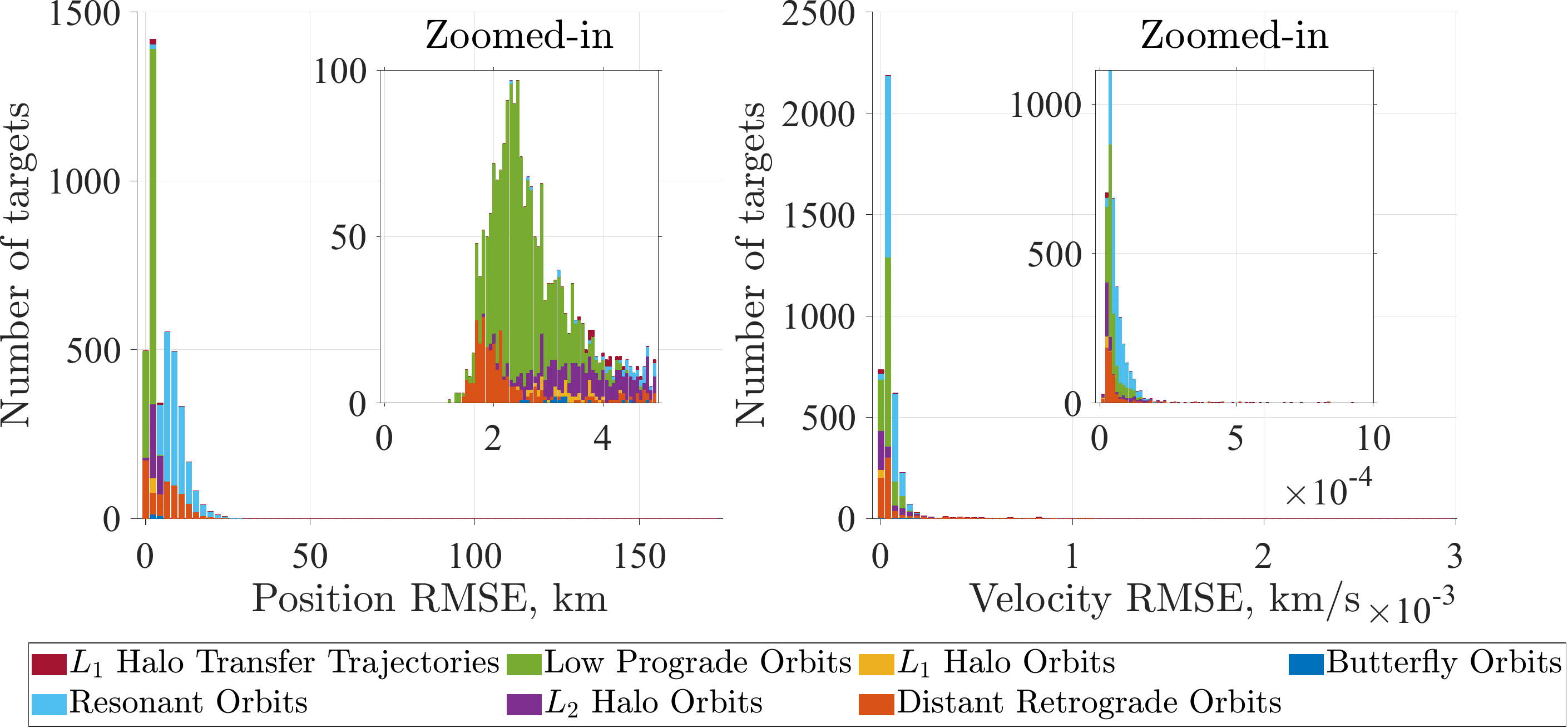}
        \caption{STP-C.}
	   \label{fig:HF_TPC_H}
    \end{subfigure}
    \caption{Histograms showing the distribution of targets separated into similar orbit groupings, for the high-fidelity observers.}
    \label{fig:hifi_histograms}
\end{figure}

Again, to look at the distribution of RMSE values, histograms of RMSE position and velocity are presented with the same grouping as before. They can be seen for the baseline, and all three STPs in Fig.~\ref{fig:hifi_histograms}. Comparing each case with their low-fidelity equivalent, there is an improvement. The grouping of resonant orbits is especially improved. While the results are not as good as the rest of the orbit groupings, they fall much closer to the rest of the targets in the validation target set. Similar to what is seen in the low-fidelity experiments, the targets that each constellation has the worst performance against are the orbits that have longer periods, traveling further than those orbits with shorter periods.

To look into this further, the DRO family is again divided into the orbits with a period greater than \SI{3.75}{TU}, and those with a period less than \SI{3.75}{TU}. The improvement in the long-period DROs is large: the averages are \SI{14.89}{km}, \SI{8.79}{km}, and \SI{10.54}{km} for STP-A, STP-B, and STP-C, respectively, across all three sensor tasking procedures, showing about a three-fold improvement over the low-fidelity sensors. For the short-period DROs, the averages are \SI{3.69}{km}, \SI{2.65}{km}, \SI{1.39}{km}, and \SI{2.47}{km} for the baseline, STP-A, STP-B, and STP-C, respectively.

These results are an improvement over the low-fidelity results, with the high-fidelity sensors providing RMSE values of less than \SI{2}{km} for 374 targets across the DRO and LPWO families using STP-A, 1722 targets across every family but the resonant orbits and transfer trajectories using STP-B, and 288 targets in the DRO, L2SHO, and LPWO families using STP-C. Even for the best-performing sensor tasking procedure, less than half of the validation target set returns a position RMSE of \SI{2}{km}.

While these values are an improvement over the results seen in the low-fidelity observer experiment, they can be further improved. The trends are the same; long-period orbits have less visibility, leading to poorer results. With this in mind, the next steps could be as follows: focus the optimization on a smaller target set so that the results of the targets in the most important location (close to the Moon) have low enough uncertainties. Across these two experiments, STP-B has performed the best, so in future experiments aiming to further improve the constellation's state estimation ability, it is the sensor tasking procedure that should be used. Additionally, more observers can be used, leading to even shorter system cadences (and hypothetically, better results). The benefit of this project is that in formulating it as a framework, all of these changes can be made quickly in response to results that do not meet a user's requirements.

\section{Conclusions}
\label{sec:conc}
This paper proposes the Cislunar SDA Constellation Design and Analysis Framework for the preliminary design of cislunar SDA constellations. Given a set of design parameters and constraints, the framework performs multi-observer placement optimization by conducting orbit determination on a small set of targets with a state estimation filter. Due to the vast size of the solution space and its non-linearity, a metaheuristic solver is used to optimize the multi-observer placement optimization problem, where the objective is to minimize the loss function, an average of the RMSEs of a set of orbits sampled from orbit families that are of interest. Through post-optimization and comparative analysis, we seek to demonstrate the efficacy of an optimized constellation in performing state estimation against a wide set of targets.

The comparative analyses provide insights into how optimized constellations change when various parameters are modified. In comparisons between a single observer and an optimized constellation, regardless of the sensor tasking procedure, the constellation consistently outperforms the single observer due to redundancy and a shorter system cadence. Although the GA is not guaranteed to find a global optimum, in both the low- and high-fidelity cases, the STP-B performs the best, suggesting that the best sensor tasking procedure is the one that enables the most measurements. High-fidelity sensors outperform low-fidelity sensors, and the optimal constellations vary between experiments. The optimized constellations perform best against orbits close to the Moon. To improve the results of orbits that travel long distances, there are several options: designing and optimizing a more focused constellation for orbit families in high-traffic areas like those close to the Moon, increasing the number of observers in the constellation, or enhancing their sensor fidelity.

The optimized constellations that these comparative analyses have found can be studied for deeper analysis of why they perform state estimation better against some target orbits rather than others. An additional analysis could be investigating the observability of the system \cite{li2019observability,salau2014observability,greaves2021relative}. Furthermore, in the state estimation problem, the model-plant dynamics mismatch could also be approached, using the aforementioned higher-level dynamics models, or an ephemeris model for reality, but having the state estimation filter use a lower-fidelity model, testing how the constellations account for this mismatch. Additionally, convergence time could be studied as a way to quantify which targets the constellations perform better or worse against, and the optimization problem could be shifted to focus on minimizing this time across a set of targets.

Future research can be conducted to enhance the fidelity of the models used in the framework. Examples include but are not limited to, modifying the measurement model to include other types of measurements, such as range, or implementing higher fidelity constraints on the measurement model such as visual magnitude, or implementing a higher-fidelity dynamics model such as the Elliptic Restricted Three-Body Problem, the Bi-Circular Restricted Four-Body Problem, or a full ephemeris model. Additionally, other state estimation filters, such as the Unscented Kalman filter, could be tested to see if they improve the results of orbit determination for the targets within cislunar space. This could begin with simply testing these optimized constellations against the same target set, but using another filter instead of the EKF. Sensor tasking optimization and tracking multiple targets simultaneously are important aspects to include as part of the optimization problem to provide better constellation designs. Determining the optimal measurement interval for different orbit families and conducting sensitivity analyses to understand how changes in various parameters such as the number of sensors affect the optimized solution are also valuable research directions.

\appendix
\section*{Appendix A: Initial Conditions of the Optimization Target Set}
This Appendix lists the initial conditions of the 39 targets in the optimization target set.
\begin{table}[!h]
\renewcommand{\arraystretch}{1}
   \caption{Optimization target set initial conditions.}
     \centering
\begin{threeparttable}
\begin{tabular}{llrrrrrrr}
\hline \hline
Target & Orbit family       & Period, TU & $x$, DU            & $y$, DU & $z$, DU & $\Dot{x}$, DU/TU & $\Dot{y}$, DU/TU & $\Dot{z}$, DU/TU \\
\hline
1      & BNO   & 3.72656     & 0.90216            & 0.01025  & 0.13052  & 0.01993           & -0.04143          & 0.12319           \\
2      & BNO   & 2.75985     & 1.01814            & -0.02810 & 0.13269  & -0.04132          & -0.04910          & -0.22452          \\
3      & BNO   & 2.75700     & 1.01044            & -0.03199 & 0.10003  & -0.05800          & -0.01358          & -0.32851          \\
4      & BSO   & 3.72656     & 1.07756            & 0.00976  & -0.12776 & 0.03207           & 0.01846           & 0.15752           \\
5      & BSO   & 2.75985     & 1.02520            & -0.00984 & -0.17047 & -0.01205          & -0.07941          & 0.06354           \\
6      & BSO   & 2.75700     & 1.01735            & 0.02667  & -0.13827 & 0.03825           & -0.05404          & -0.20711          \\
7      & L1NHO  & 2.07099     & 0.91645            & -0.08732 & 0.14063  & -0.12360          & 0.11120           & 0.22361           \\
8      & L1NHO  & 2.15460     & 0.93444            & 0.10611  & 0.09408  & 0.15246           & 0.02387           & -0.29486          \\
9      & L1NHO  & 2.23455     & 0.97875            & 0.07001  & -0.01981 & 0.12094           & -0.37603          & -0.35903          \\
10     & L1SHO  & 2.07099     & 0.89220            & 0.04742  & -0.18274 & 0.06385           & 0.19205           & 0.10694           \\
11     & L1SHO  & 2.15460     & 0.97850            & 0.07201  & 0.00754  & 0.13248           & -0.33956          & 0.39185           \\
12     & L1SHO  & 2.23455     & 0.91328            & -0.10537 & -0.12101 & -0.13910          & 0.09544           & -0.24442          \\
13     & L2NHO  & 1.38944     & 1.00411            & -0.02869 & 0.12826  & -0.04318          & -0.04435          & -0.24523          \\
14     & L2NHO  & 2.17970     & 1.06334            & 0.04467  & 0.18940  & 0.04864           & -0.16732          & 0.09728           \\
15     & L2NHO  & 2.38349     & 1.08295            & 0.00159  & 0.20231  & 0.00162           & -0.20101          & 0.00288           \\
16     & L2SHO  & 1.38944     & 0.98882            & -0.01805 & -0.01252 & -0.07145          & 0.42200           & 0.90381           \\
17     & L2SHO  & 2.17970     & 1.01834            & -0.09432 & -0.08417 & -0.11375          & 0.01927           & 0.32480           \\
18     & L2SHO  & 2.38349     & 1.06990            & -0.06801 & -0.17847 & -0.07035          & -0.16299          & 0.13273           \\
19     & R1:1O         & 6.26635     & -0.12773           & -0.89888 & 0        & -0.10385          & 0.87765           & 0                 \\
20     & R1:1O         & 6.27397     & 0.30808            & -1.37988 & 0        & -1.07683          & 0.26113           & 0                 \\
21     & R1:1O         & 6.27999     & 0.22787            & -1.39504 & 0        & -1.08227          & 0.33668           & 0                 \\
22     & R2:1O         & 5.88685     & -0.53447           & 0.06119  & 0        & -0.07469          & -0.94045          & 0                 \\
23     & R2:1O         & 6.23346     & -0.27953           & -0.31375 & 0        & 0.62241           & -1.26808          & 0                 \\
24     & R2:1O         & 6.27998     & -0.22095           & 1.01044  & 0        & 0.43976           & 0.32928           & 0                 \\
25     & R4:1O         & 6.27997     & 0.03900            & -0.36447 & 0        & 1.15062           & 0.70232           & 0                 \\
26     & R4:1O         & 6.27549     & -0.43950           & -0.44017 & 0        & 0.21767           & -0.11141          & 0                 \\
27     & R4:1O         & 6.27946     & -0.45983           & 0.44568  & 0        & -0.31789          & 0.22943           & 0                 \\
28     & DRO   & 0.15382     & 0.98707            & 0.01970  & 0        & 0.80343           & 0.03198           & 0                 \\
29     & DRO   & 5.82922     & 0.94147            & 0.77868  & 0        & 0.58313           & -0.17725          & 0                 \\
30     & DRO   & 6.27993     & 0.20056            & 1.38560  & 0        & 1.06003           & 0.36870           & 0                 \\
31     & LPEO & 1.34333     & 1.06653            & 0.00080  & 0        & -0.00294          & 0.30862           & 0                 \\
32     & LPEO & 1.67485     & 1.04983            & -0.06794 & 0        & 0.22149           & 0.05401           & 0                 \\
33     & LPEO & 2.56988     & 1.12607            & 0.02203  & 0        & -0.01593          & 0.08709           & 0                 \\
34     & LPWO & 0.19790     & 0.98893            & -0.02245 & 0        & 0.71085           & 0.03489           & 0                 \\
35     & LPWO & 1.20460     & 0.92649            & 0.04880  & 0        & -0.19084          & -0.20955          & 0                 \\
36     & LPWO & 2.14657     & 0.87494            & -0.04367 & 0        & 0.08386           & -0.03495          & 0                 \\
37     & L1TT           & N/A\tnote{\textdagger}         &  0 & -0.28642 & 0.03740  & 1.93948           & -0.26854          & -0.32641          \\
38     & L1TT           & N/A\tnote{\textdagger}        & 0 & -0.57053 & 0.04766  & 0.83103           & -0.06703          & -0.32023          \\
39     & L1TT           & N/A\tnote{\textdagger}         & 0 & -0.36456 & 0.03514  & 1.53129           & -0.19889          & -0.40470        \\
\hline \hline
\end{tabular}
\begin{tablenotes}
\item[\textdagger] As these are transfer trajectories, they do not have a period.
\end{tablenotes}
\end{threeparttable}
      \label{tab:OTS_ICs}
\end{table}

\section*{Appendix B: Initial Conditions of Optimized Constellation and their Best- and Worst-performing Targets with Low-fidelity Sensors}
This Appendix lists the initial conditions of the three tasking procedure's optimized constellations for the low-fidelity sensor constellations and their best- and worst-performing targets.

\begin{table}[!h]
\renewcommand{\arraystretch}{1}
\centering
\caption{Optimized constellation initial conditions for low-fidelity sensors.}
\label{tab:optimal_ICs}
\begin{tabular}{llrrrrrrr}
\hline \hline
STP & Orbit family                    & Period, TU & $x$, DU & $y$, DU & $z$, DU & $\Dot{x}$, DU/TU & $\Dot{y}$, DU/TU & $\Dot{z}$, DU/TU \\
\hline
Baseline & NRHO & 1.48 & 1.0192 & 0.0084 & -0.1785 & 0.0109 & -0.0964 & -0.0445 \\

\midrule
STP-A       & L2NHO & 1.42034     & 1.01059  & -0.02336 & 0.15508  & -0.03225          & -0.07038          & -0.15740          \\
       & L2SHO & 2.18236     & 1.05129  & -0.07285 & -0.16417 & -0.08092          & -0.12861          & 0.17317           \\
        & R2:1O        & 6.26478     & 0.12436  & -0.71102       & 0        & 0                 & -0.58951          & 0                 \\
        & DRO        & 0.31476     & 0.97797  & 0.03082  & 0        & 0.61175           & 0.19864           & 0                \\
        \midrule
STP-B        & L2NHO       & 1.42034     & 1.01059  & -0.02336 & 0.15508  & -0.03225          & -0.07038          & -0.15740          \\
        & L2SHO & 2.35117     & 1.00808  & 0.09471  & -0.03304 & 0.10807           & 0.15115           & -0.37438          \\
        & R4:1O              & 6.27538     & 0.18900  & 0        & 0        & 0                 & 2.50444           & 0                 \\
        & LPEO      & 1.92854     & 0.88714  & 0.04674  & 0        & -0.12273          & -0.04305          & 0                \\
        \midrule
STP-C & BSO & 3.73313                     & 0.90452                      & 0                            & -0.14384                     & 0                                     & -0.05031                              & 0                                     \\
 & L2NHO  & 1.42034                    & 1.01059                      & -0.02336                     & 0.15508                      & -0.03225                              & -0.07038                              & -0.15740                              \\
 & L2SHO  & 2.18236                     & 1.05129                      & -0.07285                     & -0.16417                     & -0.08092                              & -0.12861                              & 0.17317                               \\
 & R4:1O & 6.27542                     & 0.35512                      & -0.43599                     & 0                            & 0.43122                               & 0.08299                               & 0                                    \\
\hline \hline
\end{tabular}
\end{table}

\begin{table}[!h]
\caption{Best and worst target initial conditions for low-fidelity sensors.}
\label{tab:targ_ics_lofi}
\begin{tabular}{lllrrrrrrr}
\hline
\hline
STP      &  & Orbit family & Period, TU & $x$, DU  & $y$, DU  & $z$, DU & $\Dot{x}$, DU/TU & $\Dot{y}$, DU/TU & $\Dot{z}$, DU/TU \\
\hline
Baseline & Best                 & LPWO         & 0.22451    & 0.97333  & 0.01959  & 0       & -0.54626         & -0.40633         & 0                \\
         & Worst                & R1:1O        & 6.27578    & 0.15585  & -1.33939 & 0       & -0.97105         & 0.43428          & 0                \\
         \midrule
STP-A    & Best                 & LPWO         & 0.25211    & 0.98329  & 0.02589  & 0       & -0.64209         & -0.11315         & 0                \\
         & Worst                & R1:1O        & 6.27387    & -0.16035 & 1.00472  & 0       & 0.30638          & 0.84852          & 0                \\
         \midrule
STP-B    & Best                 & LPWO         & 0.30606    & 0.97149  & -0.02485 & 0       & 0.50671          & -0.33584         & 0                \\
         & Worst                & R1:1O        & 6.27709    & 0.36954  & 1.40995  & 0       & 1.14308          & 0.18072          & 0                \\
         \midrule
STP-C    & Best                 & LPWO         & 0.25640    & 0.96200  & 0.00591  & 0       & -0.14431         & -0.63486         & 0                \\
         & Worst                & R1:1O        & 6.27833    & -0.18097 & -1.05186 & 0       & -0.39482         & 0.84371          & 0               \\
\hline
\hline
\end{tabular}
\end{table}

\section*{Appendix C: Initial Conditions of Optimized Constellation and their Best- and Worst-performing Targets with High-fidelity Sensors}
This Appendix lists the initial conditions of the three tasking procedure's optimized constellations for the high-fidelity sensor constellations and their best- and worst-performing targets.

\begin{table}[!h]
\renewcommand{\arraystretch}{1}
\centering
\caption{Optimized constellation initial conditions for high-fidelity sensors.}
\label{tab:opt_ICs_hifi}
\begin{tabular}{llrrrrrrr}
\hline
\hline
   STP & Orbit family & Period, TU & $x$, DU & $y$, DU & $z$, DU & $\Dot{x}$, DU/TU & $\Dot{y}$, DU/TU & $\Dot{z}$, DU/TU \\
      \hline
      Baseline & NRHO & 1.48 & 1.0192 & 0.0084 & -0.1785 & 0.0109 & -0.0964 & -0.0445 \\

\midrule
STP-A & L2SHO   & 1.45344                        & 0.99772                     & -0.03692                    & -0.07986                    & -0.06613                             & 0.02110                              & 0.40552                              \\
      & R2:1O   & 6.27070                        & -0.15705                    & -0.98967                    & 0                           & -0.36541                             & 0.24333                              & 0                                    \\
      & R4:1O   & 6.27586                        & 0.40259                     & -0.44657                    & 0                           & 0.31213                              & -0.02049                             & 0                                    \\
      & LPWO    & 1.65200                        & 0.87264                     & -0.02046                    & 0                           & 0.03737                              & -0.11888                             & 0                                    \\
      \midrule
STP-B & L2SHO   & 1.45344                        & 0.99772                     & -0.03692                    & -0.07986                    & -0.06613                             & 0.02110                              & 0.40552                              \\
      & R4:1O   & 6.27650                        & -0.41703                    & 0.08243                     & 0                           & 0.41557                              & -1.06990                             & 0                                    \\
      & LPEO    & 1.67328                        & 1.09913                     & 0                           & 0                           & -1.90E-14                            & 0.16932                              & 0                                    \\
      & LPWO    & 2.06136                        & 0.88147                     & 0.04539                     & 0                           & -0.10630                             & -0.03484                             & 0                                    \\
      \midrule
STP-C & L2SHO   & 2.37810                        & 1.06040                     & -0.08494                    & -0.16070                    & -0.08875                             & -0.13340                             & 0.17647                              \\
      & R2:1O   & 6.27070                        & -0.15705                    & -0.98967                    & 0                           & -0.36541                             & 0.24333                              & 0                                    \\
      & R4:1O   & 6.27650                        & -0.41703                    & 0.08243                     & 0                           & 0.41557                              & -1.06990                             & 0                                    \\
      & LPWO    & 1.51512                        & 1.01427                     & 0                           & 0                           & -1.52E-14                            & 0.82818                              & 0  \\
      \hline
      \hline
\end{tabular}
\end{table}

\begin{table}[!h]
\caption{Best and worst target initial conditions for high-fidelity sensors.}
\label{tab:targ_ics_hifi}
\begin{tabular}{lllrrrrrrr}
\hline
\hline
STP      &  & Orbit family & Period, TU & $x$, DU  & $y$, DU  & $z$, DU & $\Dot{x}$, DU/TU & $\Dot{y}$, DU/TU & $\Dot{z}$, DU/TU \\
\hline
Baseline & Best                 & DRO          & 0.09373    & 0.98790  & -0.01407 & 0       & -0.94290         & -0.00297         & 0                \\
         & Worst                & R2:1O        & 6.19203    & -0.24127 & -0.50917 & 0       & 0.49928          & -0.74750         & 0                \\
         \midrule
STP-A    & Best                 & DRO          & 0.09798    & 0.97585  & 0.00813  & 0       & 0.52119          & 0.77041          & 0                \\
         & Worst                & R1:1O        & 6.27628    & -0.19185 & 0.73649  & 0       & -0.25470         & 0.98836          & 0                \\
         \midrule
STP-B    & Best                 & LPWO         & 0.22179    & 0.97216  & -0.01841 & 0       & 0.51976          & -0.44460         & 0                \\
         & Worst                & DRO          & 6.27971    & 0.67698  & -1.42470 & 0       & -1.26185         & -0.18278         & 0                \\
         \midrule
STP-C    & Best                 & LPWO         & 1.52186    & 0.98176  & 0.04731  & 0       & -0.45863         & 0.25174          & 0                \\
         & Worst                & R2:1O        & 6.27671    & 0.03235  & 0.44718  & 0       & -0.94585         & -0.99303         & 0               \\
\hline
\hline
\end{tabular}
\end{table}

\section*{Acknowledgment}
The authors thank the anonymous reviewers for their insightful suggestions, which have greatly contributed to improving the overall quality of the paper.

\clearpage
\newpage
\bibliography{references}

\end{document}